\documentclass[11pt, dvipsnames]{amsart}
\usepackage[utf8]{inputenc}
\usepackage{graphicx}
\usepackage{cleveref}
\usepackage{a4wide}
\usepackage{color}
\usepackage{xcolor}
\usepackage{amsmath,amssymb}
\usepackage{float}
\usepackage{soul}
\usepackage{booktabs} 
\usepackage[foot]{amsaddr}
\usepackage{amsthm}

\newtheorem{lemma}{Lemma}[section]
\newtheorem{example}{Example}[section]
\newtheorem{assumption}{Assumption}[section]
\newtheorem{prop}{Proposition}[section]
\newtheorem{remark}{Remark}[section]

\numberwithin{equation}{section} 

\title[Multi-dimensional solute transport modelling errors]{The modelling error in multi-dimensional time-dependent solute transport models}
\author{Rami Masri $^1$, Marius Zeinhofer$^1$,  Miroslav Kuchta$^1$, and Marie E. Rognes$^1$}
\email{rami@simula.no, mariusz@simula.no, miroslav@simula.no, meg@simula.no}
\address{$^1$Department of Numerical Analysis and Scientific Computing, Simula Research Laboratory}
\address{$^1$This project has received support and funding from the European Research Council (ERC) under the European Union's Horizon 2020 research and innovation programme under grant agreement 714892 (Waterscales) and from the Research Council of Norway (RCN) via FRIPRO grant agreement 324239 (EMIx).}
\date{ \today }

\newcommand{\avg}[1]{\langle #1 \rangle}
\newcommand{\bm}[1]{\boldsymbol{#1}}

\newcommand{\ds}{\,\mathrm{d}s} 
\newcommand{\dr}{\,\mathrm{d}r} 
\newcommand{\dtheta}{\,\mathrm{d}\theta} 
\newcommand{\dz}{\,\mathrm{d}z} 
\newcommand{\dt}{\,\mathrm{d}t} 

\newcommand{\R}{\mathbb{R}} 
\newcommand{\foralls}{\, \forall} 

\newcommand{\cc}{\hat{c}} 
\newcommand{\vvv}{g} 
\newcommand{\vtwo}{\vvv_s} 
\newcommand{\avgus}{\hat{u}} 

\newcommand{\diam}{\mathrm{diam}}

\begin{document}

\maketitle

\begin{abstract}
  Starting from full-dimensional models of solute transport, we derive and analyze multi-dimensional models of time-dependent convection, diffusion, and exchange in and around pulsating vascular and perivascular networks. These models are widely applicable for modelling transport in vascularized tissue, brain perivascular spaces, vascular plants and similar environments. We show the existence and uniqueness of solutions to both the full- and the multi-dimensional equations under suitable assumptions on the domain velocity. Moreover, we quantify the associated modelling errors by establishing a-priori estimates in evolving Bochner spaces. In particular, we show that the modelling error decreases with the characteristic vessel diameter and thus vanishes for infinitely slender vessels. Numerical tests in idealized geometries corroborate and extend upon our theoretical findings.
  \vspace{1em}


\end{abstract}

\section{Introduction}

We consider transport of solutes by diffusion, convection, and
exchange in a coupled system consisting of networks of slender vessels
and their surroundings. This setting is ubiquitous in the human
body~\cite{boron2012medical} as exemplified by the transport and
exchange of nutrients such as oxygen or glucose, or medical drugs in
the vasculature and surrounding tissue, e.g.~in skeletal muscle, the
liver~\cite{rohan2018modeling}, or the
placenta~\cite{wheeler2021bioengineering}; or conversely, the
transport of metabolic by-products from tissue into and through
lymphatic vessels~\cite{possenti2019numerical}. Similar structures and
processes are also fundamental in biology, think of e.g.~the roots of
vascularized plants~\cite{koch2022nonlinear}, and in geoscience
e.g.~in connection with flow and transport in reservoir
wells~\cite{gjerde2020singularity}, in the context of C02
sequestration~\cite{nordbotten2009model}, or groundwater
contamination~\cite{malenica2018groundwater}.

Of particular interest, both from a physiological and mathematical
point-of-view, is the transport of solutes in, around and out of the
human \emph{brain}. Despite decades -- even centuries -- of research,
solute transport and clearance within the human brain remain poorly
understood~\cite{zhao2022physiology, hladky2022glymphatic,
  kelley2022glymphatic}. In contrast to the rest of the body, the
brain vasculature is equipped with a blood-brain-barrier, which
carefully regulates the exchange of substances between the blood and
the surrounding tissue, while the brain parenchyma itself lacks
typical lymph vessels. Better understanding of these physiological
processes is vital for targeting brain drug
delivery~\cite{nance2022drug, lohela2022glymphatic} or for unraveling
the role of metabolic waste clearance in neurodegenerative
disease~\cite{tarasoff2015clearance,
  kelley2022glymphatic}. Concurrently, in tissue engineering, efforts
are currently underway to develop human brain cortical organoids, but
crucially rely on vascularization via e.g.~microfluidic devices for
improved oxygen and nutrient transport as well as cellular
signalling~\cite{lamontagne2022recent}.

The human brain is composed of soft tissue, is lined and penetrated by
networks of blood vessels, and is surrounded by the narrow
subarachnoid space filled with cerebrospinal fluid (CSF). The cerebral
arteries pulsate in sync with the cardiac cycle and undergo other
forms of vasomotion with variations in radii of
$\sim$1--10\%~\cite{mestre2018flow}, while the entire brain parenchyma
deforms by around 1\% as the result of a complex interplay between the
cardiac and respiratory cycles as well as
autoregulation~\cite{sloots2020cardiac,
  causemann2022human}. Perivascular (or paravascular) spaces (PVSs)
are spaces surrounding the vasculature on the brain surface or within
the brain parenchyma. On the brain surface, these spaces are clearly
visible~\cite{vinje2021brain}, and PVSs persist as the blood vessels
branch and penetrate into the brain parenchyma -- then known as
Virchow-Robin spaces. The extent to which perivascular spaces exist
along the length of the vasculature within the brain, even to the
capillary level, is debated
however~\cite{hannocks2018molecular}. Within the parenchyma,
perivascular spaces are often represented as generalized (elliptic)
annular cylinders, filled with cerebrospinal or interstitial fluid and
bounded by a nearly tight layer of astrocyte endfeet, see
e.g.~\cite{brinker2014new, wardlaw2020perivascular,
  daversin2020mechanisms} and references therein.

Solutes move by diffusion within the brain
tissue~\cite{nicholson2001diffusion}, and by diffusion and convection
within the vasculature~\cite{boron2012medical}. However, to what
extent also convection in \emph{perivascular}, \emph{intracellular} or
\emph{extracellular spaces} play a role in brain solute transport and
clearance stand as important open questions. Convective velocity
magnitudes are expected to differ by many orders of magnitude between
and within the respective compartments: blood may flow at the order of
1 m/s in major cerebral arteries~\cite{boron2012medical}, CSF flows in
surface perivascular spaces at up to 60 $\mu$m/s with P\'eclet numbers
of up to 1000~\cite{mestre2018flow}, while flow of interstitial fluid
within the tissue is unlikely to exceed 10 $\mu$m/min on
average~\cite{vinje2023human, abbott2004evidence}. Depending on their
ability to cross the blood-brain barrier, solutes may also exchange
between the vascular and perivascular spaces, as well as into the
surrounding tissue or subarachnoid space. To mathematically and
computationally study such transport at the scale of larger vascular
networks, our target here is to derive and analyze time-dependent
convection-diffusion models with a geometrically-explicit but
dimensionally-reduced representation of the (peri)vascular spaces
coupled with the full-dimensional surroundings.

As a starting point (more precise details are presented later), consider second-order elliptic equations
describing diffusion of the solute concentrations $c_v : \Omega_v
\rightarrow \R$, and $c_s : \Omega_s \rightarrow \R$:
\begin{subequations}
  \begin{align}
    - \nabla \cdot D \nabla c_s &= f \quad \text{ in } \Omega_s, \\
    - \nabla \cdot D \nabla c_v &= f \quad \text{ in } \Omega_v,
  \end{align}
  \label{eq:intro:3d-3d}%
\end{subequations}%
where $D$ is a given effective diffusion coefficient and $f$ given
sources. Assuming that the compartments are separated by a
semi-permeable membrane $\Gamma$ gives the interface condition
\begin{equation}
  D \nabla c \cdot \bm{n} = \xi [[c]] \quad \text{ on } \Gamma,
  \label{eq:intro:gamma}
\end{equation}
where $\bm{n}$ is the interface normal, $[[\cdot]]$ denotes the jump across
the interface(s), and $\xi$ is a membrane permeability parameter. Now
assuming that $\Omega_v$ can be well-represented by its centerline
$\Lambda$ with coordinate $s$ (to be made more precise later), the
coupled 3D-3D problem
of~\eqref{eq:intro:3d-3d}--\eqref{eq:intro:gamma} may be reduced to a
coupled 3D-1D problem of the form: find the solute concentrations
$\bar{c} : \Lambda \rightarrow \R$, and $c : \Omega \rightarrow \R$
\begin{subequations}
  \begin{align}
    - \nabla \cdot D \nabla c + \mathcal{C} = f \quad \text{ in } \Omega, \\
    - \partial_{s} \bar{D} \partial_s \bar{c} + \bar{\mathcal{C}} = \bar{f} \quad \text{ on } \Lambda, 
  \end{align}
  \label{eq:intro:3d-1d}%
\end{subequations}%
where $\mathcal{C}$, $\bar{\mathcal{C}}$ denote coupling terms
depending on the concentrations $c, \bar{c}$, and choice of
coupling. 
Note that flow in a porous medium (Darcy flow) can be
described with the same equation structure, with $c$ instead
representing the pore pressure and $D$ the hydraulic
conductance. Modelling, discretization, and applications of 3D-1D
problems such as~\eqref{eq:intro:3d-1d} has been the subject of active
research, especially over the last two decades, with key contributions
from e.g.~\cite{fleischman1986interaction, d2008coupling, d2012finite,
  notaro2016mixed, koppl2018mathematical, gjerde2019splitting,
  kuchta2019preconditioning, koch2020modeling, kuchta2021analysis} and references therein
to mention but a few. Notably, Laurino and
Zunino~\cite{laurino2019derivation} rigorously analyze the modelling
error associated with
replacing~\eqref{eq:intro:3d-3d}--\eqref{eq:intro:gamma}
by~\eqref{eq:intro:3d-1d}, and demonstrate that the modelling error
indeed vanishes for infinitely thin vessels.

Here, we consider a parabolic extension of the classical elliptic
3D-1D equations~\eqref{eq:intro:3d-1d} accounting also for (i)
time-evolving distributions, (ii) convective transport, (iii)
moving interfaces, and (iv) both cylindrical and non-convex
(annular) vessel networks representing e.g.~vascular and perivascular
spaces, respectively. We also derive and study a 3D-1D-1D model
representing solute transport in coupled tissue, perivascular and
vascular spaces. Previously, Possenti, Zunino and
coauthors~\cite{possenti2021mesoscale} and K\"oppl, Vidotto and
Wohlmuth~\cite{koppl20203d} have studied applications of 3D-1D models
for (oxygen) transport including convection but at steady
state. Furthermore, Formaggia et al~\cite{formaggia2001coupling}
consider coupled Navier--Stokes equations for flow problems in
compliant vessels but with a different type of mixed-dimensional
coupling. More specifically, we are interested in solute
concentrations $c_v(t) : \Omega_v(t) \rightarrow \R$, and $c_s(t) :
\Omega_s(t) \rightarrow \R$ satisfying the time-dependent diffusion
equations for a.e.~$t > 0$:
\begin{subequations}
  \begin{align}
    \partial_t c_{s} + \bm{u} \cdot \nabla c_s - \nabla \cdot D \nabla c_s &= f \quad \text{ in } \Omega_s(t), \\
    \partial_t c_{v} + \bm{u} \cdot \nabla c_v - \nabla \cdot D \nabla c_v &= f \quad \text{ in } \Omega_v(t), 
  \end{align}
  \label{eq:intro:conv-diff}%
\end{subequations}%
where now additionally $\bm{u}$ represents a convective velocity field
and the interface $\Gamma$ between $\Omega_s$ and $\Omega_t$ is
allowed to move and deform in time.  

Our main findings are as follows.
\begin{itemize}
\item
  We introduce a system of time-dependent convection-diffusion
  equations in and around embedded networks of moving vessels. Under
  suitable assumptions on the domain velocity, we prove
  well-posedness i.e.~that suitably regular weak solutions  to these
  equations exist and are unique (\Cref{sec:transport}).
  \item We derive reduced 1D equations, and we formally derive weak
    formulations of 3D-1D and 3D-1D-1D coupled models of
    time-dependent solute transport governed by convection, diffusion,
    and exchange in deforming vascular and/or perivascular networks,
    and the surrounding domain (\Cref{sec:3D1D},
    \Cref{sec:3D1D1D}). We prove well-posedness of the coupled 3D-1D
    formulation and show a regularity estimate for the 3D
    solution. These formulations are widely applicable for modelling
    transport in vascularized tissue in general and the brain in
    particular, vascular plant environments etc.
\item
  We rigorously estimate the modelling error in evolving
  Bochner spaces associated with replacing
  the time-dependent 3D-3D convection-diffusion problem by the 3D-1D
  problem via a duality argument. We show that a relevant dual problem
  is well-posed, and that the modelling error decreases with the
  characteristic vessel diameter $\epsilon$, and thus vanishes as
  $\epsilon \rightarrow 0$ (\Cref{sec:modelling_error}).
\item
  The presence of deforming networks with annular cross-sections poses
  key technical challenges relating to classical numerical analysis
  tools, such as e.g.~Poincar\'e and trace inequalities, and extension
  operators over moving, non-convex domains, which we address
  separately (\Cref{sec:inequalities}).
\end{itemize}
These points are prefaced by introducing notation and preliminary
results in \Cref{sec:notation}, while concluding remarks and outlook
relating to e.g.~the discretization errors form \Cref{sec:conclusion}.

\section{Notation and preliminaries}
\label{sec:notation}

\subsection{Function spaces, inner products and norms}

Given an open domain $O \subset \R^d$, $d \in \{1,2,3\}$ and measurable real valued functions $f,g$ , we let $(f,g)_{O}$ denote the usual $L^2$ inner product. If $O$ is the whole domain $\Omega$, the we write $(f,g) = (f,g)_{\Omega}$. The Hilbert space gernerated by this inner product is denoted by $L^2(O)$ with the usual induced norm $\|\cdot\|_{L^2(O)}$. We also use standard notation for the Sobolev spaces $W^{m,p}(O)$ and $H^{m}(O) = W^{m,2}(O)$ for $ m \in \mathbb{N}$ and $ 1 \leq p \leq \infty$. For a given weight $w \in L^{\infty}(O)$ and $w > 0$ a.e. in $O$, we define the weighted $L^2$ inner product $(f,g)_{O,w} = (f,wg)_O$ and the respective weighted $L^2$ space :
\begin{equation} \label{eq:weighted_inner_product}
\|f\|_{L^2_w(O)} =\|w^{1/2}f\|_{L^2(O)}, \,\,  L^2_{w}(O) = \{f:O \rightarrow \R \,\,  \vert \,\, \|f\|_{L_w^2(O)} < \infty\}. 
\end{equation}
The weighted Sobolev space $H^1_w(O)$ is then defined as  
\begin{equation}
    H^1_{w}(O) = \{f \in L^2_{w}(O) \,\, \vert \,\,  \|\nabla f\|_{L^2_w(O)} <  \infty\}, 
\end{equation}
and the weighted inner product and norm are
\begin{equation}
(f,g)_{H^1_{w}(O)} = (f,g)_{L^2_w(O)} + (\nabla f, \nabla g)_{L^2_{w}(O)}, \,\,\, \|f\|_{H^1_{w}(O)}^2 = \|f\|_{L^2_w(O)}^2 + \|\nabla f \|_{L^2_w(O)}^2. 
\end{equation}
We omit the subscript/weight $w$ when $w = 1$.

Given a Hilbert space $X$, we denote the dual space of $X$ by $X'$. The duality pairing between $X$ and $X'$ is denoted by
\begin{equation*}
  \langle v' , v \rangle_{X' \times X}.
\end{equation*}
For brevity in notation, we let
\begin{equation*}
  H^{-1}_w(O) = H_w^1(O)', \qquad \langle v', v \rangle_{H^{-1}_w(O)} = \langle v', v \rangle_{H^{-1}_w(O) \times H^{1}_w(O)}  .
\end{equation*}
We also recall the definition of standard Bochner type spaces. For $t,T>0$,  $f: (t, T) \rightarrow X$, we say that $f \in L^2(t,T;X)$ if 
\begin{equation}
\|f\|^2_{L^2(t,T;X)} = \int_{t}^T \|f\|_{X}^2 < \infty.
\end{equation}
If $f$ is weakly differentiable in time and $\partial_t f \in L^2(t,T;X)$, then we say $f \in H^1(t,T;X)$ with the norm:
\begin{equation}
\|f\|_{H^1(t,T;X)}^2 = \|f\|_{L^2(t,T;X)}^2 + \|\partial_t f\|_{L^2(t,T;X)}^2.
\end{equation}
Given two Hilbert spaces $V$ and $H$ with $V \subset H$, we define 
\begin{equation} \label{eq:def_bochner_space_W}
    \mathcal{W}(V,H) = \{ v: (0,T) \rightarrow V; v \in L^2(0,T;V), \,\, \partial_t v \in L^2(0,T;H) \}. 
\end{equation}
Finally, we will use the space $C^0(0,T;V)$ of continuous $V$-valued functions and the space $ \mathcal{D}(0,T;V)$ of infinitely differentiable $V$-valued functions.

\subsection{The geometrical setting}
\label{sec:geometry_specific}

\begin{figure}
  \begin{center}
   \includegraphics[width=0.3\textwidth]{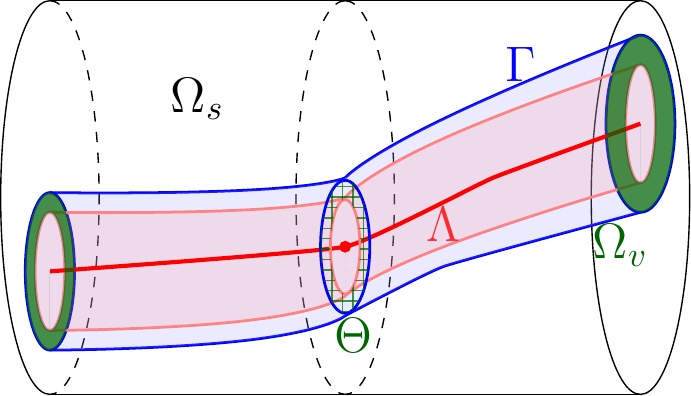}
    \hspace{20pt}
    \includegraphics[width=0.3\textwidth]{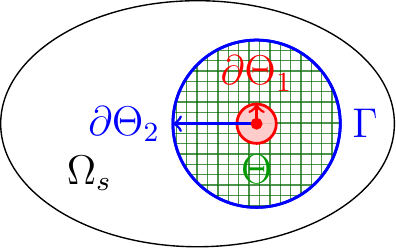}
  \end{center}
 \caption{Geometrical setting. Left: An (annular) cylinder $\Omega_v$
   (shown in green) representing a vessel parametrized in terms of
   the centerline curve $\Lambda$ (shown in red). The vessel is
   surrounded by the domain $\Omega_s$ while $\Gamma$ forms the
   lateral boundary between the domains (shown in blue). Right:
   Lateral cross-section of the domain $\Omega$ showing $\Omega_v$ in
   green, the inner boundary $\partial \Theta_1$ in red, and the outer
   boundary $\partial \Theta_2$ or $\Gamma$ in blue. Note that in the
   annular cylinder case, the inner-most cylinder, i.e.~the extrusion
   along $\Lambda$ of the domain bounded by $\partial \Theta_1$, is
   not part of $\Omega_s$ or $\Omega_v$. }
    \label{fig:geometry}
\end{figure}

We consider a generalized annular domain $\Omega_v$ (\Cref{fig:geometry}), described in
cylindrical coordinates and moving in time:
\begin{multline*}
  \Omega_v(t) = \{ \bm{\lambda}(s) + r\cos(\theta) \bm{N}(s) + r\sin(\theta) \bm{B}(s), \\
  0 < s <  L, \, 0 \leq \theta \leq 2 \pi, \,
  R_1(s,t,\theta) <  r <  R_2(s,t,\theta ) \} \subset \R^3
\end{multline*}
of length $L > 0$, inner radius $R_1 \geqslant 0$, and outer radius $R_2 > 0$. We refer to the $s$-direction as the axial direction. For $R_1 =0$, we consider a cylindrical domain where in the above definition, we let $ 0 \leq r < R_2(s,t,\theta)$. In general, $\Omega_v$ represents a vessel segment such as a perivascular space ($R_1 > 0$), or blood vessel segment, plant root or borehole ($R_1 = 0$). We assume that $\bm{\lambda}(s)  = [\lambda^1(s), \lambda^2(s), \lambda^3(s)]$ above is a parametrized $C^2$-regular curve with non-moving centerline $\Lambda$ defined as $\Lambda = \{\bm{\lambda}(s)\}$ for $s \in (0,L)$, and that $\|\bm{\lambda}'(s)\| = 1$, thus implying that $s$ is the arc length.  The vectors $\bm{N}$ and $\bm{B}$ are from the Frenet-Serret frame of $\Lambda$. Throughout the paper, $\Theta(s,t)$ denotes the cross-section of $\Omega_v(t)$ at $s \in \Lambda$. 

This domain $\Omega_v(t)$ is embedded into a fixed domain $\Omega \subset \mathbb{R}^3$ with (outer) surroundings $\Omega_s = \Omega \backslash B_{R_2}$ where $B_{R_2}$ is the outer cylinder given by: 
\begin{equation*}
B_{R_2}(t) =  \{ \bm{\lambda}(s) + r\cos(\theta) \bm{N}(s) + r\sin(\theta) \bm{B}(s), \\
0 < s <  L, \, 0 \leq  \theta \leq 2 \pi, \,
0 \leq  r <  R_2(s,t,\theta ) \}. 
\end{equation*}
We emphasize that, by construction, the surrounding domain $\Omega_s$ does not include the vessel $\Omega_v$ itself nor the inner-most generalized cylinder in the case $R_1 > 0$. We assume that for all $t \in [0,T]$, $\Omega_v(t)$ is completely embedded in $\Omega$; that is, 
\begin{equation*}
  \mathrm{dist}(\partial \Omega_v(t), \partial \Omega) > 0, \quad \forall t \in [0,T].
\end{equation*}

We denote by $\Gamma$ the lateral boundary of $\Omega_v$ intersecting the boundary of $\Omega_s(t)$, $\Gamma = \partial \Omega_v \cap \partial \Omega_s$, and by
$\Gamma_0$ and $\Gamma_L$ the vertical boundary of $\Gamma$ at $s=0$ and at $s=L$
respectively. The unit normal to $\partial \Omega_v$ is denoted by $\bm{n}_v$, and on $\Gamma$,  $\bm{n}_s = - \bm{n}_v$.
For each cross-section $\Theta(s, t)$, we label its area $A(s,t) =
|\Theta(s,t)|$ for $s \in \Lambda$ and $t \in [0,T]$. 
 We also label the boundary of the lateral cross-section of $\Theta$ by $\partial \Theta$, the boundary of the outer circle (at $r = R_2$)  by $ \partial \Theta_2 $, and (if $R_1>0$) the boundary of the inner circle (at $r= R_1$) by $\partial \Theta_1$.   
We denote by $P(s, t) = |\partial \Theta_2(s,t)|$ the perimeter of the outer circle (representing the interface between the vessel and its outer surroundings).

  Note that we consider vessels $\Omega_v$ both of cylinder-type or
  annular cylinder-type and their (outer) surroundings. The former
  case is well-suited to represent e.g.~transport in the vasculature,
  roots or geothermal wells. The latter case targets e.g.~perivascular
  transport in the brain, intracranial space, or spinal
  compartments. In the latter (annular cylinder) case, we only include
  the outer surroundings in the 3D-3D and reduced 3D-1D model
  formulations in the subsequent
  Sections~\ref{sec:transport}--\ref{sec:3D1D}. The distinction
  between inner and outer surroundings are motivated by the
  potentially large jumps in material parameters such as the
  convective velocity or diffusion coefficient between the inner and
  outer compartments in applications. Such jumps would be challenging
  to represent in an extended domain in the 3D-1D setting. These
  models are thus particularly relevant for perivascular transport with a
  vascular-perivascular barrier, such as e.g.~the blood-brain barrier
  (BBB) in the human brain. However, the case of vascular-perivascular
  exchange may also be highly relevant e.g.~in connection with a leaky
  BBB or transport of substances across the BBB. Therefore, we address
  the extended 3D-3D-3D and 3D-1D-1D problem setting representing
  coupled tissue, perivascular and vascular transport separately in
  Section~\ref{sec:3D1D1D}.

In Section~\ref{sec:networks}, we will also consider an extension of
this setting to networks of vessels. We will then consider a network
of $N$ domains $\Omega_{v,i}$ with center-curves $\Lambda_{i} =
\{\bm{\lambda}_i (s), \,\,s \in (0,L_i)\}$ for $ i = 1,\ldots,
N$. Extending upon the notation introduced above, we then denote $A_i
= |\Theta_i|$ and $P_i = |\partial \Theta_{2,i}|$ where $\Theta_i$ is
the cross-section of $\Omega_{v,i}$ and $\partial \Theta_{2,i}$ is the
outer boundary of $\Omega_{v,i}$.

\section{Transport by convection and diffusion in a moving domain}
\label{sec:transport}

We are interested in analyzing the coupled transport of a solute in a
moving domain governed by diffusion and convection in general, and in
a moving vessel and its surroundings in particular. To this end, we
introduce a system of coupled convection-diffusion equations
(Section~\ref{sec:3d3d_reference}). We may directly consider a more
general geometrical setting
(Section~\ref{setting:geometry_moving_abstract}) for the weak
formulation (Section~\ref{sec:3D-3D-weak}) to show that such solutions
exist (Section~\ref{sec:3D-3D-existence}, \Cref{prop:exist_3d_3d}).

\subsection{System of convection-diffusion equations in and around a moving vessel}
\label{sec:3d3d_reference}

We consider a moving vessel and assume that the vessel motion
$\Omega_v(t)$, and convective velocity fields $\bm{u}_v(t) :
\Omega_v \rightarrow \R^3$ and $\bm{u}_s(t) : \Omega_s
\rightarrow \R^3$ are prescribed for $t \in [0, T]$. Our coupled
three-dimensional transport boundary-value problem in an Eulerian
frame reads as: for a.e.~$t$, find the solute concentrations $c_v(t) :
\Omega_v(t) \rightarrow \R$ and $c_s (t) : \Omega_s(t) \rightarrow \R$
such that the following governing equations, interface conditions,
boundary conditions and initial conditions hold:
\begin{subequations}
  \begin{alignat}{2}
    \partial_t c_v - \nabla \cdot (D_v \nabla c_v ) + \nabla \cdot (\bm{u}_v c_v)  &= f_v, && \quad  \mathrm{in} \, \Omega_v(t) \times (0,T], \label{eq:3d_model_1} \\ 
    \partial_t c_s - \nabla \cdot (D_s \nabla c_s ) + \nabla \cdot (\bm{u}_s c_s)  &= f_s, && \quad  \mathrm{in} \, \Omega_s(t) \times (0,T], \label{eq:3d_model_2} \\ 
    (c_v \tilde{\bm{u}}_v - D_v \nabla c_v) \cdot \bm{n}_v -\xi(c_v - c_s)&= 0,    && \quad \mathrm{on} \, \Gamma(t)\times (0,T], \label{eq:3d_model_3} \\ 
    (c_v \tilde{\bm{u}}_v - D_v \nabla c_v) \cdot \bm{n}_v + (c_s \tilde{\bm{u}}_s - D_s \nabla c_s) \cdot \bm{n}_s & = 0 ,  && \quad \mathrm{on} \, \Gamma(t) \times (0,T],  \label{eq:3d_model_4}\\ 
    (c_v \tilde{\bm{u}}_v -D_v \nabla c_v )\cdot \bm{n}_v =  (c_s \tilde{\bm{u}}_s - D_s\nabla c_s )\cdot \bm{n}_s & = 0, && \quad \mathrm{on} \, (\Gamma_0(t) \cup \Gamma_L(t))\times (0,T],   \label{eq:3d_model_5}\\ 
    c_s & = 0, && \quad \mathrm{on} \,  \partial \Omega \times (0,T],  \label{eq:3d_model_6} \\
    c_v(0) = c_v^0 , \quad  \mathrm{in}\, \Omega_v(0), \quad c_s(0) &= c_s^0  && \quad \mathrm{in} \, \Omega_s(0) . \label{eq:3d_model_7}
  \end{alignat}%
 \label{eq:3d-3d}%
\end{subequations}%
If $\partial \Omega_v \backslash \Gamma \neq \emptyset$ $(R_1 > 0$, an annular domain $\Omega_v$), we also impose the boundary condition: 
\begin{alignat}{2}
(c_v \tilde{\bm{u}}_v - D_v \nabla c_v) \cdot \bm{n}_v &= 0    && \quad \mathrm{on} \, \partial \Omega_v(t) \backslash \Gamma(t) \times (0,T]. \label{eq:3d_model_3_additional}
\end{alignat}  
The time derivatives in the above formulation are the Eulerian time derivatives. The parameters $D_v$, $D_s$ are given diffusion tensors in
$\Omega_v$ and $\Omega_s$ respectively, while $f_v$ and $f_s$ are
given source functions. For $ i \in \{v,s\}$, the relative (net) velocity $\tilde{\bm{u}}_i = \bm{u}_i - \bm{w}$ accounts for the velocity of the domain $\bm{w}$, defined below cf.~\eqref{eq:domain_velocity}.  The interface condition \eqref{eq:3d_model_3}
models the lateral interface between the vessel and its surroundings
$\Gamma$ as a semi-permeable membrane with permeability $\xi$, while
the auxiliary condition \eqref{eq:3d_model_4} enforces conservation of
mass. At the vertical boundaries, the condition \eqref{eq:3d_model_5}
stipulates no flux, while we keep the concentration fixed and zero
(for simplicity) at the outermost boundary $\partial \Omega$
via~\eqref{eq:3d_model_6}. The last relations~\eqref{eq:3d_model_7}
define the initial conditions with given initial states $c_v^0 :
\Omega_v \rightarrow \R, \,  c_v^0 \in L^2(\Omega
_v)$ and $c_s^0 : \Omega_s \rightarrow \R, \, c_s^0 \in L^2(\Omega_s)$.

\subsection{Observations on the domain velocity}
\label{setting:geometry_moving_abstract}

  For the existence result we can weaken our geometrical assumptions on the domains.
  Precisely, we let $\Omega(0)\subset\mathbb R^d$ be a Lipschitz domain, i.e., open, connected
  and with a Lipschitz boundary and we assume that $\Omega_v(0) \subset \Omega(0)$ is itself a Lipschitz domain and
  compactly contained in $\Omega(0)$. In particular, it holds $\operatorname{dist}(\partial\Omega_v(0), \partial\Omega(0)) > 0$.
  We measurably partition $\partial\Omega_v(0)$ into two sets that play the role
  of $\Gamma(0)$ in \eqref{eq:3d_model_3} and \eqref{eq:3d_model_4} and $\Gamma_0(0) \cup \Gamma_L(0)$ in \eqref{eq:3d_model_5}, respectively
   (and will be denoted by the same symbols).

  We define moving domains according to the velocity method, see \cite{delfour2011shapes}. More precisely,
  assume that the domain velocity $\bm{w}:\mathbb R^d\times \mathbb R \to \mathbb R^d$ is smooth and compactly supported. We denote
  by $\bm{\psi}: \mathbb R^d\times \mathbb R \to \mathbb R^d$ the flow map of the order
  \begin{equation}
    \begin{split}        \partial_t \bm{\psi} (\bm{x},t) &= \bm{w} (\bm{\psi}(\bm{x},t), t),
    \\
    \bm{\psi}(\bm{x}, 0) &= \bm{x}.
  \end{split} \label{eq:domain_velocity}
  \end{equation}
  Standard ODE theory implies that $\bm{\psi}\in C^\infty(\mathbb R^{d+1})$ and for all fixed $t\in[0,T]$ the map
  \begin{align*}
    \bm{\psi}_t:\mathbb R^d \to \mathbb R^d, \quad \bm{x} \mapsto \bm{\psi}(\bm{x},t)
  \end{align*}
  is a diffeomorphism. In this notation, the connection between our reference domain and 
  the domains to a later time is given by
  \begin{gather*}
    \Omega_s(t) = \bm{\psi}_t(\Omega_s(0)), \quad \Omega_v(t) = \bm{\psi}_t(\Omega_v(0)).
  \end{gather*}
  As $\Omega_s(0)$ and $\Omega_v(0)$ are open, connected and Lipschitz it holds $\partial \Omega_v(t) = \bm{\psi}_t(\partial\Omega_v(0))$
   and $\partial\Omega_v(t)$ has Lipschitz boundary, see \cite{hofmann2007geometric}. An important observation which links the above definitions to the specific setting in subsection~\ref{sec:geometry_specific} is now in order. Denoting by $\mathrm{det}(D \bm{\psi}_t)$ the determinant of the Jacobian matrix of $\bm{\psi}_t$, it holds that \cite[Section 1.1.1]{nobile2001numerical}: 
   \begin{multline}
    \partial_t A(s,t) = \partial_t \int_{\Theta(s,t)} 1  = \partial_t  \int_{\Theta(s,0)} |\mathrm{det}(D \bm{\psi}_t)| \\  = \int_{\Theta(s,0)} \bm{\psi}_{-t} (\nabla \cdot \bm{w})  |\mathrm{det}(D \bm{\psi}_t) |= \int_{\Theta(s,t)} \nabla \cdot \bm{w} = \int_{\partial \Theta(s,t)} \bm{w} \cdot \bm{n}.  \label{eq:transformation_determinant}
   \end{multline}
  
\subsection{A weak formulation of the coupled 3D-3D transport model}
\label{sec:3D-3D-weak}

Let $i\in\{s,v\}$, for fixed $t\in I:=[0,T]$ we set $X_s(t) = H^1_{\partial \Omega}(\Omega_s(t)): = \{ c \in H^1(\Omega_s(t)), \, c \vert_{\partial \Omega} = 0 \}$, and $X_v(t) = H^1(\Omega_v(t))$  and $H_i(t) = L^2(\Omega_i(t))$.
Further, we abbreviate $X_i = (X_i(t))_{t\in I}$ and $H_i = (H_i(t))_{t\in I}$. To relate the function spaces at time $t$
to the reference time (and vice versa) we use the pushforward induced by $\bm{\psi}_t$, and we define:
\begin{equation*}
  \phi_t :X_i(0) \to X_i(t), \quad  \phi_t c_i = c_i \circ \psi_t^{-1}
\end{equation*}
with inverse $\phi_{-t} = \phi_t^{-1}$ given by $\phi_{-t}c_i = c_i \circ \psi_t$. By the
chain rule for Sobolev spaces it can be seen that for all $t\in[0,T]$ the maps
$\phi_t$ are linear homeomorphisms. Now, to define a function space framework we follow \cite{alphonse2015abstract}
and set 

\begin{align*}
  L^2_{X_i} &= \{ c_i:[0,T] \to \bigcup _{t\in[0,T]}X_i(t)\times \{t\}, t\mapsto (\bar c_i(t),t) \mid \phi_{-(\cdot)}\bar c_i(\cdot)\in L^2(0,T; X_i(0)) \},
  \\
  L^2_{X_i'} &= \{ f_i:[0,T] \to \bigcup _{t\in[0,T]}X'_i(t)\times \{t\}, t\mapsto (\bar f_i(t),t) \mid \phi^*_{(\cdot)}\bar f_i(\cdot)\in L^2(0,T; X'_i(0)) \},
\end{align*}
where $\phi^*_t: X_i(t)' \to X_i(0)'$ denotes the adjoint map to $\phi_t$. The above spaces are equipped with the norms: 
\begin{align*}
\forall c_i \in L^2_{X_i}, \,\, \|c_i\|^2_{L^2_{X_i}} = \int_0^T \|c_i\|^2_{X_i(t)}, \quad \forall f_i \in L^2_{X_i'}, \,\, \|f_i\|^2_{L^2_{X'_i}} = \int_0^T \|f_i\|^2_{X'_i(t)}.
\end{align*}
Next, we define a weak material derivative, where we specialize the abstract definition of \cite{alphonse2015abstract} to our case.
We say that a function $c_i \in L^2_{X_i}$ has a weak material derivative $\dot c_i \in L^2_{X_i'}$ if it holds 
\begin{equation}\label{eq:weak_material}
  \int_0^T \langle \dot c_i, \eta \rangle_{X_i'(t)} \, \mathrm dt
  = 
  - \int_0^T \int_{\Omega_i(t)}c_i\dot\eta \, \mathrm dx\mathrm dt
  -
  \int_0^T \int_{\Omega_i(t)} c_i \eta \nabla\cdot \bm{w}\,\mathrm dx \mathrm dt, 
  \quad 
  \forall \eta \in \mathcal{D}_{X_i}(0,T),
\end{equation}
where $\mathcal D_{X_i}$ is the subset of $L^2_{X_i}$ such that $t\mapsto \phi_{-t}\eta$ is a member of $\mathcal D(0,T; X_i(0))$.
We are now in a position to define the Sobolev space $W(X_i, X_i')$ used for existence theory
\begin{equation*}
  W(X_i, X_i') = \{ c_i \in L^2_{X_i} \mid \dot c_i \in L^2_{X_i'} \}.
\end{equation*}
As in the classical case, this space embeds into $C_{H_i}^0(0,T)$, which is defined similarly as $\mathcal D_{X_i}(0,T)$ above 
and thus initial value problems can be formulated meaningfully.

\begin{remark}[Connection to strong material derivative]
  For smooth functions $c$ the above definition agrees with the Arbitrary Lagrangian Eulerian (ALE) framework \cite[Section 1.1]{nobile2001numerical}, and it holds that 
  \begin{equation*}
    \dot c = \phi_t (\, \partial_t \,  (  \phi_{-t} c)).  
  \end{equation*}
  By the chain rule, it then follows for smooth functions, that 
  \begin{equation}
    \dot c(\bm{x},t) = \partial_t c(\bm{x},t) + \nabla c(\bm{x},t) \cdot  \bm{w}(\bm{x},t), \quad (\bm{x} , t) \in \Omega(t) \times (0,T).
    \label{eq:strong_mater_derivative}
  \end{equation} 
\end{remark}
Replacing the Eulerian
time derivative via the definition of the material time derivative,
and using the standard identity
\begin{equation*}
  \nabla c \cdot \bm{w} = \nabla \cdot (\bm{w} c) - (\nabla \cdot \bm{w}) c,
\end{equation*}
we can rephrase~\eqref{eq:3d_model_1}--\eqref{eq:3d_model_2} as
\begin{subequations}
  \begin{align}
    \dot{c}_v + \nabla \cdot \bm{w} c_v - \nabla \cdot ( D_v \nabla c_v) + \nabla \cdot \left ( (\bm{u}_v - \bm{w}) c_v  \right ) &= f_v,
    \quad \mathrm{in} \quad \Omega_v(t) \times (0,T], \\
      \dot{c_s} + \nabla \cdot \bm{w} c_s - \nabla \cdot (D_s \nabla c_s) + \nabla \cdot \left ( (\bm{u}_s - \bm{w}) c_s \right ) &= f_s,
      \quad \mathrm{in} \quad \Omega_s(t) \times (0,T].
  \end{align}
\end{subequations}%

To formulate a coherent weak formulation for the system of coupled equations, we introduce the following product spaces and their respective norms (written for $\bm{\phi} = (\phi_v, \phi_s)$)
\begin{subequations}
  \begin{alignat*}{2}
    \bm{V}(t) &= H^1(\Omega_v(t)) \times H^1_{\partial \Omega}(\Omega_s(t)), \,\,  &&  \,\,  \|\bm{\phi}\|_{\bm{V}(t)}^2 = \| \phi_v\|_{H^1(\Omega_v(t))}^2 +  \| \phi_s\|_{H^1(\Omega_s(t))}^2, \\
    \bm{H}(t) &= L^2(\Omega_v(t)) \times L^2(\Omega_s(t)),\,\,  && \,\, \|\bm{\phi}\|_{\bm{H}(t)}^2 = \| \phi_v\|_{L^2(\Omega_v(t))}^2 +  \| \phi_s\|_{L^2(\Omega_s(t))}^2. 
  \end{alignat*}
\end{subequations}
Similarly, we define the product space:
  \begin{alignat}{2}
    \bm{W} & =\{\bm{w} = (w_v, w_s), 
    \dot{\bm{w}} = (\dot{w}_v, \dot{w}_s): \,\, w_v \in W(X_v,X_v'), \,\, w_s \in W(X_s , X_s') \}, 
  \end{alignat}
equipped with the norm 
\begin{equation}
  \|\bm{w}\|_{\bm{W}}^2 = \sum_{i \in \{v,s\}} (\|w_i\|_{L^2_{X_i}}^2 + \|\dot{w}_i\|_{L^2_{X_i'}}^2   ).
\end{equation}
The weak formulation for \eqref{eq:3d-3d} then reads: find $\bm{c} =
(c_v, c_s) \in \bm{W}$ such that for all $\bm{\varphi} =
(\varphi_v,\varphi_s) \in \bm{V}(t)$,
\begin{multline}
    \langle \dot{\bm{c}}(t), \bm{\varphi} \rangle_{\bm{V}(t)} + \lambda(t;\bm{c}(t), \bm{\varphi})
    + \mathcal{A}(t; \bm{c}(t), \bm{\varphi}) + \mathcal{B}(t; \bm{c}(t), \bm{\varphi}) \\
    = (\bm f_v(t),\varphi_v)_{\Omega_v } + (\bm f_s(t) ,\varphi_s)_{\Omega_s},  
\label{eq:weak_form_3d3d}%
\end{multline}
complemented by the initial condition
\begin{equation*}
  \bm{c}(0) = (c_v^0, c_s^0) \in \bm{H}(0), 
\end{equation*}
where for any $\bm{c} = (c_v,c_s) \in \bm{V}(t)$ and $\bm{\varphi} = (\varphi_v, \varphi_s) \in \bm{V}(t)$ we have the bilinear forms: 
\begin{align*}
  \lambda(t;\bm{c},\bm{\varphi}) & = (\nabla \cdot \bm{w} c_v , \varphi_v)_{\Omega_v(t)} + (\nabla \cdot \bm{w} c_s ,\varphi_s)_{\Omega_s(t)},\\
  \mathcal{A}(t;\bm{c},\bm{\varphi}) & =  (D_v \nabla c_v - (\bm{u}_v - \bm{w}) c_v , \nabla \varphi_v)_{\Omega_v(t)} + (D_s \nabla c_s - (\bm{u}_s- \bm{w})c_s, \nabla \varphi_s)_{\Omega_s(t)}, \\  
\mathcal{B}(t;\bm{c}, \bm{\varphi}) &= (\xi (c_v - c_s), \varphi_v)_{\Gamma(t)} + (\xi(c_s - c_v), \varphi_s)_{\Gamma(t)}. 
\end{align*}

\subsection{Well-posedness of the convection-diffusion problem over a moving domain}
\label{sec:3D-3D-existence}


We then obtain the following result for the existence and well-posedness of weak solutions.
\begin{prop}
  \label{prop:exist_3d_3d}
  Assume the geometrical setting of Section \ref{setting:geometry_moving_abstract} and let $\xi\in L^\infty(0,T;L^\infty(\Gamma(t)))$ with $\xi \geq 0$.
   Further assume that $D_i\in L^\infty(0,T; L^\infty(\Omega_i(t), \mathbb R^{d\times d}))$ with a uniform ellipticity constant $\nu>0$ and $\bm u_i \in L^\infty(0,T; L^\infty(\Omega_i(t)))$. Then, for every 
  $\bm c_0 = (c_v^0, c_s^0)\in \bm H(0)$ and $\bm f = (f_v, f_s)\in L^2_{X_v'} \times L^2_{X_s'}$, there exists a unique 
  solution $\bm c = (c_v,c_s) \in \bm W$ to~\eqref{eq:weak_form_3d3d}. Further, there exists a constant $C$ such that 
  \[ \|\bm c\|_{\bm W} \leq C ( \| f_v \|_{L^2_{X_v'}} + \|f_s\|_{L^2_{X_s'}}  + \|\bm{c}_0\|_{\bm{H}(0)}) ) .\]
\end{prop}

\begin{proof}
  We verify the assumptions of the abstract framework given in \cite{alphonse2015abstract}. 
  These can be grouped in two sets of requirements, one set of assumptions concerns the level of 
  smoothness that must be imposed on the moving domains - in the notation of \cite{alphonse2015abstract}
  these are Assumption 2.17, 2.24 and Assumption 2.31. On the other hand we need standard assumptions on
  the involved operators which are summarized in Assumption 3.3 of \cite{alphonse2015abstract}.
   
  \paragraph*{\textit{Verifying the smoothness assumptions of the moving domains}.}
  Let $\bm c \in \bm V(0)$, then 
  by the transformation formula it holds
  \begin{equation}\label{eq:proof_well_posedness_integral_trafo}
    t \mapsto \| \phi_t \bm c \|^2_{\bm V(t)} 
    =
    \sum_{i\in \{v,s\} } \int_{\Omega_i(0)} (c_i^2 + |\nabla c_i|^2)|\det (D\psi_t)|\, \mathrm dx. 
  \end{equation}
  As $(x,t)\mapsto \psi(x,t)$ is smooth and $\psi_t$ is a diffeomorphism
  we know that $D\psi_t$ is invertible everywhere in $\overline\Omega$ and thus
  $|\det (D\psi_t)|$ is bounded away from zero. Using the smoothness of $\psi$ with respect
  to the temporal variable implies that this bound is independent of time. Hence, $t\mapsto \|\phi_t \bm c\|_{\bm V(t)}$ is continuous as required in Assumption 2.17.

  To show Assumption 2.24, we need prove that 
  \begin{equation*}
    t\mapsto \theta(t,\bm c) 
    :=
    \sum_{i\in \{v,s\} } \int_{\Omega_i(0)} c_i^2 |\det (D\psi_t)|\, \mathrm dx
  \end{equation*}
  is classically differentiable. As mentioned above, $(x,t)\mapsto \psi(x,t)$ is smooth
  and so is $(x,t)\mapsto |\det (D\psi_t)(x)|$ which allows us, resorting to Lebesgue's dominated convergence theorem, to differentiate under the integral sign.
  Further, for $\bm c^1, \bm c^2 \in \bm V(0)$ we estimate using the boundedness of $(x,t)\mapsto|\det(D\psi_t(x))|$
  \begin{equation*}
    |\theta(t, \bm c^1 + \bm c^2) - \theta(t, \bm c^1 - \bm c^2)|
    =
    \sum_{i\in \{v,s\} } \int_{\Omega_i(0)} |c^1_i c^2_i|\,|\det(D\psi_t)|\, \mathrm dx
    \leq C
    \| \bm c^1 \|_{\bm H(0)} \| \bm c^2 \|_{\bm H(0)}
  \end{equation*}
  for some constant $C$. This completes the requirements of Assumption 2.24. 
  
  Concerning Assumption 2.31 of \cite{alphonse2015abstract}, note that the map
  $T_t$ defined in equation (2.7) of this paper is in our case given by
  \begin{equation*}
    T_t: \bm H(0) \to \bm H(0), \quad \bm c = (c_v, c_s) \mapsto (c_v |\det(D\psi_t)|, c_s |\det (D\psi_t)|)
  \end{equation*}
  and as $|\det(D\psi_t)|$ is smooth and bounded away from zero it holds that
  \begin{equation*}
    \bm c \in \bm V(0) \quad \Leftrightarrow \quad T_t\bm c \in \bm V(0).
  \end{equation*}
  By Remark 2.34 in \cite{alphonse2015abstract} this guarantees that Assumption 2.31 therein holds.

  \paragraph*{\textit{Properties of the PDE Operators}}
  We now verify the coercivity and continuity properties of the bilinear forms. 
  We must show that for a.e.~$t$, there exist constants $K_1, K_2$ and $K_3$ independent of $t$ such that 
\begin{alignat}{3}
\mathcal{A}(t;\bm{c},\bm{c}) + \mathcal{B}(t;\bm{c}, \bm{c} ) &\geq K_1 \|\bm{c}\|_{\bm{V}(t)}^2 -  K_2 \|\bm{c}\|_{\bm{H}(t)}^2 && \quad \foralls \bm{c} \in \bm{V}(t),  \label{eq:coercivity_3d3d} \\ 
|\mathcal{A}(t;\bm{c},\bm{\varphi}) + \mathcal{B}(t;\bm{c}, \bm{\varphi}) |& \leq K_3 \|\bm{c}\|_{\bm{V}(t)} \|\bm{\varphi}\|_{\bm{V}(t)} && \quad \foralls \bm{c}, \bm{\phi}\in \bm{V}(t). \label{eq:cont_3d3d}
\end{alignat}
Using that $\bm u_v, \bm u_s$ and $\bm w \in L^\infty(\Omega_v(t))$
with a norm bound independent of $t\in[0,T]$ and that $D_v, D_s$ are uniformly elliptic
with ellipticity constant $\nu$ independent of time, we may estimate using Young's and
H\"older's inequality for $\bm c = (c_v, c_s)$

\begin{align*}
  \mathcal A(t; \bm c, \bm c) 
  & 
  \geq
  \nu  \sum_{i\in \{ s,v \}  } \| \nabla c_i \|^2_{L^2(\Omega_i(t))} 
  -
  \sum_{i\in \{ s,v \}  } \| \bm u_i - \bm w \|_{L^\infty(\Omega_i(t))} \| c_i \|_{L^2(\Omega_i(t))} \| \nabla c_i \|_{L^2(\Omega_i(t))}
  \\
  & 
  \geq 
  \frac{\nu}{2} \sum_{i\in \{ s,v \}  } \| \nabla c_i \|^2_{L^2(\Omega_i(t))}
  -
  \frac{1}{2\nu} \sum_{i\in \{ s,v \}  } \| \bm u_i - \bm w \|^2_{L^\infty(\Omega_i(t))} \| c_i \|^2_{L^2(\Omega_i(t))}
  \\
  &\geq
  \frac{\nu}{2}\|\bm c\|^2_{\bm V(t)} - \max_{i\in \{ s,v \}}\left( \frac{\|\bm u_i - \bm w\|_{L^\infty(\Omega_i(t))}}{2\nu}, \frac{\nu}{2} \right) \|\bm c\|^2_{\bm H (t)}.
\end{align*}
Using that $\xi \geq 0$ it is readily seen that $\mathcal B(t; \bm c, \bm c) \geq 0$, in fact it holds that 
\begin{equation*} 
  \mathcal B(t; \bm c, \bm c) = \int_{\Gamma(t)}(c_v - c_s)^2\xi \, \mathrm ds \geq 0.
\end{equation*}
For the continuity property, we note that the trace constant used to handle 
$\mathcal{B}$ is independent of $t$ since for any $c_i \in H^1(\Omega_i(t))$ it holds that 
\begin{align} 
  \|c_i\|_{L^2(\Gamma(t))} = \||\det(D(\psi^{-1}_t))|^{1/2} \phi_{-t} c_i\|_{L^2(\Gamma(0))} \leq C_0 \|\phi_{-t}c_i\|_{H^1(\Omega_i(0))} \leq C_{1}  \|c_i\|_{H^1(\Omega_i(t))} \label{eq:trace_constant_indep_t}
\end{align}
for some constants $C_0$, $C_1$. The above holds from the trace inequality on $\Gamma(0)$ and from the continuity bound 
of the map $\phi_{-t}$ which is independent of $t$. The continuity bound \eqref{eq:cont_3d3d} then immediately 
follows. Therefore, as all the assumptions of~ \cite[Theorem 3.6]{alphonse2015abstract} hold, the stated result follows.
\end{proof}

\section{Coupled 3D-1D formulations for solute transport models}
\label{sec:3D1D}

Our next objective is to derive geometrically-explicit but
dimensionally-reduced representations of the coupled solute transport
models introduced and established in the previous
(\Cref{sec:transport}). We first derive transport equations describing
the cross-section average concentration in each vessel network segment
(\Cref{subsec:derv_1d}) and their variational formulation
(\Cref{subsec:derivation_1dmodel}). Conversely, the solute transport
equations are extended accordingly; from the surrounding to the
complete domain (\Cref{subsec:weak_form_3d_derivation}). The full
coupled variational problem is well-posed (\Cref{sec:coupled_3d_1d}),
and can be extended to vascular networks (\Cref{sec:networks}).
We begin by making assumptions on the material parameters mainly to simplify the presentation. We will adopt these assumptions in the remainder of this paper.  

\subsection{Assumptions on material parameters}
The parameter $D_v$ is assumed to  be single valued function rather than a tensor, and $D_v \in L^{\infty}(0,T;L^{\infty}(\Omega_v(t)))$. The parameter  $D_s \in L^{\infty}(0,T;L^{\infty}(\Omega_s(t), \R^{d \times d}))$ with uniform ellipticity constant $\nu > 0$.  In addition,  $\xi$  and $D_v$ are assumed to be constant in each cross-section $\Theta(s,t)$, $(s,t) \in \Lambda \times (0,T)$.   Finally, we assume that the  velocity fields $\bm{u}_i \in L^{\infty}(0,T;H^2(\Omega_i)^3)$  for $i \in \{v,s\}$. 

\subsection{Derivation of a vessel-averaged (1D) transport equation}
\label{subsec:derv_1d}

The aim of this section is to derive a one-dimensional model for the
cross-section average of the concentration $c_v$. Recalling the
cross-sections $\Theta(s)$ with area $A = A(s)$, we define the
cross-section average for $s \in \Lambda$ by
\begin{equation*}
 \avg{f}(s) = \frac{1}{A(s)} \int_{\Theta(s)} f, \quad \foralls f \in L^1(\Theta(s)) .
\end{equation*}
Analogously, recalling the cross-section boundary $\partial \Theta_2(s)$ with (lateral cross-section) perimeter $P = P(s)$, we set:
\begin{equation}
  \bar{f}(s) = \frac{1}{P(s)} \int_{\partial \Theta_2(s)} f, \quad \foralls f \in L^1(\partial \Theta_2(s)).
  \label{eq:bar}
\end{equation}
For the derivation, we rely on the following assumptions on the vessel geometry and vessel deformations (adapted
from~\cite[Chapter 2]{d2007multiscale}, and \cite{laurino2019derivation}).  Assumption \ref{assumption:shape_profiles} is needed in the derivation of the reduced 1D model, see Proposition \ref{subsec:derivation_1dmodel}, and Assumption \ref{assumption:boundary_condition_A} is used in the derivation of its variational formulation, see Section~\ref{subsec:derv_1d}. 
\begin{assumption}[Averages and shape profile]
  \label{assumption:shape_profiles}
  Assume the following.
  \begin{itemize} 
    \item
      For $c_v : \Omega_v \times (0, T) \rightarrow \mathbb{R}$ and $c_s : \Omega_s \times (0, T) \rightarrow \mathbb{R}$ solving~\eqref{eq:3d-3d}, 
      the (lateral) cross-section averages are well-defined i.e.~$c_v(t) \in L^1(\Theta(s)) \cap L^1(\partial \Theta(s))$ and $c_s(t) \in L^1(\partial \Theta(s))$ for all $s \in \Lambda$ and $t \in (0, T)$.
    \item
      Further, there exists a shape function $w_c = w_c(r)$ in the radial variable $r$ only, with $\avg{w_c} = 1$ and such that the following splitting holds: for all $(s,r,\theta,t) \in \Omega_v(t) \times (0,T]$,     
    \begin{equation*}
      c_v(s, r, \theta, t) = \avg{c_v}(s,t) \, w_c(r). 
    \end{equation*}
  \end{itemize}
\end{assumption}
\begin{assumption}[Conditions on the vessel geometry and deformation]
  \label{assumption:boundary_condition_A} Assume the following: 
  \begin{equation}
    \partial_s R_2^2 = \partial_s R_1^2 = 0, \quad \mathrm{on} \,\, \Gamma_0 \cup \Gamma_L.      \label{eq:A_s}
  \end{equation}
  The above is adapted from \cite{laurino2019derivation}. In fact, if $R_1$ and $R_2$ are independent of $\theta$ or if $w_c = 1$, then we can relax the above assumption by only requiring that 
  \[ \partial_s A = 0, \quad \mathrm{on} \,\, s = 0 , L, \] 
  since it will be sufficient for the derivation of our weak formulation, see subsection \ref{subsec:derv_1d}.
\end{assumption}
The next proposition states a one-dimensional transport equation for
the average concentration $\cc$ along the vessel centerline $\Lambda$ and over
time $t \in (0, T)$.
\begin{prop}[1D transport equation]
  \label{prop:derivation_1d}
  Under Assumption \ref{assumption:shape_profiles}, the cross-section average concentration $\cc = \avg{c_v}$ satisfies the following equation in $\Lambda \times (0, T)$:
  \begin{equation}
    \partial_t (A\cc) - \partial_{s} \left( D_v A \partial_s \cc \right) + \partial_s \left( A \avg{u_{v,s} w_c} \cc \right)
    + \xi P \left( \overline{w_c} \cc - \overline{c_s} \right) +
    G(\cc) = A \avg{f_v} ,
    \label{eq:derivation_final}
  \end{equation}
  where $u_{v,s}$ is the axial component of the velocity $\bm{u}_v$, and where we have introduced the auxiliary expressions
  \begin{align}
   G(\cc) &= G(R_1, R_2, w_c)(\cc) = -  \partial_s \left (  D_v \vtwo (R_1, R_2, w_c) \cc \right ) , \label{eq:G}\\
   \vvv_s(R_1, R_2, w_c) &= \sum_{i=1}^{2} - \frac{\gamma_i}{2} \int_0^{2\pi} \partial_s R_i^2 (1-w_c(R_i)), \quad \gamma_1 = 1, \gamma_2 = -1 .
    \label{eq:def_rho1} 
  \end{align}
\end{prop}
Before proceeding with the proof of~\Cref{prop:derivation_1d}, we make two remarks.

\begin{remark}
Recall that the functions $A = A(s,t)$ and $\avg{ u_{v,s} w_c} = \avg{ u_{v,s} w_c}(s,t)$ denote the cross-sectional area and a weighted average axial velocity, respectively. These functions can be either a-priori determined or solved for via reduced flow models, such as e.g.~reduced blood flow models~\cite{vcanic2003mathematical}, perivascular fluid flow models~\cite{daversin2022geometrically}, root water uptake models~\cite{koch2018new}, or geothermal wells~\cite{gjerde2020singularity} as appropriate for the problem setting.
\end{remark}
\begin{remark}
  If $\Omega_v$ is a cylinder (representing for instance a blood vessel, reservoir well or plant root but not a perivascular space), $R_1(s, \theta, t) = 0$ and $R_2(s, \theta, t) = R(s, t)$. In this case, if we also assume that $w_c(r) = w_c(R) = 1$, then $G = 0$ and~\eqref{eq:derivation_final} simplifies to:
  \begin{equation*}
    \partial_t (A\cc ) - \partial_s(D_v A \partial_s \cc) +  \partial_s(A \avg{u_{v,s}} \cc) +  \xi P (\cc- \bar{c}_s) = A \avg{f_v}. 
  \end{equation*} 
\end{remark}

\begin{proof}[\textit{Proof of Proposition \ref{prop:derivation_1d}}] 
    We proceed via a similar approach as in~\cite{laurino2019derivation}. Namely, consider two arbitrary points $s_1$ and $s_2$ with $0 \leq s_1 < s_2 \leq L$. Let $\mathcal{P} = \mathcal{P}(t)$ denote the portion of  $\Omega_v(t)$ bounded by two cross-sections $\Theta(s_1)$ and $\Theta(s_2)$ perpendicular to $\Lambda$,  and let $\Gamma_\mathcal{P}$ denote the lateral boundary of $\mathcal{P}$. To simplify notation, we drop the subscript $v$ in \eqref{eq:3d_model_1}. We now integrate \eqref{eq:3d_model_1} over $\mathcal{P}$ omitting the integration measures when self-evident. 
\begin{equation*}
  \int_\mathcal{P} \partial_t c -  \int_\mathcal{P} \nabla \cdot (D \nabla c ) +  \int_\mathcal{P} \nabla \cdot (\bm{u} c) - \int_\mathcal{P} f
  := \mathcal{T}_1 +\mathcal{T}_2 + \mathcal{T}_3 + \mathcal{T}_4
  = 0. 
\end{equation*}
For $\mathcal{T}_1$, we have by Reynolds transport theorem accounting for the domain velocity $\bm{w}$ and by definition of the cross-section average, see e.g \cite{nobile2001numerical,alphonse2015abstract}, that
\begin{equation}
  \mathcal{T}_1
  = \int_{\mathcal{P}} \partial_t c
  = \partial_t \int_{\mathcal{P}} c - \int_{\partial \mathcal{P}} c \bm{w} \cdot \bm{n}
  = \int_{s_1}^{s_2 } \partial_t (A \avg{c}) - \int_{\partial \mathcal{P}} c \bm{w} \cdot \bm{n} .
  \label{eq:T1}
\end{equation}
For $\mathcal{T}_2$, using the divergence theorem, we have that 
\begin{multline}
  \mathcal{T}_2
  = -  \int_\mathcal{P} \nabla \cdot (D \nabla c )
  = - \int_\mathcal{\partial P}  D \nabla c  \cdot \bm{n} \\
  =  \int_{ \Theta(s_1)} D \partial_s c - \int_{ \Theta(s_2)} D \partial_s c  - \int_{s_1}^{s_2} \int_{\partial \Theta(s)}  D \nabla  c \cdot \bm{n} .
  \label{eq:T2_zero}
\end{multline}
Following \cite{laurino2019derivation}, we write the first and second terms above as follows.
\begin{multline*}
  \int_{ \Theta(s_1)} D  \partial_s c \, - \int_{ \Theta(s_2)} D \partial_s c
  = - \int_{s_1}^{s_2} \frac{\partial}{\partial s } \int_{ \Theta(s)} D \partial_s c \ds
  =  -  \int_{s_1}^{s_2}  \frac{\partial}{\partial s} \int_0^{2\pi} \int_{R_1}^{R_2} D \partial_s c r  \dr \dtheta \ds . 
\end{multline*}
Using the assumption that $D$ is constant on each cross-section, recalling that $\gamma_{1} = 1$ and $\gamma_{2}=-1$, and applying Leibniz's rule yield 
\begin{multline*}
  \int_0^{2\pi}\int_{R_1}^{R_2} D \partial_s cr \dr \dtheta
  = D  \partial_s \int_0^{2\pi} \int_{R_1}^{R_2} c r \dr \theta + \sum_{i=1}^2 \frac{\gamma_i}{2} \int_0^{2\pi} D c(R_i) \partial_s R^2_i \dtheta \\
  =  D \partial_s (A \avg{c}) + \sum_{i=1}^2 \int_0^{2\pi}\frac{\gamma_i}{2} D c(R_i) \partial_s R^2_i \dtheta.     
\end{multline*}
Thus, we write $\mathcal{T}_2$ as 
\begin{equation}
  \label{eq:T2} 
  \mathcal{T}_2
  =  \int_{s_1}^{s_2} \left(- \partial_s ( D \partial_{s}(A \avg{c})) - \sum_{i=1}^2 \int_0^{2\pi} \frac{\gamma_i}{2} \partial_s( D c(R_i) \partial_s R^2_i) \dtheta- \int_{\partial \Theta(s)} D \nabla c \cdot \bm{n} \right) \ds.  
\end{equation}
For $\mathcal{T}_3$, we proceed similarly, letting $u_{v,s}$ denote the axial component of the velocity field $\bm{u}_v$ (denoted $\bm{u}$ here). We then have 
\begin{equation}
  \label{eq:T3}
    \begin{split}
      \mathcal{T}_3
      &= \int_{\mathcal P} \nabla \cdot (\bm{u}  c ) 
      = \int_{\partial \mathcal P} (\bm{u} c) \cdot \bm{n}  
      = \int_{\Theta(s_2)} u_{v,s} c  - \int_{\Theta(s_1)} u_{v,s} c   + \int_{s_1}^{s_2} \int_{\partial \Theta (s)}(\bm{u}  c ) \cdot \bm{n} \ds \\
      & = \int_{s_1}^{s_2} \partial_s \int_{\Theta(s)} u_{v,s}c  \ds + \int_{s_1}^{s_2} \int_{\partial \Theta(s)}(\bm{u}  c) \cdot \bm{n} \ds \\
      & = \int_{s_1}^{s_2} \left( \partial_s (A \langle u_{v,s}c \rangle ) +  \int_{\partial \Theta(s)}(\bm{u} c) \cdot \bm{n}\right) \ds.
\end{split}
\end{equation}
From \eqref{eq:3d_model_3},  \eqref{eq:3d_model_3_additional} if $R_1 >0$, and the assumption that $\xi$ is constant in each cross-section,  we obtain
\begin{equation}    
  \int_{s_1}^{s_2} \int_{\partial \Theta(s)}  ( (\bm{u} - \bm{w}) c -  D \nabla c)\cdot \bm{n}
  = \int_{s_1}^{s_2}\int_{\partial \Theta_2(s)} \xi (c - c_s)
  = \int_{s_1}^{s_2} \xi  P (\overline{c} -  \overline{c_s}) .
  \label{eq:boundary_condition_derivation}
\end{equation}
 The term $\mathcal{T}_4$ simply reads 
 \begin{equation}
   \mathcal{T}_4 = - \int_{s_1}^{s_2}A\avg{f} \ds.
   \label{eq:T4}
 \end{equation}
Collecting the derivations for $\mathcal{T}_i$, $ i =1,\ldots,4$: \eqref{eq:T1},\eqref{eq:T2},\eqref{eq:T3}, and \eqref{eq:T4} and using \eqref{eq:boundary_condition_derivation} for the resulting boundary terms yield: 
\begin{multline}
\int_{s_1}^{s_2} \left(\partial_t (A\avg{c})- \partial_s (D \partial_{s}(A\avg{c}))  + \xi P (\overline{c}  - \overline{c_s}) + \partial_s (A\avg{u_{v,s} c})\right) \ds \\ 
- \int_{s_1}^{s_2} \sum_{i=1}^2 \frac{\gamma_{i}}{2} \int_0^{2\pi} \left(   \partial_s (D c(R_i)\partial_s R^2_i) \right)  \dtheta \ds =  \int_{s_1}^{s_2} A\avg{f} \ds. \label{eq:derivation_1}
\end{multline}
To make the above equation solvable, we use assumption \ref{assumption:shape_profiles} and write: 
\begin{align*} 
\partial_s(A\avg{c}) &=   A\partial_s \avg{c}+\partial_s A  \avg{c}  = A\partial_s \avg{c} +  \sum_{i=1}^2 -\frac{\gamma_i}{2} \int_0^{2\pi}  \avg{c} \partial_s R_i^2 \dtheta. 
\end{align*}
We use the above and substitute $c = c_v =  w_c(r)\avg{c_v}$ and $\cc =\avg{c_v}$ in \eqref{eq:derivation_1}.
Thus, since \eqref{eq:derivation_1} holds for any $s_1$ and $s_2$, we conclude the result.
\end{proof}

\subsubsection{Boundary conditions for the reduced transport model}

We finalize the derivation of the reduced transport model by stating boundary conditions corresponding to the cross-section average of~\eqref{eq:3d_model_5}, modified from \cite{laurino2019derivation}: 
\begin{equation}
 D_v A \partial_s \cc - A\langle u_{v,s} w_c \rangle \cc = 0 \quad \mathrm{for} \,\, s = 0, L.
  \label{eq:bc_weak_form_1D}
\end{equation}
One can see that integrating \eqref{eq:3d_model_5} over perpendicular cross-sections $\Theta(0)$ and $\Theta(L)$ and using assumption~\eqref{eq:A_s}, the above condition is recovered if $\bm{w}\cdot \bm{n}$ on $\Gamma_0(t) \cup \Gamma_1(t)$ is negligible.

\subsection{Variational formulation of the reduced transport model}
\label{subsec:derivation_1dmodel}
To formally derive a variational formulation of
\eqref{eq:derivation_final} combined with \eqref{eq:bc_weak_form_1D},
we multiply \eqref{eq:derivation_final} by $\phi \in H^1(\Lambda)$ and
integrate by parts. We first observe that
\begin{align*}
  \int_{\Lambda} G(\cc) \phi \ds &\equiv - \int_{\Lambda} \partial_s \left ( D_v \vtwo \cc \right ) \phi \ds 
  =   \int_{\Lambda} \left( D_v \vtwo \cc \right)\partial_s \phi \ds - D_v \vtwo \cc \, \phi \vert_{0}^L  .
\end{align*}
Therefore, after applying the boundary conditions \eqref{eq:A_s} and \eqref{eq:bc_weak_form_1D} and collecting terms, we obtain the variational formulation: for $t > 0$, given coefficients $D_v$, $\xi$ and functions $A \in L^\infty(0,T;L^{\infty} (\Lambda))$, $\avg{f_v} \in L^2(0, T, (H^{1}_{A}(\Lambda))')$, and $u_{v, s}$, $w_c$ such that $\avg{u_{v,s} w_c} \in L^\infty(0,T;L^{\infty} (\Lambda))$ and $\overline{w_c} \in \R $, find $\cc \in L^2(0,T;H_{A}^1(\Lambda))$ with $\partial_t \cc \in L^2(0,T;(H^{1}_{A}(\Lambda))')$ such that for all 
\begin{multline}
  \langle \partial_t \cc, \phi \rangle_{H_A^{-1}(\Lambda)}
  + \int_{\Lambda} D_v \left (A \partial_s \cc + \vtwo \cc \right )\cdot \partial_s \phi
  -  \int_{\Lambda} A \langle u_{v,s} w_c \rangle \cc \cdot \partial_s \phi \\
  + \int_\Lambda  \left(  \xi P \left (\overline{w_c} \cc - \overline{c_s} \right ) + \partial_t A \cc \right) \phi 
  = \langle \avg{f_v}, \phi \rangle_{H_A^{-1}(\Lambda)} , \,\, \forall \phi \in H^1_A(\Lambda). 
  \label{eq:weighted_weak_form_1D_pb}
\end{multline}
As mentioned, note that $\vtwo = 0$ if $w_c = 1$.
In the case of a cylindrical (vessel) domain with $R_2 = R(s, t)$ and $R_1=0$, then $\vtwo = (1-w_c(R)) \partial_s A$. 

\subsection{Variational formulation for the extended transport model}
\label{subsec:weak_form_3d_derivation}

We next formally extend the variational formulation of \eqref{eq:weak_form_3d3d} to the whole domain $\Omega$. Here, a model reduction approach is used, similar to the one by Laurino and Zunino \cite{laurino2019derivation}, to reduce the interface condition \eqref{eq:3d_model_3}. This approach uses the average operator \eqref{eq:bar} as the restriction operator to the centerline for both the trial and test functions. This is different than the approach used in D'Angelo and Quarteroni~\cite{d2008coupling} where the restriction operator for the test functions is taken as  the trace operator onto $\Lambda$ which is well-defined on  special weighted spaces that enjoy better regularity properties than $H^1(\Omega)$. As we will show, the approach  used in \cite{laurino2019derivation} and here is well-defined on functions in $H^1(\Omega)$ and yields to solutions with better regularity properties than the ones in \cite{d2008coupling}. 

From \eqref{eq:weak_form_3d3d}, we have that for $\phi \in H^1_0(\Omega)$
\begin{equation}
  \int_{\Omega_s} \dot{c}_s \phi + \int_{\Omega_s} \nabla \cdot \bm{w} c_s \phi
  + \int_{\Omega_s} D_s \nabla c_s \cdot \nabla \phi + \int_{\Gamma}  \xi (c_s - c_v) \phi
  - \int_{\Omega_s} ( \tilde{\bm{u}}_s c_s) \cdot \nabla \phi
  = \int_{\Omega_s} f_s \phi
  \label{eq:weak_from_cs}
\end{equation}
recalling that $\tilde{\bm{u}}_s = \bm{u}_s - \bm{w}$. For the first two terms, we have that
\begin{equation*}
  \int_{\Omega_s} \dot{c}_s \phi + \int_{\Omega_s} \nabla \cdot \bm{w} c_s \phi
  = \int_{\Omega_s} \partial_t c_s \phi + \int_{\Omega_s} \nabla \cdot (\bm{w} c_s) \phi. 
\end{equation*}
Consider now the fourth term in \eqref{eq:weak_from_cs}. Define an operator subtracting the
perimeter-average i.e. $\tilde{\phi} = \phi -
\overline{\phi}$. Clearly, $(\tilde{\phi}_1, \overline{\phi_2})_{\partial \Theta}=(\tilde{\phi}_2, \overline{\phi_1})_{\partial \Theta}=  0$  since $\overline{\tilde{\phi}_1} = \overline{\tilde{\phi}_2}=0$ for $\phi_1, \phi_2 \in L^1(\partial \Theta)$. We thus have that 
\begin{multline}
  \int_\Gamma \xi (c_s - c_v) \phi
  = \int_\Lambda \int_{\partial \Theta} \xi (c_s - c_v) \phi
  = \int_{\Lambda} \int_{\partial \Theta } \xi (\tilde{c}_s + \overline{c_s} - \tilde{c}_v - \overline{c_v})(\tilde{\phi} + \overline{\phi}) \\ 
  = \int_{\Lambda} \int_{\partial \Theta} \xi (\overline{c_s} - \overline{c_v}) \overline{\phi}
  + \int_{\Lambda} \int_{\partial \Theta} \xi (\tilde{c}_s - \tilde{c}_v)\tilde{\phi}.
  \label{eq:tmp1}
\end{multline}
Following \cite{laurino2019derivation,koppl2018mathematical}, we assume that the second term on the right hand side above is negligible: 
\begin{equation*}
  \int_{\partial \Theta} \xi \tilde{c}_s \tilde{\phi}  \approx 0 , \quad
  \int_{\partial \Theta} \xi \tilde{c}_v \tilde{\phi}  \approx 0.
\end{equation*}
Hence, combining~\eqref{eq:tmp1} with the assumption that $c_v = \avg{c_v} w_c = \cc w_c$ (\Cref{assumption:shape_profiles}), we obtain 
\begin{equation*}
    \int_\Gamma \xi (c_s - c_v) \phi = \int_\Lambda \xi P (\bar{c}_s - \overline{w_c} \cc) \overline{\phi}. 
\end{equation*}

Finally, we identify the domain $\Omega_s$ with $\Omega$ where we introduce the extended solution $c$. That is, we have: 
\begin{equation*}
  \int_{\Omega} \partial_t c \phi
  + \int_{\Omega} \nabla \cdot (\bm{w} c) \phi
  +  \int_{\Omega} \mathcal{E}(D_s) \nabla c \cdot \nabla \phi + \int_\Lambda \xi P (\bar{c} - \overline{w_c} \cc) \overline{\phi}
  - \int_{\Omega} ( \mathcal{E}(\bm{u}_s) - \bm{w}) c \cdot \nabla \phi
  = \int_{\Omega} \mathcal{E} (f_s) \phi.
\end{equation*}
In the above, $\mathcal{E}$ is a suitable extension operator: $\mathcal{E}: H^{1}(\Omega_s) \rightarrow H^{1}(\Omega)$. This operator will be further specified in \Cref{subsec:error_3D}. Integrating the second term above by parts, we arrive at the following weak formulation: Find $c \in L^2(0,T;H^1_0(\Omega))$ with $\partial_t c \in L^2(0,T;H^{-1}(\Omega))$ such that for all $\phi \in H^1_0(\Omega)$,
\begin{equation}
  \int_{\Omega} \partial_t c \phi + \int_{\Omega} \mathcal{E} (D_s) \nabla c \cdot \nabla \phi
  +  \int_\Lambda \xi P (\bar{c} - \overline{w_c} \cc) \bar{\phi} - \int_{\Omega} ( \mathcal E (\bm{u}_s) c)  \cdot \nabla \phi
  = \int_{\Omega} \mathcal{E} (f_s) \phi .
  \label{eq:weak_sol_tissue}
\end{equation}

\subsection{Coupled multi-dimensional variational formulation of transport model} 
\label{sec:coupled_3d_1d}

We now combine the variational formulations derived in Sections~\ref{subsec:derivation_1dmodel}--\ref{subsec:weak_form_3d_derivation}, to summarize the time-dependent coupled 3D-1D solute transport model in variational form. To this end, we introduce the following bilinear forms. First, given $\bm{u}_s$ and for all $c, v \in H^1(\Omega)$, 
\begin{equation*}
  a(c, v) = \int_{\Omega} \mathcal{E}(D_s) \nabla c \cdot \nabla v - \int_{\Omega} (\mathcal{E} (\bm{u}_s) c )\cdot \nabla v ,
\end{equation*}
where $\mathcal{E}$ is an extension operator (to be defined in \Cref{subsec:error_3D}). Second, from inspecting~\eqref{eq:weighted_weak_form_1D_pb}, and recalling the definitions of $\vtwo$ as introduced in~\eqref{eq:def_rho1}, we also define 
for all $\hat{c}, \phi \in H^1(\Lambda)$,  
\begin{align} 
  a_{\Lambda} (\hat{c}, \phi) =  \int_{\Lambda} D_v \left (A \partial_s \cc + \vtwo \cc \right )\cdot \partial_s \phi \ds -  \int_{\Lambda} A \langle u_{v,s} w_c \rangle \cc \cdot \partial_s \phi \ds + \int_\Lambda  \partial_t A \, \cc  \phi \ds . \label{eq:def_A_reduced}
\end{align}
In the above, we recall that $\vtwo$ is given  in  \eqref{eq:def_rho1} and accounts for the deviation of $c_v$ from a uniform distribution in $\Theta(s)$ for $s \in \Lambda$.  

For the coupling terms, we recall the weighted product \eqref{eq:weighted_inner_product} and define for all $v, w \in L^2_P(\Lambda)$:  
\begin{equation}
  b_{\Lambda} (v,w) =   (\xi v,w)_{L^2_P(\Lambda)} .
  \label{eq:def_b_couplingform}
\end{equation}
The coupled weak formulation reads as follows. Given $\mathcal{E}(f) \in L^2(0,T;H^{-1}(\Omega))$ and $\avg{f_v} \in L^2(0,T;H^{-1}_A(\Lambda))$, find $c \in L^2(0,T;H^1_0(\Omega)), \cc \in L^2(0,T;H^1_{A}(\Lambda))$ with $\partial_t c \in L^2(0,T; H^{-1}(\Omega))$, $\partial_t \cc \in L^2(0,T;H^{-1}_A(\Lambda))$ such that 
\begin{subequations}
  \begin{alignat}{2}
    \langle \partial_t c, v \rangle_{H^{-1}(\Omega)}  + a(c,v) +  b_\Lambda (\bar{c} - \overline{w_c}\cc, \bar{v}) &= \langle \mathcal{E}(f),v\rangle_{H^{-1}(\Omega)}, && \,\, \,  \forall v \in H^1_0(\Omega) \label{eq:coupled_pb_1},\\ 
\langle \partial_t \cc ,  \hat{v} \rangle_{H^{-1}_A(\Lambda)} + a_{\Lambda}(\cc,\hat{v}) + b_\Lambda(\overline{w_c}\cc - \bar{c}, \hat{v}) & =  \langle \langle f_v \rangle, \hat{v}\rangle_{H^{-1}_{A}(\Lambda)}, && \,\,\, \forall \hat{v} \in H_A^1(\Lambda), \label{eq:coupled_pb_2} \\ 
c^0 = \mathcal{E} (c_s^0) \in L^2(\Omega) & , \quad \hat{c}^0  = \avg{c_v^0} \in L^2_A(\Lambda). && \label{eq:coupled_pb_3}
  \end{alignat}
  \label{eq:coupled_3d_1d_weak}
\end{subequations}
Observe that the term $b_\Lambda(\bar{c},\bar{v})$ is well-defined since for $v \in H^1(\Omega)$, $\bar{v} \in L_{P}^2(\Lambda)$. Indeed, by Jensen's and trace inequality \eqref{eq:trace_constant_indep_t}, we have that  
\begin{align}
  \|\bar{v}\|_{L_P^2(\Lambda)}^2
  = \int_\Lambda \frac{1}{P}\left(\int_{\partial \Theta} v \right)^2
  \leq \int_{\Lambda} \int_{\partial \Theta } v^2
  =  \|v\|_{L^2(\Gamma)}^2 \leq C_1^2 \|v\|_{H^1(\Omega)}^2.
  \label{eq:jensen_ineq}
 \end{align}
\begin{prop}[Well-posedness and regularity of the 3D-1D problem] \label{prop:well_posedness_regularity}
  Assume that $A, \partial_t A, \\ \avg{u_{v,s} w_c}  \in  L^{\infty}(0,T;L^{\infty}(\Lambda))$, $A \geq A_0 > 0$ a.e in $\Lambda$,  $ \mathcal{E} \bm{u}_s \in L^{\infty}(0,T;L^{\infty}(\Omega, \R^{d\times d}))$, and  that $ \mathcal{E}(\mathcal{D}_s) \in L^{\infty}(0,T;L^{\infty}(\Omega))$ with uniform ellipticity constant $\tilde{\nu} > 0$. Then, the coupled weak formulation \eqref{eq:coupled_3d_1d_weak} is well-posed. 
  
  In addition, if the material parameters are H\"{o}lder continuous of index $\beta > 1/2$: 
  \[\|\mathcal{E}(D_s)(t_1) -\mathcal{E} (D_s)(t_2)\|_{L^{\infty}(\Omega, \R^{d \times d})} +  \|D_v(t_1) - D_v(t_2)\|_{L^{\infty}(\Lambda)} + \|\xi(t_1) - \xi (t_2)\|_{L^{\infty}(\Lambda)} \leq C |t_2 - t_1|^{\beta},  \]
for some constant $C$ independent of $t$, and if  
  $\partial \Omega \in C^2$, $c_v^0 \in H^1(\Omega_v),$  $c_s^{0} \in H^1(\Omega_s)$,  $\mathcal{E}(f) \in L^2(\Omega)$ and $\avg{f} \in L^2_{A}(\Lambda)$,  then $$c \in L^2(0,T;H^{3/2-\eta}(\Omega)), \quad  \eta > 0.$$
\end{prop}
\begin{proof}
  \textit{Well-posedness}. We use J.-L. Lions theorem, see e.g \cite[Theorem 10.9]{brezis2010functional}. Let $\bm{V} = H^1_0(\Omega) \times H^1_{A}(\Lambda)$ with dual $\bm{V}' = H^{-1}(\Omega) \times H^{-1}_{A}(\Lambda)$. The space $\bm{V}$ defines a Hilbert space with inner product $(\bm{u}, \bm{v})_{\bm{V}} = (u,v)_{H^1_0(\Omega)} + (\hat{u},\hat{v})_{H^1_A(\Lambda)}$, for all $\bm{u} =(u,\hat{u})$ and $\bm{v} = (v,\hat{v}) \in \bm{V}$. Further it holds that $\bm{V} \subset L^2(\Omega) \times L^2_{A} (\Lambda) \subset \bm{V}'$.  We then write \eqref{eq:coupled_3d_1d_weak} as: Find $\bm{c} = (c , \hat{c}) \in \mathcal{W}(\bm{V},\bm{V}')$ such that 
  \begin{equation*}
    \langle \partial_t \bm{c},\bm{v} \rangle_{V' \times V} + \mathcal{A}_{\Lambda} (t, \bm{c}, \bm{v}) = \ell(\bm{v}), \quad \forall \bm{v} \in \bm{V}, 
   \end{equation*}
   where for all $\bm{c} = (c, \hat{c}), \bm{v} =(v, \hat{v}) \in \bm{V}$  
   \begin{align*}
    \mathcal{A}_{\Lambda} (t, \bm{c},\bm{v}) & =  a(c,v) +  b_\Lambda (\bar{c} - \overline{w_c}\cc , \bar{v}) +  a_{\Lambda}(\cc,\hat{v}) + b_\Lambda(\overline{w_c}\cc - \bar{c}, \hat{v}), \\ 
 \ell(\bm{v}) & = (\mathcal{E}f,v)_\Omega + (A\langle f_v \rangle, \hat{v})_{\Lambda}. 
   \end{align*}
We proceed to show that the continuity and coercivity conditions of Lions' Theorem hold: 
There exist constants $M$, $\kappa$ and $\mu$ independent of $t$ such that 
\begin{alignat}{2}
  \mathcal{A}_{\Lambda} (t, \bm{c},\bm{v}) & \leq M \|\bm{c}\|_{\bm{V}} \|\bm{v}\|_{\bm{V}}, && \quad \forall \bm{c}, \bm{v} \in \bm{V}.  \label{eq:continuity_A} \\ 
  \mathcal{A}_{\Lambda}(\bm{c},\bm{c}) &\geq \kappa \|\bm{c}\|_{\bm V}^2 - \mu(\|c\|^2 + \|\hat{c}\|_{L^2_A(\Lambda)}^2),&& \quad \forall \bm{c} \in \bm{V}. \label{eq:coercivity_A}
\end{alignat}
We begin by showing \eqref{eq:continuity_A}. 
By H\"{o}lder's inequality, we immediately have that 
\[ a(c,v) \leq (\|\mathcal{E}(D_s)\|_{L^{\infty}(\Omega)}+ \|\mathcal{E} (\bm{u}_s)\|_{L^{\infty}(\Omega)}) \|c\|_{H^1(\Omega)} \|v\|_{H^1(\Omega)}. \]
Further, with H\"{o}lder's and triangle inequalities  and \eqref{eq:jensen_ineq}, we have that 
  \begin{align*}
  & b_\Lambda (\bar{c} - \overline{w_c}\cc , \bar{v})  +  b_\Lambda(\overline{w_c}\cc - \bar{c}, \hat{v}) \leq \|\xi\|_{L^{\infty}(\Lambda)} (\|\bar{c}\|_{L^2_{P}(\Lambda)} + \| \overline{w_c}\cc \|_{L^{2}_P(\Lambda)}) (\|\bar{v}\|_{L^2_P(\Lambda)} + \|\hat{v}\|_{L^2_P(\Lambda)})\\ 
& \leq \|\xi\|_{L^{\infty}(\Lambda)}\left(C_1 + (\|\overline{w_c}\|_{L^{\infty}(\Lambda)}+ 1) \|PA^{-1}\|_{L^{\infty}(0,T;L^{\infty}(\Lambda))}^{1/2}\right)^2 \|\bm{c}\|_{\bm{V}} \|\bm{v}\|_{\bm{V}}. 
\end{align*}
In the above, we note that $C_1$ is independent of $t$, see \eqref{eq:trace_constant_indep_t},  and we use the definition of weighted norms which result in the following bound. 
\begin{align}
\|P^{-1}A\|^{-1}_{L^{\infty}(\Lambda)}\|\hat{c}\|_{L^2_A(\Lambda)}^2 \leq  \|\hat{c}\|_{L^2_{P}(\Lambda)}^2 \leq \|PA^{-1}\|_{L^\infty(0,T;L^{\infty}(\Lambda))} \|\hat{c}\|_{L^{2}_A(\Lambda)}^2. \label{eq:Lp_La}
\end{align}
With \eqref{eq:Lp_La} and H\"{o}lder's inequality, the following easily follows. 
\begin{multline*}
  a_{\Lambda}(\hat{c},\hat{v}) \leq (\|D_v\|_{L^{\infty}(\Lambda)} +  \|g_s A^{-1}\|_{L^{\infty}(\Lambda)}+ \|A^{-1}\partial_t A\|_{L^{\infty}(\Lambda)}+ \|\langle u_{v,s}w_c\rangle\|_{L^{\infty}(\Lambda)}) \|\hat{c}\|_{H^1_A(\Lambda)} \| \hat{v}\|_{H^1_A(\Lambda)}. 
 \end{multline*}
 By combining the above bounds, we obtain \eqref{eq:continuity_A} for a constant $M$ independent of $t$.  We now show \eqref{eq:coercivity_A}, but we do not track constants for simplicity. It easily follows that 
\begin{equation*}
  a(v,v) + b_{\Lambda}(\overline{c}, \overline{c}) \geq \frac{\tilde{\nu}}{2} \|\nabla c\|_{L^2(\Omega)}^2 - \frac{1}{2\tilde{\nu}}\|\mathcal{E} \bm{u}_s\|_{L^\infty(\Omega)}^2 \|c\|_{L^2(\Omega)}^2 + \|\xi^{1/2}\overline{c}\|_{L^2_{P}(\Lambda)}^2 . 
\end{equation*}
With similar arguments and with \eqref{eq:Lp_La}, we also have positive constants $\kappa_1$ and  $\mu_1$ such that 
\begin{equation*}
a_{\Lambda} (\hat{c},\hat{c})+ b_\Lambda(\overline{w_c} \hat{c}, \hat{c}) \geq \kappa_1 \|\hat{c}\|_{H^1_{A}(\Lambda)}^2 - \mu_1 \| \hat{c}\|^2_{L^2_A(\Lambda)} + \| \xi^{1/2} \overline{w_c}^{1/2} \hat{c}\|^2_{L_{P}^2(\Lambda)} . 
\end{equation*} 
To handle the coupling terms, we use  
Young's inequality and  \eqref{eq:Lp_La} as follows.
\begin{multline*}
    |b_{\Lambda}(\overline{w_c} \hat{c}, \overline{c}) | + |b_{\Lambda} (\overline{c}, \hat{c}) | \leq (\|\xi^{1/2}\overline{w_c} \hat{c}\|_{L^2_P(\Lambda)}  + \|\xi^{1/2}\hat{c}\|_{L^2_P(\Lambda)})\|\xi^{1/2}\overline{c}\|_{L^2_P(\Lambda)}\\  \leq  \frac12 \|\xi^{1/2} \overline{c}\|_{L^2_P(\Lambda)}^2 + \frac12(\|\xi^{1/2} \overline{w_c}\|_{L^{\infty}(\Lambda)}^2 + \|\xi^{1/2}\|_{L^{\infty}(\Lambda)}^2)\|PA^{-1}\|_{L^{\infty}(\Lambda)}\|\hat{c}\|_{L^2_A(\Lambda)}^2.
\end{multline*}
Then, upon writing $\mathcal{A}_{\Lambda}(t, \bm{c},\bm{c}) - b_{\Lambda}(\overline{w_c} \hat{c}, \overline{c})  - b_{\Lambda} (\overline{c}, \hat{c}) = a(c,c) + b_{\Lambda} (\overline{c},\overline{c}) +  a_{\Lambda} (\hat{c},\hat{c}) + b_{\Lambda}(\overline{w_c}\hat{c}, \hat{c})$ and using the above bounds, we conclude that \eqref{eq:coercivity_A} holds. In addition, one easily sees that $\ell$ defines a bounded functional on $\bm{V}$. Therefore, all the requirements for Lions Theorem hold and existence and uniqueness of weak solutions is obtained. 

\textit{Additional regularity}. We proceed to show the stated $H^{3/2- \eta}$ regularity. The first step is to show that $\partial_t c \in L^2(0,T;L^2(\Omega))$. This is achieved by invoking  maximal regularity \cite[Theorem 7.1]{arendt2017jl}. We verify  that $\mathcal{A}_{\Lambda}(t,\bm{c},\bm{v})$ is H\"{o}lder continuous of index $\beta >1/2$: there exists a  constant $C$  independent of $t$ such that 
\begin{equation}
  |\mathcal{A}_{\Lambda}(t_2, \bm{c}, \bm{v}) - \mathcal{A}_{\Lambda}(t_1, \bm{c}, \bm{v}) | \leq K|t_2 - t_1|^{\beta} \|\bm{c}\|_{\bm{V}} \|\bm{v}\|_{\bm{V}}, \quad \forall \bm{c}, \bm{v} \in \bm{V}. \label{eq:Holder_cont_bound}
\end{equation}  
The delicate terms in $\mathcal{A}_{\Lambda}$ are the ones involving $b_{\Lambda}(\cdot,\cdot)$ as the bounds for all the other terms follow directly from the assumptions on the material parameters. We provide some details for showing \eqref{eq:Holder_cont_bound} in Appendix~\ref{appendix:handling_b}. Under the additional assumption that $c_s^0 \in H^1(\Omega_s), c_v^0 \in H^1(\Omega_v)$, we have that $c^0 = \mathcal{E} c_s \in H^1(\Omega)$ and $\hat{c}^0 = \avg{c_v^0} \in H^1_{A}(\Lambda)$. 
Thus, since $\mathcal{E} f \in L^2(\Omega)$ and $\avg{f_v} \in L_A^2(\Lambda)$, we have verified the assumptions of   \cite[Theorem 7.1]{arendt2017jl} and $\partial_t c \in L^2(0,T;L^2(\Omega))$.  

We now use the fractional space $H^{1/2 + \eta}(\Omega)$ normed by
\[\|v\|_{H^{1/2+\eta}(\Omega)}^2 =  \|v\|^2_{L^2(\Omega)} + \int_\Omega \int_\Omega \frac{|v(x) - v(y)|^2}{|x-y|^{4+2\eta}},  \] and we define the linear functional $\mathcal{F}(v)$: 
\begin{equation}
  \mathcal{F}(v)
  = - \int_{\Omega} \partial_t c v + \int_\Gamma \xi (\overline{w_c}\cc - \bar{c}) v
  + \int_{\Omega} \nabla \cdot(\mathcal{E} \bm{u}_s c) v + \int_\Omega \mathcal{E} f v.
  \label{eq:def_F_v}
\end{equation}
The trace theorem yields for a positive constant $K_{\Gamma}$ \cite{di2012hitchhikers}: 
\begin{equation}
  \|v\|_{L^2(\Gamma)} \leq K_{\Gamma} \|v\|_{H^{1/2+\eta}(\Omega)}, \quad \forall v \in H^{1/2+\eta}(\Omega).
  \label{eq:trace_theorem}
\end{equation}
With the above, $\mathcal{F}(v)$ is a bounded linear functional on $H^{1/2+\eta}(\Omega)$. Indeed, with Cauchy-Schwarz inequality and \eqref{eq:trace_theorem},   we have 
\begin{multline}
  \label{eq:bounding_negative_norm_F_0}
  \sup_{v \in H^{1/2+\eta}(\Omega), v \neq 0} \frac{|\mathcal{F}(v)|}{\|v\|_{H^{1/2+\eta}(\Omega)}}  \leq  \|\partial_t c\|_{L^2(\Omega)}
  +  K_\Gamma \|\xi (\overline{w_c}\cc  - \bar{c})\|_{L^2(\Gamma)} \\
  + \|\mathcal{E}\bm{u}_s\|_{L^\infty(\Omega)} \|\nabla c\|_{L^2(\Omega)} 
  + \|\nabla (\mathcal{E}\bm{u}_s)\|_{L^6(\Omega)} \|c\|_{L^3(\Omega)} + \|\mathcal{E} f\|_{L^2(\Omega)} .
\end{multline}
The second term above is further bounded as follows: 
\begin{equation*}
  \|\xi (\overline{w}_c \cc  - \overline{c})\|_{L^2(\Gamma)} \leq \|\xi  \overline{w}_c \hat{c}\|_{L^2_{P}(\Lambda)} + \|\xi  \overline{c}\|_{L^2_P(\Lambda)} \leq K ( \|\hat{c}\|_{L^2_{P}(\Lambda)}  + C_1 \|c\|_{H^1(\Omega)}), 
\end{equation*} 
where we used \eqref{eq:jensen_ineq} and \eqref{eq:trace_constant_indep_t}. For the third and fourth terms in \eqref{eq:bounding_negative_norm_F_0}, Sobolev embedding results yield: 
\begin{equation*}
  \|\mathcal{E}\bm{u}_s\|_{L^{\infty}(\Omega)}\|\nabla c\|_{L^2(\Omega)} + \|\nabla (\mathcal{E}\bm{u}_s)\|_{L^6(\Omega)} \|c\|_{L^3(\Omega)} \leq K \|\mathcal{E}\bm{u}_s\|_{H^2(\Omega)}\|c\|_{H^1(\Omega)}.
\end{equation*}
Thus, \eqref{eq:bounding_negative_norm_F_0} becomes:
\begin{multline}
  \|\mathcal{F}\|_{H^{-1/2-\eta}(\Omega)} \leq \|\partial_t c\|_{L^2(\Omega)} + K\|\hat{c}\|_{L^2_{P}(\Lambda)} \\
  + K(C_1 + \|\mathcal{E}\bm{u}_s\|_{H^2(\Omega)}) \|c\|_{H^1(\Omega)} + \|\mathcal{E}f\|_{L^2(\Omega)}.
\label{eq:bounding_negative_norm_F}
\end{multline}
For a.e. $t \in (0,T)$,  $c(t) \in H^1_0(\Omega)$ solves 
\begin{equation}
  \int_{\Omega} D_s \nabla c \cdot \nabla v = \mathcal{F}(v), \quad \forall v \in H^1_0(\Omega).
\end{equation}
It then follows from Lemma 3.10 in \cite{laurino2019derivation}, see also \cite{gong2014approximations}, and the principle of superposition that a.e.~in time 
\begin{equation}
\|c\|_{H^{3/2-\eta}(\Omega)} \leq K \|\mathcal{F}\|_{H^{-1/2-\eta}(\Omega)}. \label{eq:h32_estimate}
\end{equation}
The above can also be deduced from interpolation theory, see Chapter 14 in  \cite{brenner2008mathematical}. 
Integrating the above bound over $t$,  using \eqref{eq:bounding_negative_norm_F} and the regularity properties of $c$ and $\cc$ as discussed above, we have  that $c \in L^2(0,T; H^{3/2 - \eta}(\Omega))$. 
\end{proof}

\subsection{Extension to vascular networks} 
\label{sec:networks}

\newcommand{\uvsi}{\avg{u_{v, s, i}}}
\newcommand{\uvsk}{\avg{u_{v, s, k}}}


Up til now, we have considered a representation of a single vessel and
its surroundings. However, in applications such as for transport in
the human (peri)-vasculature or root networks, each vessel is but a
segment of a larger (peri)vascular network. To extend our setting,
consider now a network of $N$ domains $\Omega_{v,i}$ with
center-curves $\Lambda_{i} = \{\bm{\lambda}_i (s), \,\,s \in
(0,L_i)\}$ for $ i = 1,\ldots, N$. Denote by $\Lambda_{\mathrm{graph}}
= \cup_{i} \Lambda_i$. We use a similar notation and approach as
Laurino and Zunino~\cite{laurino2019derivation}. By direct extension from one to several vessels, letting $w_{c} = 1$
in $\Omega_{v, i}$ for all $i$ for clarity, we have that in each $\Lambda_i$, the 1D
concentration $\cc_i$ solves:
\begin{equation}
  \partial_t (A_i \cc_i) - \partial_{s}( D_{v,i} A_i \partial_s \cc_i)  + \partial_s (A_i \uvsi \cc_i) + \xi P_i ( \cc_i  - \overline{c_s})   = A_i\avg{f_{v,i}}.
    \end{equation}

The key next step is to specify interface and inlet/outlet
conditions. Let $Y$ denote the collection of bifurcation points,
i.e.~vertices that are shared between two or more curves: $y \in Y$ if
there exists at least one pair $(i,j)$ such that $y = \lambda_i(0) =
\lambda_j(L_j)$. The set of curves with inlet nodes is denoted by $I$
and the set of curves with outlet nodes is denoted by $O$.  At the
level of one node $y_j \in Y$, we separate the connecting curves as
follows.
\begin{align*}
  I_j &= \{ i \in \{1, \ldots, N\}: \,\, \bm{\lambda}_i(0) = y_j\} \quad \text{(curves having $y_j$ as an inlet node)}, \\
  O_j &= \{ i \in \{1, \ldots, N\}: \,\, \bm{\lambda}_i(L_i)= y_j\} \quad \text{(curves having $y_j$ as an outlet node)}.
\end{align*}
Now, at every bifurcation point, we enforce conservation of fluxes
and continuity or instantaneous mixing of the solute:
\begin{align*}
\sum_{k \in I_j} (A_k \uvsk \cc_k - D_{v,k}  \partial_{s} \cc_k)(0) & = \sum_{k \in O_j} (A_k \uvsk \cc_k - D_{v,k}  \partial_{s} \cc_k)(L_k), \text{ and }   \\ 
\cc_k(0) & = \cc_{i}(L_i), \quad \forall k \in I_j, \,\, \forall i \in O_j. 
\end{align*}
For inlet and outlet curves, we set 
\begin{align*}
   (A_k \uvsk \cc_k -  D_{v,k} \partial_{s} \cc_k) (0) &= 0 , \quad  \foralls k \in I \\ 
     (A_k \uvsk \cc_k -  D_{v,k} \partial_{s} \cc_k)(L_k) &= 0 , \quad  \foralls k \in O.
\end{align*} 

Let
\begin{equation}
  H^1(\Lambda_{\mathrm{graph}}) = \bigoplus H^1_{A_i}(\Lambda_i) \cap C^0(\Lambda_{\mathrm{{graph}}}), 
\end{equation}
consist of functions that are locally in $H_{A_i}^1(\Lambda_i)$ for each $i$ (which implies continuity in each $\Lambda_i$ since $\Lambda_i$ is 1D) and that are continuous across bifurcation points. A natural weak formulation for the coupled network with the 3D surroundings now follows: Given $\mathcal{E}f \in L^2(0,T;H^{-1}(\Omega))$ and $(\avg{f_{v,1}}, \ldots, \avg{f_{v,n}}) \in L^2(0,T;H^1(\Lambda_{\mathrm{graph}})')$, find $c \in L^2(0,T;H^1_0(\Omega)), \cc = (\hat{c}_1 , \ldots, \hat{c}_N)  \in L^2(0,T;H^1(\Lambda_{\mathrm{graph}}))$ with $\partial_t c \in L^2(0,T; H^{-1}(\Omega)), \partial_t \cc \in L^2(0,T;H^1(\Lambda_{\mathrm{graph}})')$ such that for all $v \in H^1_0(\Omega)$ and $\hat{v} \in H^1(\Lambda_{\mathrm{graph}})$:
\begin{subequations}
  \begin{align}
    \langle \partial_t c, v \rangle_{H^{-1}(\Omega)}  + a(c,v) + \sum_{i=1}^N b_{\Lambda_i} (\bar{c} - \cc_i, \bar{v})
    &= \langle \mathcal{E}f,v \rangle_{H^{-1}(\Omega)}  \label{eq:coupled_pb_1_graph} \\
    \sum_{i=1}^N \left(  \langle \partial_t \cc_i, \hat{v} \rangle_{H^{-1}_{A_i}(\Lambda_i)} +  a_{\Lambda_i}(\cc_i,\hat{v}) +  b_{\Lambda_i} (\cc_i - \bar{c}, \hat{v})\right)
    &= \sum_{i=1}^N \langle \langle f_{v,i} \rangle, \hat{v} \rangle_{H^{-1}_{A_i}(\Lambda_i)} 
    \label{eq:coupled_pb_2_graph}
  \end{align}
  \label{eq:coupled_pb_3_graph}%
\end{subequations}%
with the initial conditions
\begin{equation}
  c^0 = \mathcal{E} c_s^0, \quad \hat{c}_i^0  = \avg{c_{v,i}^0} \quad i \in \{1, \ldots, N \}.
\end{equation}
In the above, the forms $b_{\Lambda_i}(\cdot,\cdot)$ and $a_{\Lambda_i}(\cdot, \cdot)$ are obtained by naturally modifying the form $b(\cdot,\cdot)$  \eqref{eq:def_b_couplingform} and  the form $a_{\Lambda}(\cdot, \cdot)$ \eqref{eq:def_A_reduced} respectively.  At the cost of only additional notation and conditions similar to \eqref{eq:A_s}, the above can be easily extended for the case of non-uniform concentration profiles $w_c \neq 1$ in each vessel.

\section{Coupled 3D-1D-1D models of solute transport} 
\label{sec:3D1D1D}

In this section, we focus on a case of particular neurological
relevance, namely the case of a vascular network surrounded by a
perivascular network and embedded in brain tissue with semi-permeable
and moving membranes.  
From the 3D-3D-3D
equations, we derive a coupled 3D-1D-1D model formulation allowing for
strong jumps between the vascular and tissue domains in terms of
material parameters (e.g.~diffusion coefficient, velocity). We do not
analyze this model further here, beyond stating a weak
formulation. However, noting the similarity between the 3D-1D and
3D-1D-1D models, we expect that their well-posedness and model error
analysis would follow from applications of the same techniques.

\subsection{A coupled 3D-3D-3D model of vascular-perivascular-tissue transport}
\label{sec:3d3d3d}

We now consider the case of $\Omega_v(t)$ representing a cylindrical
blood vessel
\begin{equation*}
  \Omega_v(t) = \{ \bm{\lambda}(s) + r\cos(\theta) \bm{N}(s) + r\sin(\theta) \bm{B}(s), 0 < s < L, \, 0 \leq \theta \leq 2 \pi, \, 0 \leq r < R_1(s,t,\theta ) \},
\end{equation*}
and introduce an intermediate annular domain $\Omega_p(t)$ representing
a perivascular space along the centerline $\bm{\lambda}$ surrounding
the blood vessel $\Omega_v(t)$:
\begin{multline*}
  \Omega_p(t)
  = \{ \bm{\lambda}(s) + r\cos(\theta) \bm{N}(s) + r\sin(\theta) \bm{B}(s), \\ 0  < s  <  L, \, 0 \leq \theta \leq 2 \pi, \,  R_1(s,t,\theta) < r < R_2(s,t,\theta ) \} .
\end{multline*} 
The domain $\Omega_p(t)$ is further surrounded by a domain
$\Omega_s(t) \subset \R^d$, and the fixed domain $\Omega$ is defined
such that $\Omega =  (\Omega_p \cup \Omega_v \cup \Omega_s)$. In
each domain $\Omega_i$, for $i \in \{v,p,s\}$ and $t \in (0, T]$, we
  assume that we are given a velocity field $\bm{u}_i$ and diffusion
  coefficient $D_i$, and we are interested in finding the 
  concentration $c_i : \Omega_i \times (0, T] \rightarrow \R$ such that
\begin{equation*}
  \partial_t c_i - \nabla \cdot (D_i \nabla c_i ) + \nabla \cdot (\tilde{\bm{u}}_i c_i) = f_i.
\end{equation*}
As before, $\tilde{\bm{u}}_i = \bm{u}_i - \bm{w}$ with $\bm{w}$ representing the domain velocity.

We assume that the interfaces $\Gamma_v$ (separating the vasculature
$\Omega_v$ and perivasculature $\Omega_p$) and $\Gamma_s$ (separating
the perivasculature $\Omega_p$ and tissue $\Omega_s$) are
semi-permeable:
\begin{align*}
    (c_v \tilde{\bm{u}}_v - D_v \nabla c_v) \cdot \bm{n} -\xi_{v}(c_v - c_p ) &= 0  &&  \mathrm{on} \, \Gamma_{v}   \times (0,T],\\ 
    (c_p \tilde{\bm{u}}_p - D_p \nabla c_p) \cdot \bm{n} -\xi_{s} (c_p - c_s) &= 0  &&  \mathrm{on} \, \Gamma_{s} \times (0,T],
\end{align*}
where $\bm{n}$ denotes a consistently-oriented normal at the interfaces, and $\xi_v$ and $\xi_s$ are the membrane permeabilities of the concentrations. In addition, conservation of mass is enforced on $\Gamma_{v}$ and $\Gamma_{s}$ with conditions similar to~\eqref{eq:3d_model_4}. At the sides $\Gamma_v^0, \Gamma_s^0$ and $\Gamma_v^L, \Gamma_s^L$, where by the superscripts $0$ and $L$ we denote the cross-sections of any
interface $\Gamma$ at $s = 0$ and $s = L$, respectively, we apply no flux boundary conditions:
\begin{align*}
  (c_i \tilde {\bm{u}}_i + \nabla c_i) \cdot \bm{n} = 0, \quad \mathrm{on} \,\, \Gamma_i^0 \cup \Gamma_i^{L} \times (0,T], \,\,\, i \in \{v, s\} .
\end{align*}
On $\partial \Omega \times (0,T]$, we set $c_s = 0$. 

\subsection{Derivation of 1D averaged equations}

We now aim to derive coupled cross-section averaged equations for the  vascular and perivascular concentrations. Let $\Theta_v(s)$ and $\Theta_p(s)$ be the cross-sections of $\Omega_v$ and $\Omega_p$ respectively at $ s \in \Lambda$. We also denote by $\partial \Theta_v(s)$ and $\partial \Theta_p(s)$ the inner and outer boundaries of the cross-section $\Theta_p(s)$, respectively. Note that $\partial \Theta_v(s)$ is also the boundary of $\Theta_v(s)$. We let $A_i(s,t) = |\Theta_i(s,t)|$ and $P_i = |\partial \Theta_i(s,t)|$ for $i \in \{v, p\}$. We introduce the following cross-sectionally averaged quantities 
\begin{alignat*}{2}
  \cc_v(s,t)  & = \frac{1}{A_{v}(s,t)} \int_{\Theta_v(s,t)} c_v,&& \quad 
  \cc_p(s,t)  =  \frac{1}{A_{p}(s,t)} \int_{\Theta_p(s,t)} c_p,   \quad \forall (s,t) \in \Lambda \times (0,T).   
\end{alignat*}

The reduced 1-D equations for $\cc_v$ and $\cc_p$ are presented in the next proposition. To clarify the presentation, we consider the constant cross-section case (rather than allowing radially-varying weights $w_c$).
\begin{prop}[1D-1D vascular-perivascular transport equations]
  Assume that the vascular and perivascular concentrations $c_v, c_p$ solve the equations of~\Cref{sec:3d3d3d}, are constant on each cross-section:
  \begin{equation*}
    c_v(s,r,\theta, t) = \cc_v, \quad 
    c_p(s,r,\theta, t) = \cc_p, 
  \end{equation*}
  and are sufficiently regular in the sense that
  $c_v \in L^1(\Theta_v(s)) \cap L^{1}(\partial \Theta_v(s))$,  $c_p \in  L^1(\Theta_p(s)) \cap L^{1}(\partial \Theta_p(s))$ for all $s \in \Lambda$. Also assume that $c_s \in  L^{1}(\partial \Theta_p(s))$ for all $s \in \Lambda$. Then, the vascular cross-section averaged concentration $\cc_v$ satisfies the following in $\Lambda$:  
  \begin{equation}
    \label{eq:3d_1d1d_1}
       \partial_t (A_v \cc_v)  -  \partial_{s}(D_v A_v \partial_s \cc_v) + \partial_s (A_v \avg{u_{v,s}} \cc_v) + \xi_{v}P_v (\cc_v - \cc_p ) = A_v \avg{f_v}.  
  \end{equation}
  In addition, the perivascular cross-section averaged concentration $\cc_p$ satisfies the following also in $\Lambda$:
  \begin{equation}
    \label{eq:3d_1d1d_2}
    \partial_t (A_p \cc_p) - \partial_{s} (D_p A_p \partial_s \cc_p)
    + \partial_s (A_p \avg{u_{p,s}} \cc_p) + \xi_{v} P_v (\cc_p - \cc_v) + \xi_{s} P_p (\cc_p - \overline{c_s} ) = A_p \avg{f_p},
  \end{equation}
  where $\overline{c_s}$ is the lateral average of $c_s$ over $\partial \Theta_p$.  
\end{prop}
\begin{proof} We provide a brief proof sketch. For deriving \eqref{eq:3d_1d1d_1}, we follow the same arguments as the proof of~\Cref{prop:derivation_1d}. In particular, with the notation of \Cref{prop:derivation_1d}, the same equations hold with $R_1 = 0$, $R_2 = R_1$, and $\overline{c_s} = \overline{c_p} = \cc_p$. For \eqref{eq:3d_1d1d_2}, the same arguments also hold. The main difference is in the step~\eqref{eq:boundary_condition_derivation}. We now have by the stated interface and boundary conditions that
  \begin{multline*}  
  \int_{\partial \Theta_v(s)} ( \tilde{\bm{u}}_p c_p -  D_p \nabla c_p ) \cdot \bm{n}_p + \int_{\partial \Theta_p(s)} ( \tilde{\bm{u}}_p c_p -  D_p \nabla c_p )\cdot \bm{n}_p  \\ 
 = \int_{\partial \Theta_v(s)} \xi_{v} (c_p - c_v) +  \int_{\partial \Theta_p(s)} \xi_{s} (c_p - c_s) =  \xi_{v} P_v (\overline{c_p} - \overline{c_v}) + \xi_{s} P_s(\overline{c_p} - \overline{c_s}),
  \end{multline*}
  where the overlines denote context-dependent lateral averages (defined relative to the respective interfaces). Now, invoking the cross-section average assumptions, we adopt all the remaining arguments in the proof of \Cref{prop:derivation_1d} to arrive at the stated equations.
\end{proof}

\subsection{Coupled 3D-1D-1D formulation}

A similar approach as in~\Cref{subsec:weak_form_3d_derivation} is adopted to extend the solution $c_s$ to the whole domain $\Omega$. The coupled 3D-1D-1D perivascular-vascular-tissue weak formulation then reads: find $c \in L^2(0,T;H_0^1(\Omega))$ with $\partial_t c \in L^2(0,T;H^{-1}(\Omega))$ such that 
\begin{equation}
  \langle \partial_t c, v\rangle_{H^{-1}(\Omega)} + a(c, v) +  b_{\Lambda}^p (\xi_{s}(\overline{c} - \cc_p), \overline{v})
  = (\mathcal{E}f,v)_{\Omega}, \quad \foralls v \in H_0^1(\Omega).
  \label{eq:coupled_3d3d1d_1}
\end{equation}
In addition, find $\cc_i \in L^2(0,T;H^1_{A_i}(\Lambda))$ with $\partial_t \cc_i \in L^2(0,T;H^{-1}_{A_i}(\Lambda))$ for $i \in \{v, p\}$ such that $\foralls \hat{v} \in H^1_{A_p}(\Lambda)$, 
\begin{align}
  \label{eq:coupled_3d3d1d_2} 
  \langle \partial_t \cc_p, \hat{v} \rangle_{H^{-1}_{A_p}(\Lambda)}
  +  a_{\Lambda}^p(\cc_p, \hat{v})  +b_{\Lambda}^p ( \xi_{s} (\cc_p - \overline{c_s}), \hat{v}) 
  +  b_{\Lambda}^v(\xi_v(\cc_p - \cc_v), \hat{v}) &= (A_p \avg{f_p},\hat{v})_{\Lambda},
\end{align}
and $\foralls \hat{v} \in H^1_{A_v}(\Lambda)$:
\begin{equation}
\label{eq:coupled_3d3d1d_3}
  \langle \partial_t \cc_v, \hat{v}\rangle_{H^{-1}_{A_v}(\Lambda)} + a_{\Lambda}^v(\cc_v, \hat{v})  +  b_{\Lambda}^v(\xi_{v}(\cc_v -\cc_p) , \hat{v})
  = (A_p \avg{f_v},\hat{v})_{\Lambda} .
\end{equation}
In the above, the forms $a^p_{\Lambda}, a^v_{\Lambda}$ are given by~\eqref{eq:def_A_reduced} where $A$ is taken to be either $A_v$ or $A_p$, and we have defined
\begin{equation*}
  b_{\Lambda}^i( \hat{v},\hat{w}) = ( \hat{v},\hat{w})_{\Lambda,P_i}, \quad \forall \hat{v},\hat{w} \in L^2_{P_i}(\Lambda), \,\,  i \in \{p,v\} .  
\end{equation*}

\section{Inequalities for Sobolev spaces over annular and moving domains}
\label{sec:inequalities}

To estimate the modelling error induced by the model reduction
introduced in~\Cref{sec:3D1D}, we expect to rely on typically standard
inequalities such as the Poincar\'e and trace inequalities on
$H^1(\Omega_v)$. However, since generally $\Omega_v$ is non-convex and
allowed to move in time, these inequalities require some
attention. Moreover, a key question is how the inequality constants
depend on the (inner and) outer radii. In this section, we address
these theoretical questions separately. Here and in what follows, we assume that $R_1$ and $R_2$ are independent of $\theta$. 

We define the maximal cross-section diameter
$\epsilon_{\max}$ and axial radius variation $\epsilon_s$: 
\begin{align}
  \epsilon_{\max} &= \max_{t \in [0,T]} \max_{s \in \Lambda} \epsilon(s,t), 
  \quad \mathrm{where} \quad
  \epsilon(s,t) =  \mathrm{diam}(\Theta(s,t)) = \max_{\bm{x}, \bm{y}\in\Theta(s,t)} |\bm{x}- \bm{y}|, \label{eq:def_epsilon} \\ 
  \epsilon_s & = \|\partial_s R_1\|_{L^{\infty}(0,T;L^{\infty}(\Lambda))} + \|\partial_s R_2\|_{L^{\infty}(0,T;L^{\infty}(\Lambda))}.  \label{eq:def_epsilon_t}
\end{align}
We assume that as $\epsilon_{\max} \rightarrow 0$, 
\begin{align}
  \epsilon(s,t) \lesssim R_1(s, t) \lesssim \epsilon (s,t)  \,\,\, \mathrm{and} \,\,\,\epsilon(s,t) \lesssim R_2(s, t) \lesssim \epsilon (s,t).   \label{eq:assump_radii}
\end{align}  
This implies that: 
\begin{align}
  \epsilon(s,t) \lesssim P(s,t)\lesssim \epsilon (s,t) , \,\,\, \mathrm{and} \,\,\, \epsilon^2 (s,t) \lesssim A(s,t) \lesssim \epsilon^2(s,t),
  \label{eq:relation_area_diam}
\end{align}
with (implicit) inequality constants independent of $\Omega_v, \Omega_s$. 
In what follows, $K$ will denote a generic constant independent of $\epsilon_{\max}$ and of the norms of $c, \cc, c_v,$ and $c_s$. This generic constant $K$  may take different values when used in different places and may depend on the final time and on the material parameters. Hereinafter, we will use $A \lesssim B$ if there exists a generic constant $K$ as defined above such that $A \leq K B$. 
\begin{lemma}[Poincar{\'e} inequality over $\Theta$]
  \label{assump:poincare}
    For a.e.~$s \in \Lambda$ and $t \in (0,T]$, the following
      Poincar\'{e} inequality holds with inequality constant $K_p$
      independent of $\epsilon(s,t)$:
      \begin{equation}
        \|v - \avg{v}\|_{L^2(\Theta(s,t))} \leq K_p \epsilon(s,t)\|\nabla v\|_{L^2(\Theta(s,t))},
        \quad v \in H^1(\Theta(s,t)).
        \label{eq:Poincare_inequality}
      \end{equation}
\end{lemma}
\begin{proof}
  See for example \cite[Section 3.3]{guermond2021finite}. The dependence on the diameter is recovered from standard scaling arguments where the constant $K_p$ depends on a scaled annulus which is independent of $t$, $R_1,$ and $R_2$.
\end{proof}

\begin{example}[{Poincar{\'e} inequality over $\Theta$}]\label{ex:poincare}
  We can also numerically study the behaviour of the constant $K_p$
  in~\eqref{eq:Poincare_inequality} via the following eigenvalue
  problem: find $u\in H^1(\Theta)$ and $\lambda > 0$ such that
  \begin{equation}
    \label{eq:poincare_gevp}
    \begin{aligned}
      -\Delta u &= \lambda (u - \avg{u}) &\quad\text{in }\Theta,\\
      -\nabla u \cdot\bm{n} &= 0 &\quad \text{on }\partial\Theta.
    \end{aligned}
  \end{equation}
  Denoting by $\lambda_1$ the smallest eigenvalue of \eqref{eq:poincare_gevp},
  it follows that $\lambda^{-1/2}_1 = K_p \epsilon$.

  To investigate how $K_p$ varies with the sizes of annular domains,
  we let $\Theta$ be an annulus with inner and outer radii $R_1$ and
  $R_2$, respectively. We are interested in studying the cases where
  (i) $R_2$ is fixed while $R_1 \rightarrow 0$, and (ii) $R_2$ is
  fixed and $R_1 \rightarrow R_2$ (corresponding to $\epsilon
  \rightarrow 0$, and covered by the theoretical result). We solve the
  eigenvalue problem~\eqref{eq:poincare_gevp} numerically via
  continuous linear finite elements defined relative to uniform meshes
  of the annuli using the FEniCS finite element
  software~\cite{logg2012automated} and the SLEPc
  eigen solvers~\cite{slepc}, and a relative difference between
  smallest eigenvalue approximations on consecutive meshes of
  0.1\%. The smallest approximate eigenvalue is denoted $\tilde
  \lambda_1$, For both cases, we observe that $\tilde \lambda_{1}$
  scales linearly with the diameter $\epsilon = 2 R_2$ of $\Theta$
  (Figure~\ref{fig:poincare}, left). Denoting the estimated slope by
  $\tilde{K}_p \approx K_p$, we further observe that $K_p$ remains
  bounded, both as $R_1 \rightarrow 0$ and $R_1 \rightarrow R_2$
  (Figure~\ref{fig:poincare}, right).
\
  

  \begin{figure}
    \centering
    \includegraphics[height=0.35\textwidth]{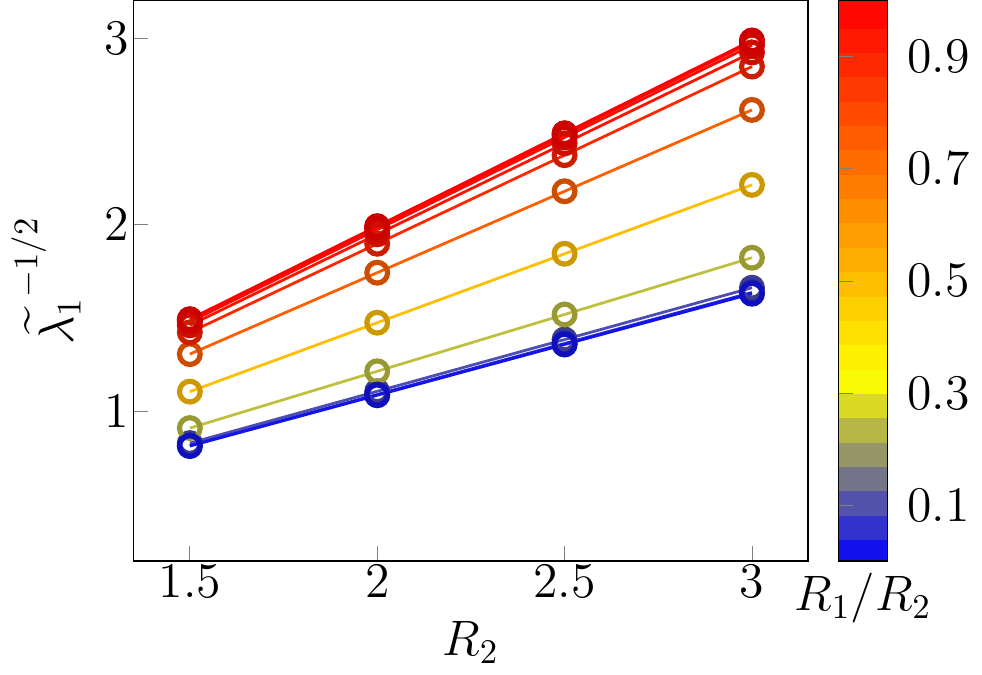}
    \hfill
    \includegraphics[height=0.35\textwidth]{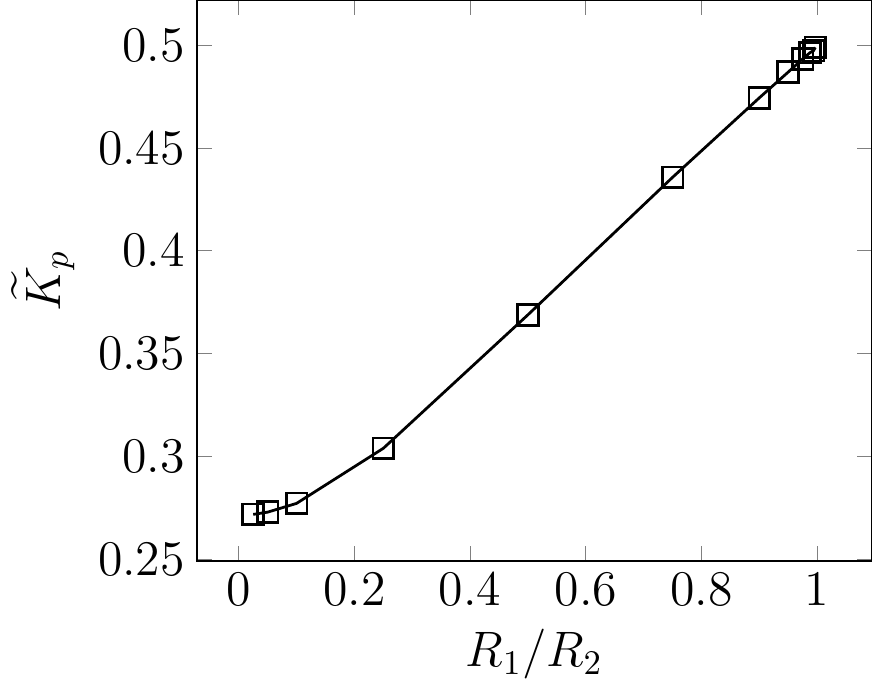}
    \vspace{-10pt}    
    \caption{
      Numerical investigation of the Poincar{\'e} inequality on annular domains (\Cref{ex:poincare}).
      Left: Linear scaling of the (approximate) smallest eigenvalue of \eqref{eq:poincare_gevp}
      related to the constant $K_p\epsilon$ in \eqref{eq:Poincare_inequality}. Right:
      Dependence of $K_p$ on the ratio of radii reveals that both limits
      $R_1\rightarrow 0$ and $R_1\rightarrow R_2$ lead to bounded $K_p$. }
    \label{fig:poincare}
  \end{figure}
\end{example}

For a convex and regular domain such as a circle or ellipse, it is well-known that a Sobolev trace inequality holds~\cite{brenner2008mathematical}. However, is this also the case for (nearly) annular domains? The subsequent~\Cref{lem:trace} addresses this question affirmatively.
\begin{lemma}[Trace inequality over $\Theta$]
  For a.e.~$s \in \Lambda$ and $t \in (0,T]$, the following trace inequality holds with $K_{\mathrm{tr}}$ independent of $\epsilon(s, t)$:  
    \begin{equation}
      \|v\|_{L^2(\partial \Theta(s,t))}^2 \leq K_{\mathrm{tr}} \left( \epsilon (s,t)^{-1} \|v\|_{L^2( \Theta(s,t))}^2 + \epsilon(s,t) \|\nabla v\|_{L^2(\Theta(s,t))}^2 \right)
      \quad v \in H^1(\Theta(s,t)),
      \label{eq:trace_estimate_general}
    \end{equation}
    where $\epsilon(s, t) = \diam(\Theta(s, t))$.
    \label{lem:trace}
\end{lemma} 

\begin{proof}
  For the circular case, if $\Theta(s) = \{(r\cos(\theta),r\sin(\theta)), 0 \leq r <
  R_2(s,\theta), 0 \leq \theta \leq 2 \pi \} $, this inequality is
  well known, Section 1.6 in \cite{brenner2008mathematical}. We use
  similar arguments to extend the proof to an annulus.

  Suppose now that $R_1 > 0$ and let $\Theta(s) = \{(r \cos(\theta), r \sin(\theta)),
  R_1(s,\theta) < r < R_2(s,\theta), 0 \leq \theta \leq 2 \pi \} $. We
  omit $s, t$ in the notation for the sake of brevity. Let $ (r ,
  \theta) \in (R_1,R_2) \times [0,2 \pi]$. We write
\begin{align*}
  R_1^2 u^2(R_1,\theta) - r^2 u^{2}(r,\theta) = - \int_{R_1}^r \partial_{z} (z^2 u^2(z, \theta))\dz. 
\end{align*}
Here, for simplicity, we write $u(r,\theta) = u(r\cos(\theta), r\sin(\theta))$. Thus, we have that 
\begin{align*}
 R_1^2 u^2(R_1,\theta)  \leq  r^2 u^{2}(r,\theta) + \int_{R_1}^r |\partial_{z} (z^2 u^2(z,\theta))| \dz. 
\end{align*}
Integrating over $(R_1,R_2)$ and over $[0,2\pi]$, we find that 
\begin{multline*}
  R_1 (R_2 - R_1) \|u\|^2_{L^2(\partial \Theta_1)} \leq R_2 \|u\|^2_{L^2(\Theta)} \\
  + \int_{0}^{2\pi} \int_{R_1}^{R_2}\int_{R_1}^r  2 |u(z, \theta)| \, |\partial_{z} u(z,\theta)| z^2 \dz \dr \dtheta 
  +  \int_{0}^{2\pi} \int_{R_1}^{R_2}\int_{R_1}^r  2z u^2(z,\theta) \dz \dr \dtheta ,
 \end{multline*}
where $\partial \Theta_1$ denotes the inner circle of $\Theta$.  Simplifying the last term and using Cauchy--Schwarz inequality for the penultimate, we obtain: 
 \begin{align*}
  R_1 (R_2 - R_1) \|u\|^2_{L^2(\partial \Theta_1)} \leq (R_2 + 2(R_2-R_1)) \|u\|^2_{L^2(\Theta)} + 2(R_2 -R_1) R_2 \|u\|_{L^2(\Theta)} \|\nabla u\|_{L^2(\Theta)}. 
 \end{align*} 
 With assumption~\eqref{eq:assump_radii} and Young's inequality, we obtain:
 \begin{align*}
  \|u\|^2_{L^2(\partial \Theta_1)} \lesssim   \epsilon^{-1} \|u\|^2_{L^2(\Theta)} + \epsilon \|\nabla u\|^2_{L^2(\Theta)}.
 \end{align*} 
 Similar arguments yield the same bound over the outer circle $\partial \Theta_2$ for $\|u\|_{L^2(\partial \Theta_2)}^2$. Adding the two bounds gives the result. The above computations are for smooth functions. The result for functions in $H^1(\Theta)$ follows by density.  \end{proof}

We now turn to consider a trace inequality for the surrounding domain
$\Omega_s$ (Lemma~\ref{lemma:trace_ineq_omegas}) by way of an
extension operator (Lemma~\ref{lemma:extension}) first introduced and
studied in~\cite{sauter1999extension}.

\begin{lemma}[Extension operator]
  \label{lemma:extension}
  For any $t \in [0,T]$ and $k \in \{1,2\}$, there exists an extension operator 
$\mathcal{E}(t): H^k(\Omega_s(t)) \rightarrow H^k(\Omega)$  satisfying $ \mathcal{E}(t)v \vert_{\Omega_s(t)} = v \vert_{\Omega_s(t)}, \,\, \mathcal{E}(t)v \vert_{\Gamma(t)} = v \vert_{\Gamma(t)}$ and such that
  \begin{equation} \| \mathcal{E}(t) v \|_{H^k(\Omega)} \leq K_{\mathcal{E}} \|v\|_{H^k(\Omega_s(t))},   \quad  \foralls v \in H^k(\Omega_s(t)),
    \label{eq:extension_operator} 
  \end{equation}
  with a constant $K_{\mathcal{E}}$ independent of $\epsilon_{\max}$ and $t$. 
\end{lemma}
\begin{proof}
The construction of the extension operator and the proof of the continuity bound are very similar to \cite[Theorem 2.1]{sauter1999extension}. For completeness, we provide some details adapted to our geometrical setting. First, we define the extension from a  fixed domain $\tilde{B}= B_1 \backslash B$ to $B_1$ where $B$ and $B_1$ are cylindrical domains of radii $1$ and $2$ respectively. 
 Let $\mathcal{E}_0: H^1(\tilde{B})\rightarrow H^1(\R^3)$ be the extension operator as defined in \cite[Section 5.4]{evans2009partial}. We have the following two bounds: 
 \begin{align*}
\|\mathcal{E}_0 u\|_{H^1(\R^d)} \leq  K_1 \|u\|_{H^1(\tilde{B})}, \,\,\,\, \|\mathcal{E}_0 u\|_{H^2(\R^d)} \leq  K_2 \|u\|_{H^2(\tilde{B})}.
 \end{align*}
 Let $H^k_{0}(B) = \{v \in H^k(B), \partial^{\alpha} v = 0, \,\, |\alpha | < k \,\, \mathrm{on} \,\, \partial B\}$ and define $z \in H^k_0(B)$ such that
 \begin{equation}
 \sum_{|\alpha| = k} (\partial^{
  \alpha} z , \partial^{\alpha} q)_{B} = \sum_{|\alpha| = k} (\partial^{\alpha} (\mathcal{E}_0 u), \partial^{\alpha} q)_{B}, \,\,\, \forall q \in H^k_0(B). 
 \end{equation} 
One can show that $z$ is well-defined by the Lax-Milgram theorem since a Poincar{\'e} inequality holds in $H^k_0(B)$, and we have that 
\begin{equation}
  \|z\|_{H^k(B)} \leq \tilde{K}_k \|\mathcal{E}_0 u\|_{H^k(B)}.  
\end{equation}
The extension operator $\mathcal{E}_{\tilde{B}}:H^1(\tilde{B}) \to H^1(B_1)$ is then defined as follows:
\begin{align}
  \mathcal{E}_{\tilde{B}} u\, (\bm{x}) = \begin{cases}
    u(\bm{x}), & \bm{x} \in \tilde{B} \\ 
    \mathcal{E}_0 u(\bm{x})- z(\bm{x}), & \bm{x} \in B 
  \end{cases}. 
\end{align}
To show continuity of $\mathcal{E}_{\tilde{B}}$, we have that for $k \in \{1,2\}$: 
\begin{multline}
  \|\mathcal{E}_{\tilde{B}}u \|^2_{H^{k}(B_1)}
  = \|u\|^2_{H^k(\tilde{B})} + \|  \mathcal{E}_0 u- z\|_{H^k(B)}^2 
  \leq \|u\|^2_{H^k(\tilde{B})} +( \|  \mathcal{E}_0 u\|_{H^k(B)}+ \|z\|_{H^k(B)})^2  \\ 
  \leq \|u\|^2_{H^k(\tilde{B})}+ (1 +\tilde{K}_k)^2 K^2_k \|u\|^2_{H^k(\tilde{B})} \leq K_k^f\|u\|^2_{H^k(\tilde{B})} .
\end{multline}
In the above, we let $K^f_k = 1 + (1 +\tilde{K}_k)^2 K^2_k$ which clearly depends on $B, \tilde{B}$ and $B_1$. A key property of this extension is that $\mathcal{E}_{\tilde{B}} p = p $ for all polynomials $p$ of degree less than $k$, see \cite[Lemma 2.1]{sauter1999extension}. By choosing $p$ as the average of $u$ for $k =1$ or the Lagrange interpolant of degree $1$ for $k=2$,   this observation yields the following bounds on the semi-norms:
\begin{equation}
|\mathcal{E}_{\tilde{B}}u |^2_{H^k(B_1)} = |\mathcal{E}_{\tilde{B}}(u - p) |^2_{H^k(B_1)} \leq K_k^f \|u -p\|^2_{H^{k}(\tilde{B})} \leq K_k^f K_2 |u|^2_{H^{k}(\tilde{B})}, 
\end{equation}
for some constant $K_2$. Now, we define the extension operator from $ H^k(\Omega_s(t)) \to H^k(\Omega)$ as:
\begin{equation}
\mathcal{E}(t) v = \begin{cases}  \mathcal{E}_{\tilde{B}}(v \circ \chi_{R_2(t)}^{-1}) \circ \chi_{R_2(t)},   &\mathrm{in} \,\,  B_{2 \, R_2(t)} \\ 
  v, & \mathrm{in} \,\,  \Omega_s(t) \backslash B_{2 \, R_2(t)}
\end{cases},
\end{equation}
where $B_{2 R_2(t)}$ is the cylinder surrounding $B_{R_2}$ of radius $2R_2$ and $$\chi_{R_2(t)} ((s, R_2 \cos(\theta), R_2 \sin(\theta)) )= (s, \cos(\theta), \sin(\theta)), \,\,\, \forall (s,\theta)  \in \Lambda \times [0,2\pi].$$ The continuity  of $\mathcal{E}$ then follows from  a scaling argument and \eqref{eq:assump_radii} which yield that 
\begin{align*}
|v \circ \chi_{R_2(t)}^{-1}|^2_{H^i(\tilde{B})} \lesssim \epsilon_{\max}^{-3+2i} |v|^2_{H^i(B_{2R_2(t) }\backslash B_{R_2(t)})}, \,\,\, |\hat{v}\circ \chi_{R_2(t)} |^2_{H^{i}(B_{2R_2(t)})} \lesssim \epsilon_{\max}^{3-2i}| \hat{v}|^2_{H^{i}(B_1)}, 
\end{align*}
for $i \in \{0,1,2\}.$
Thus, we obtain the following:
\begin{equation}
\begin{split}
  \|\mathcal{E}(t)v\|_{H^k(\Omega)}^2
  &= \|v\|_{H^k(\Omega_s(t) \backslash B_{2\,R_2(t)})}^2 + \|\mathcal{E}_{\tilde{B}}(v \circ \chi_{R_2(t)}^{-1}) \circ \chi_{R_2(t)}\|_{H^k(B_{2 \, R_2(t)})}^2 \\ 
  &\lesssim  \|v\|_{H^k(\Omega_s(t))}^2 + \sum_{i=0}^k \epsilon_{\max}^{3-2i} |\mathcal{E}_{\tilde{B}}(v \circ \chi_{R_2(t)}^{-1})|^2_{H^i(B_1)} \\ 
  &\leq  \|v\|_{H^k(\Omega_s(t))}^2  +  \sum_{i=0}^k \epsilon_{\max}^{3-2i} K_k^f K_2 |v \circ \chi_{R_2(t)}^{-1}|^2_{H^i(\tilde{B})}\\ 
  &\lesssim  \|v\|^2_{H^{k}(\Omega_s(t))}. \qedhere
 \end{split} 
\end{equation}
\end{proof}

\begin{lemma}[Trace inequality over $\Omega_s$]
  \label{lemma:trace_ineq_omegas} 
  There exists a constant $K_{\Gamma}$ independent of $t$ and of $\epsilon_{\max}$ such that 
  \begin{align}
    \|v\|_{L^2(\Gamma(t))} \leq K_{\Gamma} (\epsilon_{\max}\, |\ln \epsilon_{\max}|)^{1/2} \,  \|v\|_{H^1(\Omega_s(t))},  \quad \foralls v \in H^{1}(\Omega_s(t)).
    \label{eq:trace_eq_Omegas}
  \end{align}
\end{lemma}

\begin{proof}
  Without loss of generality, we consider the case of $\Omega_v$ being an annular cylinder domain and $\Omega_s$ its outer surroundings. We have for $v \in H^{1}(\Omega_s(t))$: 
  \begin{align}
    \|v\|^2_{L^2(\Gamma(t))} = \|\mathcal{E} v\|_{L^2(\Gamma(t))}^2
    = \int_{\Lambda} \|\mathcal{E} v\|_{L^2(\partial \Theta_2(s,t))}^2 \ds.
    \label{eq:writing_norm}
  \end{align} 
  We use ideas from the  proofs of \cite[Lemma 2.1 and Lemma 2.2]{koppl2018mathematical} where we adapt the arguments to 3D. We  write for a.e.~$s \in \Lambda, t \geq 0$,
  \begin{equation}
    \|\mathcal{E}v \|_{L^2(\partial \Theta_2(s,t))} \leq \|\mathcal{E}v
    - \overline{\mathcal{E}v }\|_{L^2(\partial \Theta_2(s,t))}
    + \|\overline{\mathcal{E}v }\|_{L^2(\partial \Theta_2(s,t))}.
    \label{eq:triang_trace_pf}
  \end{equation}
The first term is bounded by a Stekloff type inequality \cite{kuttler1969inequality}:
\begin{equation}
  \|\mathcal{E}v - \overline{\mathcal{E}v }\|_{L^2(\partial \Theta_2(s,t))}  \leq K_{\mathrm{st}} \epsilon(s,t)^{1/2} \|\nabla (\mathcal{E} v)\|_{L^2(\Theta_2(s,t))}. 
  \label{eq:trace_pf_0}
\end{equation}
For the second term in \eqref{eq:triang_trace_pf}, observe that by definition of the perimeter average
\begin{equation}
  \|\overline{\mathcal{E}v }\|_{L^2(\partial \Theta_2(s,t))} = |\partial \Theta_2(s,t)|^{1/2} |\overline{\mathcal{E}v }| .
\end{equation}
From the proof of \cite[Lemma 2.1]{koppl2018mathematical},  we further have for $p > 2$
\begin{equation}
  | \overline{\mathcal{E}v}| \leq  \left( \pi R_2(s,t)^{2} \right )^{-1/p} \|\mathcal{E} v\|_{L^{p}(\Theta_2(s,t))}
  + \frac{1}{2\sqrt{\pi}} \|\nabla \mathcal E {v} \|_{L^2(\Theta_2(s,t))}. 
\end{equation}
Hence, we obtain:
\begin{equation}
  \| \overline{\mathcal{E} v }\|_{L^2(\partial \Theta_2(s,t))}
  \leq K \left( \epsilon(s,t)^{1/2-2/p} \|\mathcal{E} v \|_{L^p(\Theta_2(s,t))} +  \epsilon(s,t)^{1/2}\|\nabla \mathcal{E} v\|_{L^2(\Theta_2(s,t))} \right ).
\end{equation} 
Upon substituting in \eqref{eq:writing_norm}, we have that 
\begin{align}
  \int_{\Lambda} \|\mathcal{E} v\|_{L^2(\partial \Theta_2(s,t))}^2 \leq K \int_{\Lambda} (\epsilon_{\max}^{1-4/p} \|\mathcal{E} v\|_{L^p(\Theta_2(s,t))}^2 + \epsilon_{\max} \|\nabla \mathcal{E} v\|^2_{L^2(\Theta_2(s,t))} ).  
  \label{eq:bound_1_trace}
\end{align}
Consider now a fixed cylindrical domain $B_{\Lambda}$ around the centerline $\Lambda$ with cross-sections $\Theta_{\Lambda}(s,t)$. We emphasize that $B_{\Lambda}$ does not depend on $\epsilon_{\max}$. Observe that for $\epsilon_{\max}$ small, $\Omega_v(t) \subset B_{\Lambda}$.  Let $\tilde{B}_{\Lambda}$ be another cylinder such that $B_{\Lambda} \subset \tilde{B}_{\Lambda} \subset \Omega$ with cross-sections $\tilde{\Theta}_{\Lambda}(s,t) \supset \Theta_{\Lambda}(s,t)$.  Define $\chi$ to be a smooth cut-off function on $B_{\Lambda}$ such that $\chi = 1$ in $B_{\Lambda}$ with compact support in $\tilde{B}_{\Lambda}$. 
By construction, we have 
\begin{align*}
\|\mathcal{E}v\|_{L^{p} (\Theta_2(s,t))}   \leq \|\mathcal{E}v\|_{L^p (\Theta_{\Lambda} (s,t))}= \|\chi (\mathcal{E} v) \|_{L^p (\Theta_{\Lambda} (s,t))} \leq \|\chi (\mathcal{E} v) \|_{L^p (\tilde{\Theta}_{\Lambda} (s,t))}.
\end{align*}
Since $\chi( \mathcal{E}v )\in H^1_0(\tilde{\Theta}_{\Lambda}(s,t))$, we apply the Sobolev embedding result  in 2D which gives a constant with an  explicit dependence on $p$ \cite[eq (6.20)]{thomee2007galerkin}: 
\begin{equation*}
  \|\chi (\mathcal{E} v) \|_{L^p (\tilde{\Theta}_{\Lambda} (s,t))} \leq  K p^{1/2} \|\nabla ( \chi (\mathcal{E} v))\|_{L^2 (\tilde{\Theta}_{\Lambda} (s,t))}  \leq K p^{1/2} \|\mathcal{E} v \|_{H^1(\tilde{\Theta}_{\Lambda} (s,t))}.
\end{equation*}
The above constant depends on $\tilde{\Theta}_{\Lambda} (s,t)$ and on $\chi$ but not $\epsilon_{\max}$. Substituting in \eqref{eq:bound_1_trace}, and choosing $p = |\ln \epsilon_{\max}|$ yields: 
\begin{align*}
  \|\mathcal{E} v\|_{L^2(\Gamma(t))}^2
  &\leq K  \int_{\Lambda} \left ( \epsilon_{\max} |\ln \epsilon_{\max} |\|\mathcal{E} v\|_{H^1(\tilde{\Theta}_{\Lambda}(s,t))}^2
  + \epsilon_{\max} \|\nabla \mathcal{E} v\|^2_{L^2(\Theta_2(s,t))} \right ) \\ 
  &\leq K \epsilon_{\max} \left (
  |\ln \epsilon_{\max} | \|\mathcal{E} v\|_{H^1(\tilde{B}_{\Lambda})}^2 + \|\mathcal{E}v\|_{H^1(\Omega)}^2 \right )
  \leq K \epsilon_{\max} |\ln \epsilon_{\max} |\|\mathcal{E}v\|_{H^1(\Omega)}^2.
\end{align*}
Using \eqref{eq:extension_operator} in the above concludes the proof.
\end{proof}


\begin{example}[{Trace inequality over $\Omega_s$}]
  \label{ex:trace}
  %
  The scaling law of \Cref{lemma:trace_ineq_omegas} can be demonstrated
  numerically by considering the following Stekloff eigenvalue problem \cite{Stekloff1902}: find
  $u\in H^1(\Omega_s)$ and $\lambda > 0$ such that
  \begin{equation}\label{eq:stekloff_gevp}
    \begin{alignedat}{2}
    \Delta u &= u && \quad \text{ in }\Omega_s,\\
    \nabla u\cdot\bm{n} &= \lambda u &&\quad\text{ on }\Gamma,\\
    \nabla u\cdot\bm{n} &= 0 &&\quad\text{ on }\partial\Omega_s\setminus\Gamma.\\
    \end{alignedat}
  \end{equation}
  More precisely, for $\lambda_1$ being the smallest non-zero eigenvalue of \eqref{eq:stekloff_gevp},
  there holds that
  \begin{align}\label{eq:stekloff_eigw}
    \|u\|_{L^2(\Gamma)} \leq \lambda_1^{-1/2}  \|u\|_{H^1(\Omega_s)},  \quad \foralls v \in H^{1}(\Omega_s) .
  \end{align}
  Thus, approximations to $\lambda_1$ in \eqref{eq:stekloff_gevp} can
  be used to estimate the bound in \eqref{eq:trace_eq_Omegas}. Now,
  consider an embedding domain $\Omega_s$ (also) in the shape of an
  cylinder with unit height and unit outer radius. Consider an inner
  cylinder $\Omega_v$ with diameter $\epsilon = 2 R_1$ (and unit
  height) and consider a decreasing sequence of $R_1$s. As
  in~\Cref{ex:poincare}, we approximate this smallest non-zero
  eigenvalue using the discretization of~\eqref{eq:stekloff_gevp} by
  continuous linear elements defined relative to a series of uniformly
  refined meshes, and deem the eigenvalues converged when the relative
  difference between refinements is less then 5\%. Clearly, the trace
  constant decreases with decreasing $\epsilon_{\max}$ (\Cref{fig:trace}).
  We note that the data are well-fitted by the theoretically established
  $\epsilon_{\max}^{1/2}\lvert \ln \epsilon_{\max} \rvert^{1/2}$ expression especially for small radii.
  \begin{figure}[H]
    \centering
    \includegraphics[width=0.45\textwidth]{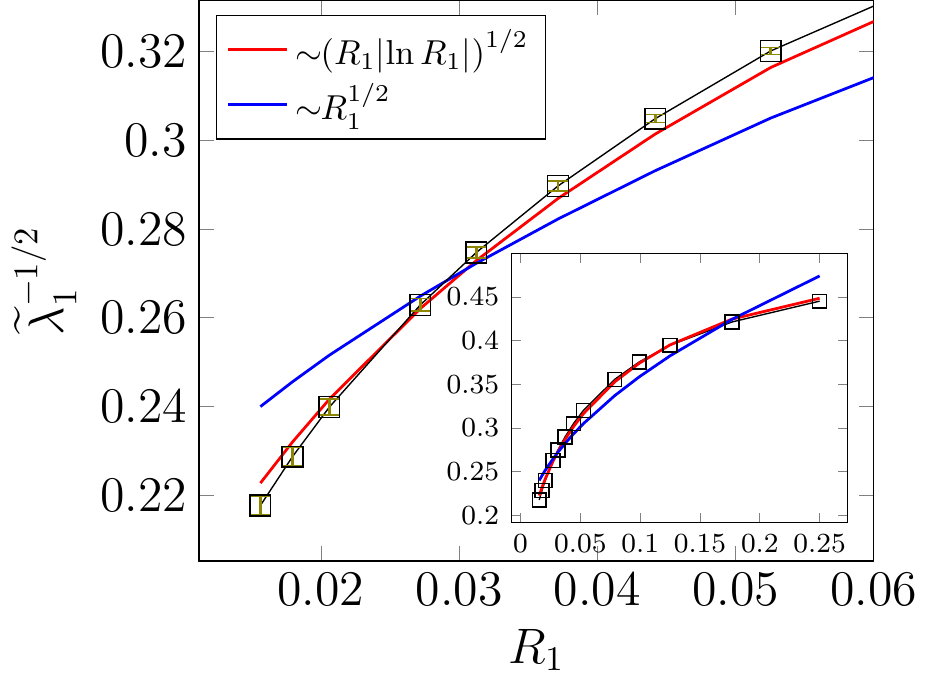}
    \vspace{-10pt}
    \caption{ Numerical investigation of the trace inequality
      (\Cref{ex:trace}) for vessels of decreasing diameters
      $\epsilon = 2 R_1$.}
    \label{fig:trace}
  \end{figure}
\end{example}

\section{Analysis of the modelling error}
\label{sec:modelling_error}

We next turn to address the following question: how large of a
modelling error has been introduced by the derivation and associated
assumptions of the coupled 3D-1D model
in~\Cref{subsec:derivation_1dmodel}? We begin by considering the
modelling error associated with the cross-section average
concentration in the vessel, before turning to the modelling error in
the surroundings. 
We will make use of a duality argument, and therefore introduce and
analyze the stability of an associated dual problem before turning to
the modelling error estimates. 

\subsection{Well-posedness and stability of a dual transport problem}

We consider the properties of a backward-in-time dual transport
problem, defined in association with the forward vessel transport
problem of \Cref{prop:exist_3d_3d},
in~\Cref{lemma:backward_parabolic}. A key aspect is to appropriately
account for the moving domain and time-derivatives with respect to
moving frames. We therefore explicitly track the domain dependence on
time $t$.
\begin{lemma}[Backward-in-time dual problem]
  \label{lemma:backward_parabolic}
  The following problem is well posed: given $g \in L^2(0,T;L^2(\Omega_v(t)))$, find $h \in W = \{ L^2(0,T; H^1(\Omega_v(t))) \mid  \dot{h} \in L^2(0,T; H^{-1}(\Omega_v(t)))\}$ and $h(T) = 0$ in $ \Omega_v(T)$ such that for $t \in (0,T)$ and for all $\varphi \in H^1(\Omega_v(t))$:
\begin{multline}
  - \langle \dot{h}(t), \varphi \rangle_{H^{-1}(\Omega_v(t))}
  + (D_v \nabla h(t), \nabla \varphi)_{\Omega_v(t)} + (\xi h(t), \varphi)_{\Gamma(t)}\\  - ((\bm{u}_v - \bm{w})\cdot \nabla h(t),  \varphi)_{\Omega_v(t)} 
  = (g(t),\varphi)_{\Omega_v(t)}.
  \label{eq:back_time_dual}
\end{multline}
In addition, the following stability bound holds 
\begin{multline}
  \|h\|_{L^{\infty}(0,T;L^2(\Omega_v(t)))} + 
  \|D_v^{1/2} \nabla h\|_{L^2(0,T;L^2(\Omega_v(t)))}\\  + \|\xi^{1/2} h\|_{L^2(0,T;L^2(\Gamma(t)))} 
  \leq K_{b} \|g\|_{L^2(0,T;L^2(\Omega_v(t)))} ,
  \label{eq:regularity_backward_pb}
  \end{multline}   
  where $K_b$ is independent of $\epsilon_{\max}$ but depends on the final time $T$, on $\|  \nabla \cdot \bm{w}\|_{L^{\infty} (0,T;L^{\infty}(\Omega))}$,  and on $ \|D_v^{-1/2}\tilde{\bm{u}}_v\|_{L^\infty(0,T;L^{\infty}(\Omega_v(t)))}$. 
\end{lemma}
\begin{proof}
We consider the forward-in-time solution $z \in W$ and $z(0) = 0$ in $ \tilde{\Omega}_v(0)$ solving for all $\varphi \in H^1(\tilde{\Omega}_v(t))$ 
\begin{equation*}
  \langle \dot{z}, \varphi \rangle_{H^{-1}(\tilde{\Omega}_v(t))}
  + (\overset{\leftarrow}{D_v} \nabla z, \nabla \varphi)_{\tilde{\Omega}_v(t)} + (\overset{\leftarrow}{\xi} z, \varphi)_{\tilde{\Gamma}(t)} - ((\overset{\leftarrow}{{\bm{u}}_v} - \overset{\leftarrow}{\bm{w}}) \cdot  \nabla z,  \varphi )_{\tilde{\Omega}_v(t)} = (\overset{\leftarrow}{g},\varphi)_{\tilde{\Omega}_v(t)},  
  \end{equation*}
where a $``\leftarrow"$ over a function indicates that we reverse the time, e.g., $\overset{\leftarrow}{g}(t) = g(T-t)$. The domain $\tilde{\Omega}_v(t) = \Omega_v(T-t)$ (similarly $\tilde{\Gamma}(t) = \Gamma(T-t)$) is given by the flow map $\overset{\leftarrow}{\bm{\psi}}(\bm{x}, t) = \bm{\psi}(\bm{x},T-t)$. 
Setting $ h(t) = z(T-t)$ for $t \in [0,T]$, we recover the solution to \eqref{eq:back_time_dual} since $\dot{h} = -\dot{z}$ and $h(T) = z(0) = 0$. Verifying the existence and uniqueness of $z$ then follows from the abstract framework in \cite{alphonse2015abstract} and from very similar arguments to the proof of Proposition~\ref{prop:exist_3d_3d}.

Moreover, choose $\varphi = h \in L^2(0,T;H^1(\Omega_v(t)))$ in \eqref{eq:back_time_dual}, integrate over time $\tau \in [t,T]$  and use the following formula \cite[Theorem 2.40 and Corollary 2.41]{alphonse2015abstract}: 
\begin{align}
  - 2 \int_t^T \langle \dot h, h  \rangle_{H^{-1}(\Omega_v(\tau))} \mathrm{d} \tau
  = \|h(t)\|^2_{L^2(\Omega_v(t))} + \int_{t}^T (h, h \nabla \cdot \bm{w})_{\Omega_v(\tau)} \mathrm{d} \tau. 
\end{align} 
Along with H\"{o}lder's inequality, this yields: 
\begin{multline}
\frac12  \|h(t)\|^2_{L^2(\Omega_v(t))}  + \|D_v^{1/2} \nabla h\|^{2}_{L^2(t,T;L^2(\Omega_v))} + \|\xi^{1/2} h\|^{2}_{L^2(t,T;L^2(\Gamma))} \\
  \leq \|g\|_{L^2(t,T;L^2(\Omega_v(t)))} \|h\|_{L^2(t,T;L^2(\Omega_v(t)))} + \frac12 \|\nabla \cdot \bm{w}\|_{L^\infty(0,T;L^\infty(\Omega_v(t)))} \|h\|^2_{L^2(t,T;L^2(\Omega_v(t)))} \\
  + \|D_v^{-1/2}\tilde{\bm{u}}_v\|_{L^\infty(0,T;L^\infty(\Omega_v(t)))} \|h\|_{L^2(t,T;L^2(\Omega_v(t)))} \|D_v^{1/2}\nabla h\|_{L^2(t,T;L^2(\Omega_v(t)))}.
\end{multline}
Applying Young's inequality for the first and last term on the right hand side above results in:
\begin{multline*}
  \frac12 \|h(t)\|^2_{L^2(\Omega_v(t))}
  + \frac12 \|D_v^{1/2} \nabla h\|^{2}_{L^2(t,T;L^2(\Omega_v))}
  + \|\xi^{1/2} h\|^{2}_{L^2(t,T;L^2(\Gamma))} 
  \leq \frac12 \|g\|^2_{L^2(0,T;L^2(\Omega_v(t)))} \\
  + \frac12 \left (1 +\frac12
  \|\nabla \cdot \bm{w}\|_{L^\infty(0,T;L^\infty(\Omega_v(t)))} 
  + \|D_v^{-1/2}\tilde{\bm{u}}_v\|^2_{L^\infty(0,T;L^\infty(\Omega_v(t)))} \right ) \|h\|^2_{L^2(t,T;L^2(\Omega_v(t)))}.  
\end{multline*}
The result can then be concluded by Gr\"onwall's inequality, see e.g~\cite[Appendix B.k]{evans2009partial}. 
\end{proof}

\subsection{Model error introduced in the derivation of the 1D model}

With this dual stability result at hand, we now turn to our first main
modelling error estimate, namely comparing (the 
extension of) the cross-section average vessel solution $\cc : \Lambda \times (0, T)
\rightarrow \R$ with its reference solution $c_v : \Omega_v \times (0,
t) \rightarrow \R$, the weak solution of~\eqref{eq:3d_model_1}. More
specifically, we aim to quantify the modelling error
\begin{equation*}
  \|c_v - E \cc \|_{L^2(0, T; L^2(\Omega_v))} .
\end{equation*}  
The (constant cross-section) extension $E$ from $H^1(\Lambda)$ to
$H^1(\Omega_v)$ is given by
\begin{equation}
  E \cc(s,r,\theta,t) = \cc(s,t), \quad
  \foralls (r, \theta) \in \Theta(s,t), \quad \foralls (s,t) \in \Lambda \times (0,T).
  \label{eq:natural}
\end{equation}
We will frequently write $\cc = E \cc$ when context allows this simplification.
\Cref{prop:first_model_err} gives the main modelling error estimate 
for the solutions in the vessel. 

\begin{prop}[Model error in the vessel]
  \label{prop:first_model_err}
  Let $c_v, c_s$ be weak solutions to the coupled 3D-3D transport problem~\eqref{eq:weak_form_3d3d} and assume that $c_v(0) \in H^{1}(\Omega_v)$.  
   Let $c, \cc$ be the weak solutions to the reduced coupled 3D-1D problem~\eqref{eq:coupled_3d_1d_weak} with $w_c = 1$. 
 Then, 
 \begin{multline}
  \label{eq:model_err_1d-3d} 
  \|c_v - \cc  \|_{L^2(0,T;L^2(\Omega_v))}   \\ 
     \lesssim   K_b  \, K_p \left ( \|f_v\|_{L^2(0,T;H^1(\Omega_v))} + \|\nabla c_v(0)\|_{L^2(\Omega_v)} \right ) \epsilon_{\max}  + (K_b \, K_{\Gamma} \,C_2) \, \, \epsilon_{\max}^{1/2} |\ln \epsilon_{\max}|^{1/2} \\ 
     + K_b(K_p + 1)(K_{\mathrm{tr}} + 1) C_1 \left(
     \epsilon_{\max}^{1/2} + \|u_{v,r}\|_{L^\infty(0,T;L^{\infty}(\Omega_v(t)))} + \|u_{v,\theta} \|_{L^\infty(0,T;L^{\infty}(\Omega_v(t)))} + \epsilon_s
     \right ).
 \end{multline}
Here, $C_1$ and $C_2$ depend on the material parameters and the solutions $c, \cc,$ and $c_s$ as 
  \begin{multline*}
    C_1 = \left ( \max_{s\in \Lambda, t \in [0,T]} \|u_{v,s}\|_{H^1(\Theta)} + \|\bm{w}\|_{L^{\infty}(\Omega)} \right ) \|\cc\|_{L^2(0,T;L^2(\partial \Omega_v (t)))} \\
    + \|\cc\|_{L^2(0,T;H^1(\Omega_v(t)))}  + \|D_v \partial_s \cc + \cc \avgus\|_{L^2(0,T;L^2( \Omega_v(t)))} + \|\xi (\cc - \overline{c})\|_{L^2(0,T;L^2(\Gamma(t)))} ,
  \end{multline*}
  and
  \begin{equation*}
   C_2 = \|\xi^{1/2}\|_{L^{\infty}(0,T;L^{\infty}(\Gamma(t)))} \left ( \|c_s\|_{L^2(0,T;H^1(\Omega_s(t)))} + \|c\|_{L^2(0,T;H^1(\Omega))} \right).   
  \end{equation*}
   In addition, there exists a constant $K$ depending only on the material parameters and the final time $T$ but not on $\epsilon_{\max}$ such that $C_1 + C_2 \lesssim K$. Under the additional  assumption that 
   \begin{equation} \|\overline{c} - \overline{c_s}\|_{L^2(0,T;L^2(\Gamma(t)))} \lesssim \epsilon_{\max}^{1/2}, \label{eq:epsilon:better}\end{equation} bound \eqref{eq:model_err_1d-3d}   can be improved by  replacing its last term by $K_b (K_{\mathrm{st}} K_{\mathcal{E}} + 1)C_2 \,\,\epsilon_{\max}^{1/2}$. 
   \end{prop}

Before presenting the proof, we remark that this proposition and in
particular~\eqref{eq:model_err_1d-3d} provides a rigorous bound on the
error in the vessel introduced by the derivation of the 3D-1D
model. For the error to converge to $0$ as $\epsilon_{\max}
\rightarrow 0$, one needs to assume that $u_{v,r}$ and $u_{v,\theta}$
are negligible at least for small $\epsilon_{\max}$.  In addition, if
$\epsilon_{s} \lesssim \epsilon_{\max}$ as $\epsilon_{\max} \rightarrow 0$, then one
recovers a convergence rate of $1/2$ (up to a log factor) with respect
to $\epsilon_{\max}$.  The additional
assumption~\eqref{eq:epsilon:better} essentially leads to an estimate
for the error induced by Assumptions \ref{assumption:shape_profiles}
and \ref{assumption:boundary_condition_A} of \Cref{subsec:derv_1d}
alone, without those of~\Cref{subsec:weak_form_3d_derivation}; i.e.~an
estimate for the error between $E\cc$ as given in
\eqref{eq:weighted_weak_form_1D_pb} and $c_v$ the solution of
\eqref{eq:weak_form_3d3d}. In this case, we can remove the log factor.


    \begin{proof}(Proposition \ref{prop:first_model_err})
     We proceed in three main steps to (I) derive a first identity for
     the modelling error by using a duality argument, (II) manipulate this identity by deriving the weak form satisfied by the extended solution $E \cc$,  and (III)
     bound its terms via Poincar\'e, trace,  Stekloff inequalities, and the regularity bound derived in Lemma \ref{lemma:backward_parabolic}. 

     \emph{Step I.}
     We first recall from~\eqref{eq:weak_form_3d3d} that the reference solution $c_v$ satisfies  
   \begin{equation}
     \langle \dot{c}_v, \phi \rangle_{H^{-1}(\Omega_v(t))} +  (\nabla \cdot \bm{w} c_v, \phi)_{\Omega_v(t)} + a_\mathrm{ref}(c_v,\phi)  = \ell_\mathrm{ref}(\phi), \quad \foralls \phi \in H^1(\Omega_v(t)), \label{eq:reference_weak_form}
   \end{equation}
   where we have introduced the two forms 
   \begin{alignat*}{2}
     a_{\mathrm{ref}} (c,\phi) & = (D_v \nabla c, \nabla \phi)_{\Omega_v(t)} + (\xi c, \phi)_{\Gamma(t)} - ((\bm{u}_v - \bm{w})c, \nabla \phi )_{\Omega_v(t)}
     && \quad \foralls c,\phi \in H^1(\Omega_v(t)),  \\ 
     \ell_\mathrm{ref}(\phi) & =  (\xi c_s, \phi)_{\Gamma(t)} + (f_v,\phi)_{\Omega_v(t)}
     && \quad \foralls \phi \in H^1(\Omega_v(t)). 
   \end{alignat*}
   To estimate the error $e \equiv c_v - E\hat{c}$, we proceed by duality. Namely, let $h$ be the solution of \eqref{eq:back_time_dual} with $g = e \in L^2(0,T;L^2(\Omega_v(t)))$. From \cite[Corollary 2.41]{alphonse2015abstract} and the fact that $h(T)=0$,
   the following integration by parts formula holds: 
   \begin{equation*}
     \int_0^T - \langle \dot{h}, e\rangle_{H^{-1}(\Omega_v(t))}
     = \int_0^T \langle \dot{e}, h \rangle_{H^{-1}(\Omega_v(t))} + \int_0^T ( e, h \nabla \cdot \bm{w})_{\Omega_v(t)} + ( e(0), h(0))_{\Omega_v(0)}.
   \end{equation*}
   With this identity, \eqref{eq:back_time_dual} tested with $e \in L^2(0,T;H^1(\Omega_v(t)))$ and integrated over $(0, T)$ reads:
   \begin{equation*}
     \int_0^T \langle \dot{e}, h \rangle_{H^{-1}(\Omega_v(t))}
     + \int_0^T (e, h \nabla \cdot \bm{w})_{\Omega_v(t)} + ( e(0), h(0))_{\Omega_v(0)} + \int_0^T a_{\mathrm{ref}} (e,h)= \int_{0}^T \|e\|_{L^2(\Omega_v(t))}^2. 
   \end{equation*}
   Subtracting the time-integrated~\eqref{eq:reference_weak_form}, combined with the observation that indeed $\dot{\hat{c}} \in L^2(0,T;L^2(\Omega_v(t)))$ (where we write $\hat{c}$ in place of $E \hat{c}$ here and in the following), we obtain the following identity for the modelling error $e$:
   \begin{multline}
     \int_{0}^T \|e\|_{L^2(\Omega_v(t))}^2 
     = - \int_0^T (\dot{\hat{c}}, h)_{\Omega_v(t)} -  \int_0^T (\hat{c}, h \nabla \cdot \bm{w})_{\Omega_v(t)} 
     -  \int_0^T a_{\mathrm{ref}} (\hat{c}, h) \\
     +  \int_0^T \ell_{\mathrm{ref}} (h)  + ( e(0), h(0))_{\Omega_v(0)}.
     \label{eq:dual_interm_0}
   \end{multline}

   \emph{Step II.} Next, we aim to derive an alternative expression for this error identity. By definition of the strong material derivative cf.~\eqref{eq:strong_mater_derivative}:
   \begin{equation}
     \int_0^T (\dot{\hat{c}}, h)_{\Omega_v(t)}
     +  \int_0^T (\hat{c}, h \nabla \cdot \bm{w})_{\Omega_v(t)}
     = \int_{0}^T (\partial_t \hat{c}, h)_{\Omega_v(t)}
     + \int_0^T (\nabla \cdot(\hat{c} \bm{w}), h)_{\Omega_v(t)}
     .
   \end{equation} 
   Note that this definition holds for $\hat{c} \in  \{ v \in L^2(0,T; H^1(\Omega_v(t)))\, \vert \, \dot{v}  \in L^2(0,T; L^2(\Omega_v(t))) \}$ by density of $\mathcal{D}(0,T; H^1(\Omega_v(t)))$ in such spaces, \cite[Lemma 2.38]{alphonse2015abstract}. Further, integrating by parts gives
   \begin{equation*}
     (\nabla \cdot(\hat{c} \bm{w}), h)_{\Omega_v(t)}
     = (\hat{c} \bm{w} \cdot \bm{n}, h)_{\partial \Omega_v(t)} - (\hat{c} \bm{w}, \nabla h)_{\Omega_v(t)} ,
   \end{equation*}
   while the cross-section average definitions combined with the chain rule yield
   \begin{equation*}
     (\partial_t \hat{c}, h)_{\Omega_v(t)} 
     =  (\partial_t \hat{c}, A \avg{h})_{\Lambda}
     = (\partial_t (A\hat{c}),\avg{h})_{\Lambda} - (\hat{c} \partial_t A,\avg{h})_{\Lambda}  .
    \end{equation*}

   We will derive equivalent expressions for the two terms on the right hand side. First for the last term, by definition of the area $A$ and \eqref{eq:transformation_determinant}, we have that: 
    \begin{equation*}
      (\hat{c}\partial_t A, \avg{h})_{\Lambda}
      = \int_{\Lambda} \hat{c} \avg{h} \partial_t \left( \int_{\Theta(t) } 1  \right)
      =  \int_{\Lambda} \int_{\partial \Theta(t)} \hat{c} \avg{h} \bm{w}\cdot \bm{n}
      = ( \hat{c} \bm{w}\cdot \bm{n},  \avg{h} )_{\partial \Omega_v(t)}. 
    \end{equation*}
    Second, we will address the former term in combination with other terms from~\eqref{eq:dual_interm_0}. To this end, denote by $\hat{\bm{u}}_v = (\langle u_{v,s} \rangle, 0 , 0)$. Note  by the definition of $a_{\mathrm{ref}}$, the cross-section, and perimeter averages, and by adding and subtracting, that, 
    \begin{equation}
      \begin{aligned}
          \label{eq:ref_reduced_cc}
      a_{\mathrm{ref}}(\cc,& h) - (\cc \bm{w}, \nabla h)_{\Omega_v(t)} 
      = (D_v \partial_s \cc , \partial_s h)_{\Omega_v(t)} + (\xi \cc, h)_{\Gamma(t)} - (\bm{u}_v \cc, \nabla h)_{\Omega_v(t)} \\ 
      &= (D_v A \partial_s \cc, \avg{\partial_s h})_{\Lambda} + (\xi P \cc, \overline{h})_{\Lambda} - (\bm{u}_v \cc, \nabla{h})_{\Omega_v(t)} \\ 
      &= (D_v A \partial_s \cc, \avg{\partial_s h})_{\Lambda} + (\xi P \cc, \overline{h})_{\Lambda} - ((\bm{u}_v - \hat{\bm{u}}_v) \cc, \nabla{h})_{\Omega_v(t)}  -  (A \langle u_{v,s} \rangle \cc , \avg{\partial_s h})_{\Lambda}. 
    \end{aligned}
    \end{equation}
   We proceed by returning to the weak formulation of the coupled 3D-1D problem~\eqref{eq:coupled_pb_2} with $\hat{v} = \avg{h} \in H^1_{A}(\Lambda)$. We now invoke the assumption that $w_c = 1$; then $\vtwo = 0$ and $\overline{w_c} = 1$.
   Let  $\avgus = \avg{u_{v,s}}$. Then, \eqref{eq:coupled_pb_2}, after combining time-integration terms, gives that
   \begin{equation}
     ( \partial_t (A\cc),  \avg{h})_{\Lambda} + ( D_v  A\partial_s \cc, \partial_s \avg{h})_{\Lambda}  -   (A \avgus \cc, \partial_s  \avg{h})_{\Lambda} + (\xi P(\cc - \bar{c}), \avg{h})_{\Lambda}= (A\avg{f_v}, \avg{h})_{\Lambda}.
     \label{eq:rewriting_1d_weak}
   \end{equation}
   Observe that \footnote{With Leibniz integration rule, we have 
     \begin{align*} 
       \partial_s (A \avg{h}) &= \partial_s \left(\int_{R_1}^{R_2} \int_{0}^{2\pi} h r \dr\dtheta\right)  = \int_{R_1}^{R_2} \int_{0}^{2\pi} \partial_s h r  + \int_{0}^{2\pi} ( h(R_2) R_2 \partial_s R_2 - h(R_1) R_1 \partial_s R_1) \\ 
       & = A \avg{\partial_s h} + \int_{\partial  \Theta_2} h \partial_s R_2 -  \int_{\partial \Theta_1} h \partial_s R_1
     \end{align*}
   }
   \begin{equation} 
     A \partial_s \avg{h} = \partial_s(A \avg{h}) - \avg{h} \partial_s A =  A \avg{\partial_s h} +  \int_{\partial \Theta_2} (h - \avg{h}) \partial_s R_2 - \int_{\partial \Theta_1} (h - \avg{h}) \partial_s R_1 .
     \label{eq:product_rule_avg}
   \end{equation}
   Using~\eqref{eq:product_rule_avg} and  \eqref{eq:ref_reduced_cc} in \eqref{eq:rewriting_1d_weak}, we obtain: 
   \begin{equation}
     (\partial_t(A\cc), \avg{h})_{\Lambda} + a_{\mathrm{ref}}(\cc, h) - (\cc \bm{w}, \nabla h)_{\Omega_v(t)} 
     = - ((\bm{u}_v - \hat{\bm{u}}_v) \cc, \nabla{h})_{\Omega_v(t)} +  \ell(h) - \ell_1(h), 
   \end{equation}
   where we have introduced the short-hand
   \begin{align*}
     \ell(h) & =  \int_\Lambda  \xi P  \overline{c}\,  \overline{h} + \int_\Lambda A \avg{f_v} \avg{h}
     =   (\xi \overline{c}, h)_\Gamma + (\avg{f_v},h)_{\Omega_v}, \\ 
     \ell_1(h) & = \int_{\Lambda}(-D_v  \partial_s \cc +  \cc \avgus) \left(  \int_{\partial \Theta_1}  (h - \avg{h}) \partial_s R_1 - \int_{\partial \Theta_2} (h - \avg{h}) \partial_s R_2  \right)  - (\xi (\avg{h} - \bar{h}) , \overline{c}- \cc )_{\Gamma(t)}.
     \end{align*}
   Collecting all the above expressions in \eqref{eq:dual_interm_0} yields:
   \begin{multline}
     \label{eq:equation_with_terms_to_bound}
     \int_{0}^T \|e\|_{L^2(\Omega_v(t))}^2 
     = \int_0^T ( (\bm{u}_v - \hat{\bm{u}}_v)\cc, \nabla h )_{\Omega_v(t)} 
     - \int_0^T (\hat{c} \bm{w} \cdot \bm{n}, h - \avg{h})_{\partial \Omega_v (t)}  \\ +  \int_0^T \left ( \ell_{\mathrm{ref}} (h) - \ell(h) \right )
     + \int_0^T \ell_1(h) + ( e(0), h(0))_{\Omega_v(0)} 
     := \sum_{i=1}^5 W_i.  
     \end{multline}
   
     \emph{Step III.}  We now bound each term $W_i$ ($i = 1, \dots,
     5$) on the right-hand side of
     \eqref{eq:equation_with_terms_to_bound}. For brevity, we omit the
     time-dependence of the domains in the notation in the below.
     For $W_1$, write 
     \begin{equation*}
      W_1 = \int_0^T ( (0, u_{v,r}, u_{v,\theta})\cc, \nabla h )_{\Omega_v}+ \int_0^T (u_{v,s} - \avg{u_{v,s}}) \cc \partial_s  h = W_{1,1} + W_{1,2}.
     \end{equation*}
     An application of H\"{o}lder's inequality yields
    \begin{equation}
      W_{1,1}
      \equiv \int_0^T ( (0, u_{v,r}, u_{v,\theta})\cc, \nabla h )_{\Omega_v}
      \leq \int_0^T  ( \|u_{v,r}\|_{L^{\infty}(\Omega_v)} + \|u_{v,\theta}\|_{L^{\infty}(\Omega_v)})\|\cc\|_{L^2(\Omega_v)} \|h\|_{H^1(\Omega_v)} . \nonumber
    \end{equation}
For $W_{1,2}$, with H\"{o}lder's and Poincar{\'e}'s inequality~\eqref{eq:Poincare_inequality}, we have that 
   \begin{align*}
     W_{1,2}
     &= \int_0^T \int_{\Lambda} \cc \int_{\Theta} (u_{v,s} - \avg{u_{v,s}}) \partial_s h \leq \int_0^T \int_{\Lambda}|\cc| \| u_{v,s} - \avg{u_{v,s}}\|_{L^2(\Theta)} \|\partial_s h\|_{L^2(\Theta)} \\
     & \leq K_p \int_0^T \int_{\Lambda} |\cc| \epsilon(s,t) \|\nabla u_{v,s}\|_{L^2(\Theta)} \|\partial_s h\|_{L^2(\Theta)} \\
     &= K_p \int_0^T \int_{\Lambda} \epsilon(s, t) P^{-1/2}\|\cc\|_{L^2( \partial \Theta_2)}\|\nabla u_{v,s}\|_{L^2(\Theta)} \|\partial_s h\|_{L^2(\Theta)}. 
   \end{align*}
   Thus, with H\"{o}lder's inequality and the assumption that the vessel area and outer perimeter are both bounded in terms of $\epsilon$ but with (implicit) inequality constants independent of $\Omega_v, \Omega_s$~\eqref{eq:relation_area_diam}, 
   \begin{align}
   W_{1,2} \, \lesssim \, K_p \int_0^T  \epsilon_{\max}^{1/2} \left ( \max_{s \in \Lambda} \|\nabla u_{v,s}\|_{L^{2}(\Theta)} \right )\|\cc\|_{L^2(\Gamma)}\|h\|_{H^1(\Omega_v)}.
   \end{align}
   Hence, with H\"{o}lder's inequality again, we obtain the following for $W_1$: 
   \begin{multline}
     \label{eq:bound_first_term}
 W_1 =   W_{1,1} + W_{1,2}
     \, \lesssim \,  K_p \,\, \epsilon_{\max}^{1/2} \max_{s \in \Lambda, t \in [0,T]}\|u_{v,s}\|_{H^1(\Theta)} \|\cc\|_{L^2(0,T;L^2(\Gamma))} \|h\|_{L^2(0,T;H^1(\Omega_v))} 
    \\   + ( \|u_{v,r} \|_{L^{\infty}(0,T; L^{\infty}(\Omega_v))} + \| u_{v,\theta}\|_{L^{\infty}(0,T; L^{\infty}(\Omega_v))})
 \|\cc\|_{L^2(0,T;L^2(\Omega_v))} \|h\|_{L^2(0,T;H^1(\Omega_v))}  . 
   \end{multline}

   Continuing, we bound $W_2$ and $W_4$ by first obtaining a bound on $\|\overline{v} - \avg{v}\|_{L^2(\Gamma)}$ for any $v \in H^1(\Omega_v)$. First note that
    \begin{multline*}
      \| \avg{v} - \overline{v} \|_{L^2(\partial \Theta)}^2  
      = \int_{\partial \Theta} ( \avg{v} - \overline{v}) ( \avg{v} - \overline{v})
      = \int_{ \partial \Theta} ( \avg{v} - v  ) ( \avg{v} - \overline{v}) + \int_{\partial \Theta} ( v - \overline{v}) ( \avg{v} - \overline{v}) \\
      = \int_{\partial \Theta} ( \avg{v} - v  ) ( \avg{v} - \overline{v})  
      \leq \| \avg{v} - v \|_{L^2(\partial \Theta)} \| \avg{v} - \overline{v} \|_{L^2(\partial \Theta)} .
    \end{multline*}    
    Using this observation, the trace inequality~\eqref{eq:trace_estimate_general}, and Poincare's inequality \eqref{eq:Poincare_inequality}, we have that for any $v \in H^1(\Theta)$
    \begin{multline}
      \label{eq:trace_then_poincare}
      \|\avg{v} - \bar{v}\|_{L^2(\partial \Theta)} \leq \|\avg{v} - v\|_{L^2(\partial \Theta)} \\
      \leq K_{\mathrm{tr}} \left ( \epsilon(s,t)^{-1/2} \|\avg{v} - v\|_{L^2(\Theta)} + \epsilon (s,t)^{1/2} \|\nabla v\|_{L^2(\Theta)} \right )
      \leq K_{\mathrm{tr}}( K_p + 1) \epsilon_{\max}^{1/2} \|\nabla v \|_{L^2(\Theta)}.
    \end{multline} 
    With Cauchy--Schwarz inequality, we have that: 
    \begin{equation*}
      \begin{split}
        &W_2 + W_4
        = \int_0^T (\hat{c} \bm{w} \cdot \bm{n}, h - \avg{h})_{\partial \Omega_v} + \int_0^T \ell_1(h) \\
        &\leq \int_0^T  (\|\bm{w}\|_{L^{\infty}(\Omega)}\|\cc\|_{L^2(\partial \Omega_v)} + \epsilon_s \|D_v \partial_s \cc + \cc \avgus\|_{L^2(\partial \Omega_v)}  + \| \xi (\cc-\overline{c})\|_{L^2(\Gamma)} )  \| \avg{h} - h\|_{L^2(\partial \Omega_v )}  \\ 
      & \leq  K_{\mathrm{tr}}(K_p + 1)  \epsilon_{\max}^{1/2} \int_0^T (\|\bm{w}\|_{L^{\infty}(\Omega)}\|\cc\|_{L^2(\partial \Omega_v)} + \epsilon_s \|D_v \partial_s \cc + \cc \avgus\|_{L^2(\partial \Omega_v)}  + \|\xi(\cc - \overline{c})\|_{L^2(\Gamma)}) \|h\|_{H^1(\Omega_v)}  \\ 
      & \lesssim K_{\mathrm{tr}}(K_p + 1) (\epsilon_{\max}^{1/2} \|\bm{w}\|_{L^{\infty}(0,T;L^{\infty}
      (\Omega))} \|\cc\|_{L^2(0,T;L^2(\partial \Omega_v))} + \epsilon_s \|D_v \partial_s \cc + \cc \avgus\|_{L^2(0,T;L^2(\Omega_v))} \nonumber \\ & \quad  + \epsilon_{\max}^{1/2}   \|\xi (\cc - \overline{c})\|_{L^2(0,T;L^2(\Gamma))} ) \|h\|_{L^2(0,T;H^1(\Omega_v))}. 
      \end{split}
      \end{equation*}
    In the above, we used that $\epsilon_{\max}^{1/2} \|D_v \partial_s \hat{c}+ \hat{c}\avgus\|_{L^2(\partial \Omega_v )} \lesssim \|D_v \partial_s \hat{c}+ \hat{c}\avgus\|_{L^2(\Omega_v)}. $ This follows from the observation that $D_v, \cc,$ and $ \avgus$ are uniform on each cross-section and from \eqref{eq:assump_radii}. 
    
Consider now the definition of $W_3$ in combination with Cauchy-Schwarz:
\begin{equation}
  \begin{split}
    \label{eq:W4}
    W_3
    &= \int_0^T (\xi(c_s - \bar{c}),h)_{\Gamma} + (f_v - \avg{f_v}, h)_{\Omega_v}  \\  
    &\leq \int_0^T \|\xi^{1/2} (c_s - \overline{c})\|_{L^2(\Gamma)}\|\xi^{1/2}h\|_{L^2(\Gamma)}
    + \|f_v -\avg{f_v}\|_{L^2(\Omega_v)} \|h\|_{L^2(\Omega_v)}  .
  \end{split}
\end{equation}
For the first integrand term of the previous line, we may use the observation that $\|\overline{c}\|_{L^2(\Gamma)} = \|\overline{c}\|_{L^2_P(\Lambda)} \leq \|c\|_{L^2(\Gamma)}$ and the trace Lemma \ref{lemma:trace_ineq_omegas} over $\Omega_s$. 
\begin{equation*}
  \|\xi^{1/2} (c_s - \overline{c})\|_{L^2(\Gamma)} \leq K_\Gamma (\epsilon_{\max} |\ln \epsilon_{\max}|)^{1/2} \|\xi\|^{1/2}_{L^{\infty}(\Gamma)}(\|c_s\|_{H^1(\Omega_s)} + \|c \|_{H^1(\Omega_s)}). 
\end{equation*}
For the last integrand in \eqref{eq:W4}, we use the Poincar\'{e} inequality~\eqref{eq:Poincare_inequality}. Combining with Hölder's inequality, we obtain 
\begin{multline*}
  W_3 \leq K_\Gamma (\epsilon_{\max} |\ln \epsilon_{\max}|)^{1/2} \|\xi\|^{1/2}_{L^{\infty}(\Gamma)} \left ( \|c_s\|_{L^2(0,T;H^1(\Omega_s))} + \|c\|_{L^2(0,T;H^1(\Omega))} \right ) \|\xi^{1/2} h\|_{L^2(0, T; L^2(\Gamma))}  \\ + K_p \epsilon_{\max} \|f_v\|_{L^2(0,T;H^1(\Omega_v))} \|h\|_{L^2(0,T;L^2(\Omega_v))}. 
\end{multline*}
Alternatively, if the sharper bound~\eqref{eq:epsilon:better} holds, we then first use the triangle inequality for bounding the first term in  \eqref{eq:W4}:
\begin{equation*}
  \|\xi^{1/2} (c_s - \overline{c})\|_{L^2(\Gamma)}
  \leq \|\xi^{1/2} (c_s - \overline{c_s})\|_{L^2(\Gamma) } + \|\xi^{1/2} (\overline{c_s} - \overline{c}) \|_{L^2(\Gamma)} ,
\end{equation*}
and then a Stekloff-type inequality along with the boundedness of the extension operator $\mathcal{E}$~\eqref{eq:extension_operator}, giving: 
\begin{equation*}
  \|c_s - \overline{c_s}\|^2_{L^2(\Gamma)}
  = \int_{\Lambda} \|\mathcal{E}c_s - \overline{\mathcal{E}c_s}\|^2_{L^2(\partial \Theta_2)}
  \leq K_{\mathrm{st}}\epsilon_{\max}  \int_{\Lambda} \|\nabla \mathcal{E} c_s\|^2_{L^2(\Theta_2)}
  \leq  K_{\mathrm{st}}K_{\mathcal{E}} \epsilon_{\max} \| c_s\|^2_{H^1(\Omega_s)} . 
\end{equation*}
Then $W_3$ can instead be bounded by:
\begin{multline*}
  W_3 \lesssim (K_{\mathrm{st}} K_{\mathcal{E}} + 1) \epsilon_{\max}^{1/2}\|c_s\|_{L^2(0,T;H^1(\Omega_s))} \|\xi h\|_{L^2(0, T; L^2(\Gamma))}  \\ + K_p  \epsilon_{\max} \|f_v\|_{L^2(0,T;H^1(\Omega_v(t)))} \|h\|_{L^2(0,T;L^2(\Omega_v(t)))}. 
    \end{multline*}
The term $W_5$ involving the modelling error associated with the initial condition is handled by the Poincar{\'e}
inequality~\eqref{eq:Poincare_inequality},and that $\cc(0) = \avg{c_v(0)}$:
       \begin{multline*}
         \|(c_v - \cc)(0)\|^2_{L^2(\Omega_v(0))} = \int_{\Lambda} \int_{\Theta} (c_v(0)-\cc(0))^2
         \leq  K_p \epsilon_{\max}^2 \int_{\Lambda} \|\nabla c_v(0)\|_{L^2(\Theta)}^2 \\
         = K_p \epsilon_{\max}^2 \|\nabla c_v(0)\|^2_{L^2(\Omega_v(0))}. 
       \end{multline*}
    This implies that 
    \begin{equation}
        W_5 \leq K_p \epsilon_{\max} \|\nabla c_v(0)\|^2_{L^2(\Omega_v(0))} \|h(0)\|_{L^2(\Omega_v(0))}.
    \end{equation}
    Collecting all the above bounds in \eqref{eq:equation_with_terms_to_bound} and using \eqref{eq:regularity_backward_pb} yields the estimate. The proof of the boundedness of $C_1$ and $C_2$ by a constant $K$ independent of  $\epsilon_{\max}$ is given in the Appendix, section \ref{sec:bound_W}. \qedhere
   \end{proof}

\subsection{Model error introduced in the surrounding 3D domain}
\label{subsec:error_3D}
In this subsection, we study the error introduced in the model derivation of the extended transport model (\Cref{subsec:weak_form_3d_derivation}). In particular, we aim to study the difference $(c_s - c)$ between the reference solution $c_s \in L^2(0,T; H_{\partial \Omega}^1(\Omega_s(t)))$ satisfying the weak solute transport equations defined over $\Omega_s(t)$~\eqref{eq:weak_form_3d3d} and the reduced (or perhaps more aptly, extended) solution $c \in L^2(0,T; H_0^1(\Omega))$ satisfying the weak solute transport equations defined over $\Omega$~\eqref{eq:weak_sol_tissue}. Here, we will assume that $\mathcal{E} D_s \in L^{\infty}(0,T;L^{\infty}(\Omega, \R^{3\times 3}))$ with a uniform ellipticity constant $\tilde{\nu} > 0$.   

We start by recalling the relevant equations and that $\tilde{\bm{u}}_s = \bm{u}_s - \bm{w}$, we have that $c_s$ and $c$ satisfy
\begin{multline}
  \langle \dot{c}_s,  \phi \rangle_{H^{-1}(\Omega_s(t))} + \int_{\Omega_s(t)}( \nabla \cdot \bm{w} c_s \phi + D_s \nabla c_s \cdot \nabla \phi 
  -  (\tilde{\bm{u}}_s c_s) \cdot \nabla \phi ) \\
  +  \int_{\Gamma(t)} \xi (c_s - c_v) \phi
  = \int_{\Omega_s(t)} f_s \phi
  \quad \foralls \phi \in H^1_{\partial \Omega} (\Omega_s(t)) ,
  \label{eq:recalling_ref_3d}
\end{multline}
and
\begin{multline}
  \int_{\Omega} \partial_t c  \phi + \int_{\Omega} (\mathcal{E}D_s \nabla c \cdot \nabla \phi - (\mathcal{E}\bm{u}_s c) \cdot \nabla \phi )
  +  \int_{\Gamma(t)} \xi(\overline{c} - \cc )\phi
  = \int_{\Omega} \mathcal{E} f_s \phi 
  \quad \foralls \phi \in H_0^1(\Omega) .
  \label{eq:derived_3dmodel_with_extension}
\end{multline}
In \eqref{eq:derived_3dmodel_with_extension}, the Eulerian derivative is used since now $\partial \Omega$ is independent of $t$.  
In~\eqref{eq:derived_3dmodel_with_extension}, we  used that 
\begin{equation}
  \int_{\Lambda} \xi P (\overline{c} - \hat{c})  \overline{\phi} = \int_{\Gamma(t)} \xi (\overline{c} - \hat{c}) \phi.  
 \end{equation}

As a step on the way towards quantifying $c_s - c$ over the whole domain, we introduce an intermediate solution $c_r$ solving~\eqref{eq:derived_3dmodel_with_extension} but without the coupling terms and aim to bound $(c_s - c_r)$ and $(c_r - c)$. More precisely, let $c_r \in L^{2}(0,T;H^1_0(\Omega))$ with $c_r(0) = c(0)=\mathcal{E} c_s(0)$ solve
\begin{equation}
  \int_{\Omega} \partial_t c_r  \phi + \int_{\Omega} \mathcal{E}D_s \nabla c_r \cdot \nabla \phi - (\mathcal{E}\bm{u}_s c_r) \cdot \nabla \phi = \int_{\Omega} \mathcal{E}f_s \phi
  \label{eq:regular_ignore_coupling}
\end{equation}
for all $t > 0$ and for all $\phi \in H^{1}_0(\Omega)$.
From standard parabolic regularity results, see e.g~\cite[Chapter 7]{evans2009partial}, and from the continuity of the extension operator~\eqref{eq:extension_operator}, we have for a convex  domain $\Omega$ that: 
\begin{multline}
  \|\partial_t c_r\|_{L^2(0,T;L^2(\Omega))} + \|c_r\|_{L^2(0,T;H^2(\Omega))} \leq K ( \|\mathcal{E} f_s\|_{L^2(0,T;L^2(\Omega))} + \|c_r(0)\|_{H^1(\Omega)}) \\  \leq K_r (\|f_s\|_{L^2(0,T;H^1(\Omega_s(t)))} + \|c_s(0)\|_{H^1(\Omega_s(0))}).  \label{eq:standard_parab_reg}
\end{multline}
Here $K_r$ depends on $\mathcal{E}D_s$, $\mathcal{E}\bm{u}_s$, and the final time $T$. 

We proceed by first bounding $c - c_r$ in~\Cref{lemma:difference_bw_c_cr}, and then consider $c_r - c_s$ and $c - c_s$ in \Cref{lemma:model_error_3d}.
  \begin{lemma}[Estimating $c-c_r$]
    \label{lemma:difference_bw_c_cr}
    For $c$ and $c_r$ defined by~\eqref{eq:derived_3dmodel_with_extension} and~\eqref{eq:regular_ignore_coupling} respectively, there holds that
\begin{equation}
  \|c - c_r\|_{L^{\infty}(0,T;L^2(\Omega))} \leq K_{e_1} \epsilon_{\max}^{1/2} |\ln \epsilon_{\max} |^{1/2} \left ( \|c\|_{L^2(0,T;H^1(\Omega))} +  \|\xi^{1/2} \hat{c}\|_{L^{2}(0,T;L^2(\Gamma(t)))} \right ) ,
 \end{equation}
where $\epsilon_{\max}$ is the maximal vessel cross-section diameter as defined by~\eqref{eq:def_epsilon} and $K_{e_1}$ depends on $T$, $\tilde{\nu}^{-1/2}$, and $\bm{u}_s$,  but  not on $\epsilon_{\max}$.
  \end{lemma}
\begin{proof}
  Define $e_1 \equiv c - c_r$. Subtracting \eqref{eq:regular_ignore_coupling} from \eqref{eq:derived_3dmodel_with_extension}, choosing $\phi = e_1$, integrating over time, and using standard arguments, we obtain: 
  \begin{multline*}
    \|e_1(t)\|_{L^2(\Omega)}^2 + \frac{\tilde{\nu}}{2}  \|\nabla  e_1\|_{L^2(0,t;L^2(\Omega))}^2 
    \leq \frac{1}{2 \tilde \nu}\| \mathcal{E} \bm{u}_s\|_{L^{\infty}(0,T;L^{\infty}(\Omega))}^2 \|e_1\|^2_{L^2(0,t;L^2(\Omega))} \\
    + \|\xi^{1/2}(\overline{c} - \hat{c})\|_{L^2(0,t;L^2(\Gamma(t)))} \| \xi^{1/2} e_1\|_{L^2(0,t;L^2(\Gamma(t)))} \equiv L_1 + L_2.
  \end{multline*}
For the last term $L_2$, we use the trace inequality over $\Omega_s$ (\Cref{lemma:trace_ineq_omegas}) since $e_1= c-c_r \in L^2(0,T;H_0^1(\Omega))$ and thus $e_1 \in L^2(0,T;H^1(\Omega_s(t)))$. Along with Young's inequality, we derive   
\begin{align*}
  L_2
  &\leq K_{\Gamma}\epsilon_{\max}^{1/6}\|\xi^{1/2}\|_{L^{\infty}(0,T;L^{\infty}(\Gamma(t)))}\|\xi^{1/2}(\overline{c} - \hat{c})\|_{L^2(0,t;L^2(\Gamma(t)))}  \| e^{1}\|_{L^2(0,t;H^1(\Omega_s(t)))} \\
  &\leq K_{\Gamma}^2 \left(\frac{1}{\tilde \nu} + 1 \right) \|\xi\|_{L^{\infty}(0,T;L^{\infty}(\Gamma(t)))} \|\xi^{1/2}(\overline{c} - \hat{c})\|_{L^2(0,T;L^2(\Gamma(t)))}^2\epsilon_{\max}|\ln \epsilon_{\max}| \\ & \quad + \frac{\tilde{\nu}}{4} \| \nabla e_1\|^2_{L^2(0,t;L^2(\Omega))}  + \frac{1}{4} \|e_1\|_{L^2(0,t;L^2(\Omega))}^2.
\end{align*}
The first term in the last line above can be further bounded as follows: 
\begin{multline}
  \|\xi^{1/2}(\overline{c} - \cc)\|_{L^2(0,T;L^2(\Gamma(t)))} \leq \|\xi^{1/2}\overline{c}\|_{L^2(0,T;L^2(\Gamma(t)))} + \|\xi^{1/2}\cc\|_{L^2(0,T;L^2(\Gamma(t)))} \\
  \leq K_{\Gamma} \epsilon_{\max}^{1/2}|\ln \epsilon_{\max}|^{1/2}\|\xi\|_{L^{\infty}(0,T;L^\infty(\Gamma))}^{1/2}\|c\|_{L^2(0,T;H^1(\Omega))} + \|\xi^{1/2}\cc\|_{L^2(0,T;L^2(\Gamma(t)))}. 
\end{multline}
The above holds by first noting that $\|\overline{c}\|_{L^2(\Gamma(t))} = \|\overline{c}\|_{L^2_P(\Lambda)}$ and then using Jensen's inequality as in \eqref{eq:jensen_ineq} followed by~\Cref{lemma:trace_ineq_omegas} and $\Omega_s \subset \Omega$. With the above and using $\epsilon_{\max} \lesssim 1$, we obtain that 
\begin{align*}
  \|e_1(t)\|_{L^2(\Omega)}^2 + \frac{1}{4} \|D_s^{1/2}\nabla  e_1\|_{L^2(0,t;L^2(\Omega))}^2 
  \lesssim \frac{1}{2}\left(\frac12 + \frac{1}{\tilde{\nu}} \|\mathcal{E} \bm{u}_s\|_{L^{\infty}(0,T;L^{\infty}(\Omega))}^2\right) \|e_1\|^2_{L^2(0,t;L^2(\Omega))} \\
  +  \epsilon_{\max} |\ln \epsilon_{\max}|(  \|c\|_{L^2(0,T;H^1(\Omega))} + \|\xi^{1/2}\cc\|_{L^2(0,T;L^2(\Gamma(t)))})^2.
\end{align*}
With Gr\"onwall's inequality, we can conclude the result.
\end{proof}

\begin{prop}[Model error in the surroundings]
  \label{lemma:model_error_3d}
  Assume 
that $\Omega$ is convex. 
  Let $c_v, c_s$ be the weak solutions of the coupled 3D-3D transport problem~\eqref{eq:weak_form_3d3d}, and $\cc, c$ be the weak solutions of the reduced 3D-1D problem~\eqref{eq:coupled_3d_1d_weak} with $w_c = 1$.
  Then, there holds that 
    \begin{multline}
      \|c_s - c\|_{L^2(0,T;L^2(\Omega_s(t)))}
      \lesssim 
      N_1 \left((1+\|\bm{u}_s\|_{L^{\infty}(0,T;H^2(\Omega_s(t)))}) \epsilon_{\max}^{2/3}  + \epsilon_{\max} |\ln \epsilon_{\max}| \right) \\ + N_2 (\epsilon_{\max} |\ln \epsilon_{\max}|)^{1/2}. 
      \label{eq:bounding_error_3dmodel}
\end{multline}  
Here, $N_1$ and $N_2$ are given by: 
\begin{align*}
  N_1 &= \|f_s\|_{L^2(0,T;H^1(\Omega_s(t)))} + \|c_s(0)\|_{H^1(\Omega_s(0))},   \\  
  N_2 &=  \|\xi^{1/2} c_v\|_{L^2(0,T;L^2(\Gamma(t)))} + \|c\|_{L^2(0,T;H^1(\Omega)) } + \|\xi^{1/2} \hat{c}\|_{L^2(0,T;L^2(\Gamma(t)))}.
 \end{align*}
In addition, $N_2$ is bounded independently of $\epsilon_{\max}$. 
\end{prop}
\begin{proof}
  Considering Lemma \ref{lemma:difference_bw_c_cr}, it suffices to estimate $\| c_s - c_r\|_{L^2(0,T;L^2(\Omega_s(t)))}$  as the final result  follows by the triangle inequality. The derivation also follows by duality arguments. 
  Define $\psi$ as the solution of the following backward-in-time problem: find $\psi \in L^2(0,T;H_{\partial \Omega}^1(\Omega_s(t)))$ with $\dot{\psi} \in  L^2(0,T;H^{-1}(\Omega_s(t)))$ and $\psi(T) = 0$ in $ \Omega_s(T)$ such that for a.e.~$t$ in $(0,T)$ and for  all  $ v \in H_{\partial \Omega}^1(\Omega_s(t))$:
  \begin{equation}
- \langle \dot{\psi}, v \rangle_{H^{-1}(\Omega_s(t))} + (D_s \nabla \psi, \nabla v)_{\Omega_s(t)} + (\xi \psi, \phi)_{\Gamma(t)} - (\tilde{\bm{u}}_s\cdot \nabla \psi,  v)_{\Omega_s(t)} = (c_s - c_r,v)_{\Omega_s(t)}. \label{eq:back_time_dual_2}
\end{equation}
Then, using similar arguments as in Lemma~\ref{lemma:backward_parabolic}, we have 
\begin{multline}
    \|\psi\|_{L^{\infty}(0,T;L^2(\Omega_s(t)))} + 
    \nu \|\nabla \psi\|_{L^2(0,T;L^2(\Omega_s(t)))} + \|\xi^{1/2} \psi\|_{L^2(0,T;L^2(\Gamma(t)))} \\
    \leq K \left (1 + \|\nabla \cdot \bm{w}\|_{L^{\infty}(0,T;L^\infty(\Omega))} + \| \tilde{\bm{u}}_s\|_{L^\infty(0,T;L^{\infty}(\Omega_s(t)))} \right ) \|c_s - c_r\|_{L^2(0,T;L^2(\Omega_s(t)))} . \label{eq:regularity_backward_pb_2}
    \end{multline} 
Testing \eqref{eq:back_time_dual_2} with $v = e \equiv c_r - c_s \in H_{\partial \Omega}^1(\Omega_s(t))$ for a.e.~$t$, integrating from $0$ to $T$, using the integration by parts rule \cite[Corollary 2.41]{alphonse2015abstract}, and using that $(c_r-c_s)(0) = 0 $ in $\Omega_s(0)$ and $\psi(T) = 0$ in $\Omega_s(T)$ yield 
\begin{multline*}
  L = \int_0^T \|c_r - c_s\|_{L^2(\Omega_s(t))}^2  = 
    \int_0^T \langle \dot{e}, \psi \rangle_{H^{-1}(\Omega_s(t))} +  (e \psi, \nabla \cdot \bm{w})_{\Omega_s(t)} + 
    (D_s \nabla \psi, \nabla e)_{\Omega_s(t)} \\ + (\xi \psi,e)_{\Gamma(t)} - (e (\bm{u}_s -\bm{w}), \nabla \psi)_{\Omega_s(t)} \dt .
\end{multline*}
Next, we expand $e = c_r - c_s$, replace $\psi$ by $\mathcal{E} \psi$ its extension from $\Omega_s(t)$ to $\Omega$, use the equations for the weak solution $c_s$ recalled in~\eqref{eq:recalling_ref_3d} (with $\bm{u}_s - \bm{w} = \tilde{\bm{u}}_s$), the relation between the material and partial time derivative~\eqref{eq:strong_mater_derivative} in combination with the product rule to find:
\begin{multline*}
  L
  \equiv \int_0^T (\partial_t c_r, \mathcal{E} \psi)_{\Omega_s(t)}
  +  (\mathcal{E} \psi, \nabla \cdot (c_r\bm{w}))_{\Omega_s(t)}
  + (D_s \nabla (\mathcal{E} \psi), \nabla c_r)_{\Omega_s(t)} \\
  - ((\bm{u}_s - \bm{w} ) c_r, \nabla \mathcal{E} \psi)_{\Omega_s(t)}  
  - (f_s, \psi)_{\Omega_s(t)} +  ( \xi(c_r- c_v), \psi)_{\Gamma(t)} \dt .
\end{multline*}
Now, we use the definition of $c_r$~\eqref{eq:regular_ignore_coupling}, expand terms involving $\bm{w}$ and use integration by parts to be left with terms over $\Omega_v$ and $\Gamma$:
\begin{multline*}
  L = \int_0^T - (\partial_t c_r, \mathcal{E} \psi)_{\Omega_v(t)} - (\mathcal{E}D_s \nabla (\mathcal{E} \psi) , \nabla c_r)_{\Omega_v(t)}
  + (\mathcal{E} \bm{u}_s c_r, \nabla \mathcal{E} \psi)_{\Omega_v(t)} + (\mathcal{E} f, \mathcal{E} \psi)_{\Omega_v(t)} \dt \\
  + \int_0^T ( \psi,  c_r \bm{w} \cdot \bm{n})_{\Gamma(t)}  + ( \xi(c_r- c_v), \psi)_{\Gamma(t)} \dt   
  \equiv T_1 + \ldots + T_6.  
\end{multline*}

Our next task is to bound each term $T_i$ for $i =1, \dots, 6$. Hereinafter, we omit writing $t$ for the sake of brevity. To bound $T_1$, we first apply Cauchy-Schwarz inequality to have that 
\begin{equation*}
  T_1 \leq \|\partial_t c_r \|_{L^2(0,T;L^2(\Omega_v))} \|\mathcal{E} \psi \|_{L^2(0,T;L^2(\Omega_v)) }. 
\end{equation*}
With H\"{o}lder's inequality, a Sobolev embedding, and the continuity of the extension operator~\eqref{eq:extension_operator}, we obtain
\begin{multline}
  \|\mathcal{E} \psi\|_{L^2(\Omega_v)}
  \leq |\Omega_v|^{1/3} \|\mathcal{E}\psi\|_{L^6(\Omega_v)}
  \leq |\Omega_v|^{1/3} \|\mathcal{E}\psi\|_{L^6(\Omega)} \\
  \leq K |\Omega_v|^{1/3} \|\mathcal{E}\psi\|_{H^1(\Omega)}
  \leq K |\Omega_v|^{1/3} \|\psi\|_{H^1(\Omega_s)}.
  \label{eq:bounding_extended_psi} \end{multline} 
In the above bound, $K$ depends on $\Omega$ but not on $\Omega_v$. Hence, 
\begin{align}
  T_1 \leq K \max_{t \in [0,T]} |\Omega_v|^{1/3}\|\partial_t c_r \|_{L^2(0,T;L^2(\Omega_v))} \| \psi \|_{L^2(0,T;H^1(\Omega_s)) }.
  \label{eq:bound_A1}
\end{align}

To handle $T_2$, we use a similar approach. Since $c_r \in L^2(0,T;H^2(\Omega))$, $\nabla c_r \in L^2(0,T;H^{1}(\Omega)^3)$, a continuous Sobolev embedding yields: 
\begin{equation}
  \|\nabla c_r \|_{L^q(\Omega)} \leq K \|\nabla c_r\|_{H^{1}(\Omega)},  \quad q \in [1,6].
  \label{eq:embed_Hs}
\end{equation}
Hence, with H\"{o}lder's inequality and the above bound~\eqref{eq:embed_Hs}, we have
\begin{align*}
\|\nabla c_r \|_{L^2(\Omega_v)} \leq |\Omega_v|^{1/3} \|\nabla &c_r\|_{L^{6}(\Omega_v)} \leq |\Omega_v|^{1/3} \|\nabla c_r\|_{L^6(\Omega)} \leq K  |\Omega_v|^{1/3} \|c_r\|_{H^{2}(\Omega)}.
\end{align*} 
Then, with the continuity of $\mathcal{E}$~\eqref{eq:extension_operator}, it follows that
\begin{multline*}
  T_2
  \leq \|\mathcal{E}D_s \nabla \mathcal{E} \psi \|_{L^2(0,T;L^2(\Omega_v))} \|\nabla c_r\|_{L^2(0,T;L^2(\Omega_v))} \\
  \leq K\max_{t \in [0,T]} |\Omega_v|^{1/3} \| \psi \|_{L^2(0,T;H^1(\Omega_s))} \|c_r\|_{L^2(0,T;H^{2}(\Omega))},
\end{multline*}
where $K$ depends on $\mathcal{E}D_s$ and again $\Omega$, but not on $\Omega_v$.

For $T_3$, we again use similar arguments as for $T_1$ cf.~\eqref{eq:bounding_extended_psi} to obtain that 
$$ \|c_r\|_{L^2(\Omega_v)} \leq |\Omega_v|^{1/3} \|c_r\|_{L^{6}(\Omega_v)} \leq|\Omega_v|^{1/3} \|c_r\|_{L^{6}(\Omega)}  \leq K |\Omega_v|^{1/3} \|c_r\|_{H^1(\Omega)}.  $$
Further, by the Sobolev embedding $H^2(\Omega) \subset L^\infty(\Omega)$, 
the following bound holds
\begin{align*}
  T_3 
  & \leq K\max_{t\in [0,T]} |\Omega_v|^{1/3} \|\mathcal{E}\bm{u}_s\|_{L^{\infty}(0,T;L^{\infty}(\Omega))}\|c_r\|_{L^2(0,T;H^1(\Omega))}\| \mathcal{E} \psi \|_{L^2(0,T;H^1(\Omega_v))}\\ 
  &\leq K \max_{t\in [0,T]} |\Omega_v|^{1/3} \|\mathcal{E}\bm{u}_s\|_{L^{\infty}(0,T;H^2(\Omega))}\|c_r\|_{L^2(0,T;H^1(\Omega))}\| \psi \|_{L^2(0,T;H^1(\Omega_s))}  \\ 
  & \leq K \max_{t\in [0,T]}  |\Omega_v|^{1/3}\|\bm{u}_s\|_{L^{\infty}(0,T;H^2(\Omega_s))}\|c_r\|_{L^2(0,T;H^1(\Omega))}\| \psi \|_{L^2(0,T;H^1(\Omega_s))}. 
\end{align*}
With \eqref{eq:bounding_extended_psi}, 
the term $T_4$ is bounded as follows. 
\begin{multline*}
  T_4
  \leq \|\mathcal{E} f_s\|_{L^2(0,T;L^2(\Omega_v))}\|\mathcal{E}\psi\|_{L^2(0,T;L^2(\Omega_v))} \\
  \leq K \max_{t\in [0,T]}  |\Omega_v|^{2/3} \|f_s\|_{L^2(0,T;H^1(\Omega_s))} \| \psi \|_{L^2(0,T;H^1(\Omega_s))}. 
\end{multline*}  
For the remaining $T_5$ and $T_6$, we use Cauchy-Schwarz and the trace inequality over $\Omega_s$ (\Cref{lemma:trace_ineq_omegas}, \eqref{eq:trace_eq_Omegas}) to arrive at
\begin{multline*}
  T_5 \leq \| \bm{w}\|_{L^{\infty}(0,T;L^{\infty}(\Gamma))} \|c_r\|_{L^2(0,T;L^{2}(\Gamma(t)))} \|\psi\|_{L^{2}(0,T;L^2(\Gamma(t)))} \\
  \leq K \epsilon_{\max} |\ln \epsilon_{\max}|\|c_r\|_{L^2(0,T;H^1(\Omega))} \|\psi\|_{L^{2}(0,T;H^1(\Omega_s(t)))},
\end{multline*}
and
\begin{align*}
  T_6
  &\leq K \epsilon^{1/2}_{\max}|\ln \epsilon_{\max}|^{1/2}\|\xi^{1/2}(c_r -c_v)\|_{L^2(0,T;L^2(\Gamma(t)))} \|\psi\|_{L^2(0,T;H^1(\Omega_s(t)))} \\
  &\leq K  \epsilon_{\max}^{1/2}|\ln \epsilon_{\max}|^{1/2}( \|\xi^{1/2}c_v\|_{L^2(0,T;L^2(\Gamma(t)))} + \epsilon_{\max}^{1/2} |\ln \epsilon_{\max}|^{1/2}\|c_r\|_{L^2(0,T;H^1(\Omega_s(t)))}) \|\psi\|_{L^2(0,T;H^1(\Omega_s(t)))} .
\end{align*}
Now, having bounded $T_1, \dots, T_6$, we use the regularity bound of the backward in time problem ~\eqref{eq:regularity_backward_pb_2} and  that $|\Omega_v(t)| \leq K \epsilon_{\max}^2$, to obtain
\begin{multline*} 
  \|e_2\|_{L^2(0,T;L^2(\Omega))} 
  \leq  K \epsilon_{\max} |\ln \epsilon_{\max}| \|c_r\|_{L^2(0,T;H^1(\Omega))}
  + K (\epsilon_{\max} |\ln \epsilon_{\max}| )^{1/2} \|\xi^{1/2}  c_v\|_{L^2(0,T;L^2(\Gamma(t)))} \\
  +  K \epsilon_{\max}^{2/3} \left( \|\partial_t c_r \|_{L^2(0,T;L^2(\Omega_v))} + (\|\bm{u}_s\|_{L^\infty(0,T;H^2(\Omega_s))}+1) \|c_r\|_{L^2(0,T;H^2(\Omega))} + \|f_s\|_{L^2(0,T;H^1(\Omega_s))} \right).
  \end{multline*}
The proof is concluded by the triangle inequality, \eqref{eq:standard_parab_reg}, and Lemma \ref{lemma:difference_bw_c_cr}. 
The boundedness of $N_2$ is shown in Appendix~\ref{sec:bound_W}. 
\end{proof}

\section{Numerical results}

In this section, we consider two numerical examples to demonstrate the analysis presented in the previous sections. The two examples correspond to the 3D-1D model of \Cref{sec:coupled_3d_1d} and to the 3D-1D-1D model of \Cref{sec:3D1D1D}. Our implementation uses the FEniCS finite element framework \cite{alnaes2015fenics} and the $(\mathrm{FEniCS})_{\mathrm{ii}}$ module \cite{kuchta2020assembly}.

\subsection{A coupled 3D-1D solute transport finite element example} 
\label{sec:example1}

\begin{figure}[H]
  \begin{minipage}{0.48\linewidth}
    \centering
    \includegraphics[width=\textwidth]{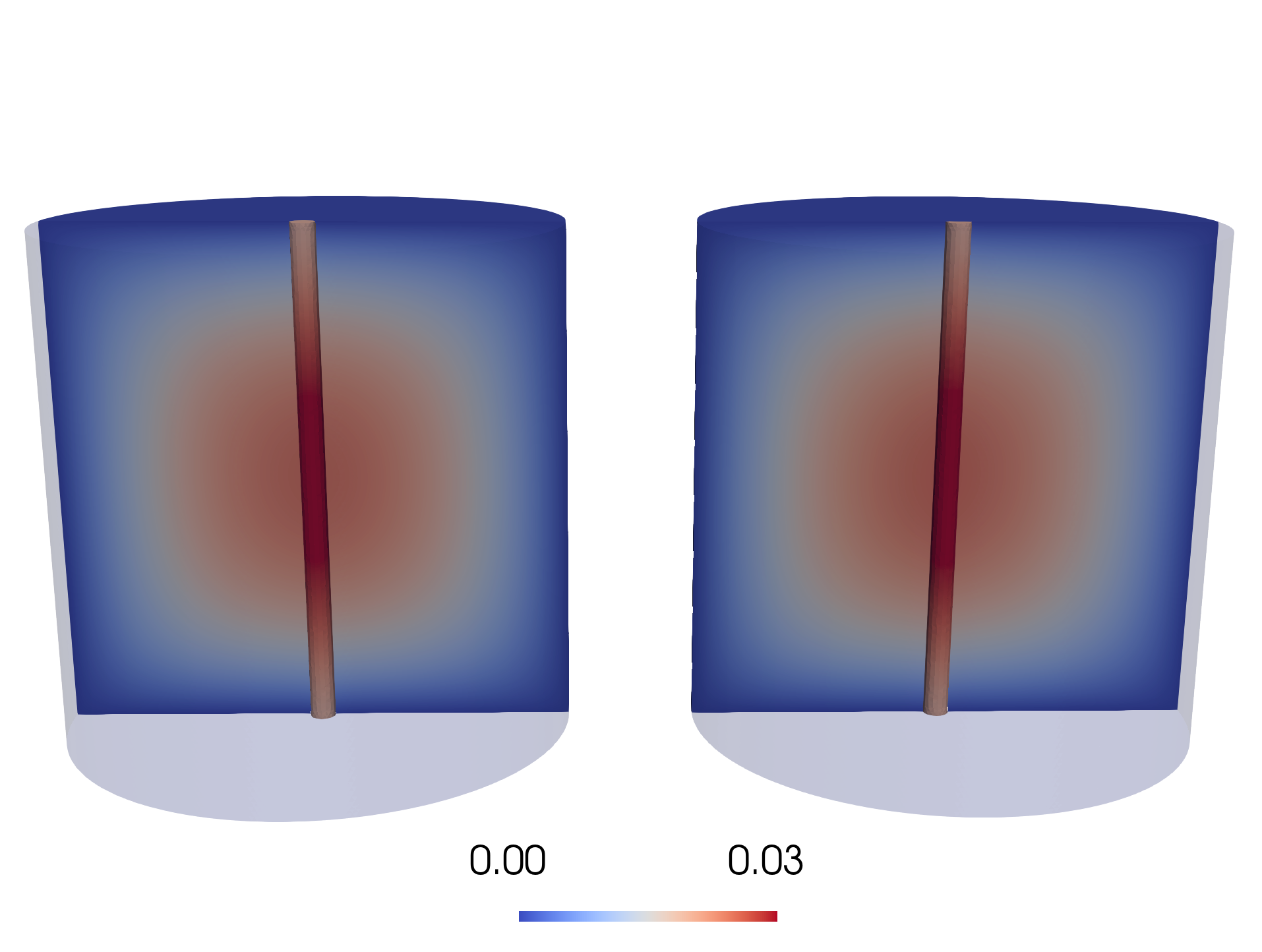}
  \end{minipage}
  \hfil
  \begin{minipage}{0.48\linewidth}
    \centering
    \includegraphics[width=\textwidth]{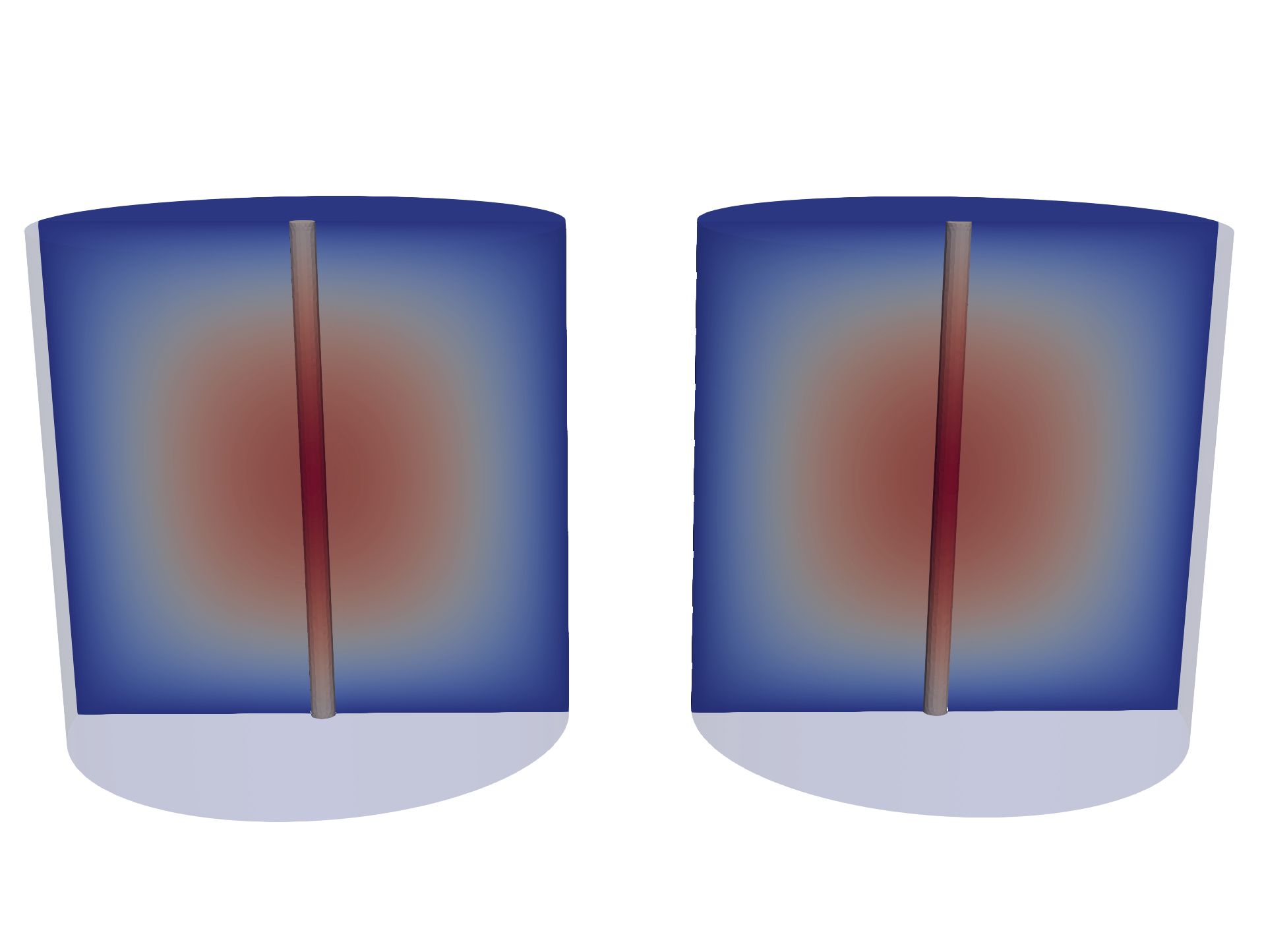}
  \end{minipage}
  \vspace{-10pt}    
  \caption{Plot of the numerical solutions for the first example with $R = 0.0125$. The 3D-3D solution $(c_{s,\tau h}, c_{v,\tau h})$ is plotted next to the 3D-1D model $(c_{\tau,h}, E \cc_{\tau h})$ with the 1D solution extended to the inner cylinder. The outer cylinder is clipped at the plane intersecting the center line $\Lambda$. (Left) Solutions shown at $t= 0.1$. (Right) Solutions shown at $ t = T=0.2$.}
  \label{fig:first_num_eg}
\end{figure}

We let the surrounding domain $\Omega$ (also) take the form of
cylinder with radius $0.5$ and length $L$ containing an inner cylinder
$\Omega_v$ of radius $R = R_2 < 0.5$ with centerline $\Lambda$. Using
a Galerkin finite element method in space with continuous piecewise
linear polynomials defined relative to conforming meshes of $\Omega_s
= \Omega \backslash \Omega_v$, $\Omega_v$ and an implicit Euler
discretization in time with time step $\tau$, we compute approximate
3D-3D solutions $c_{v, \tau h}, c_{s, \tau h}$ of \eqref{eq:3d-3d}. On
the same meshes of $\Omega$ with centerline meshes $\Lambda_h$, we
compute approximate solutions to \eqref{eq:coupled_3d_1d_weak}, again
using continuous piecewise linear finite elements defined relative to
$\Omega$ for $c_{\tau h}$ and relative to $\Lambda_h$ for
$\hat{c}_{\tau h}$ (\Cref{fig:first_num_eg}).  We set $D_v= D_s= \xi = 1$, $f_s= f_v = 0.5$,
$\hat{c}_h^0 = 1.0$, $c^0 = 0.0$, $\bm{u}_v = (0.5,0,0)$ and $\bm{u}_s
= (0.1, 0, 0)$, and $T = 0.2$. 
To numerically explore the modelling error for decreasing radii ($R$,
$\epsilon_{\max} \rightarrow 0$), we consider a series of experiments
with different radii $R \in \{0.2, 0.1, 0.05, 0.025, 0.0125\}$, a relatively small, fixed mesh size
$(h_{\min}, h_{\max})\vert_{\Omega_v} = (0.009,0.014)$ and $(h_{\min}, h_{\max})\vert_{\Omega_s}= (0.01, 0.024)$, and small, fixed time step $\tau = 0.01$. In practice, we compute the
discrepancy between the approximate solutions in the 3D and 1D vessels:
\begin{equation*}
  \|c_{v,\tau h}(T)- E \cc_{\tau h}(T)\|_{L^2(\Omega_v)}  \approx \|c_{v}(T)- E \cc (T)\|_{L^2(\Omega_v)} 
\end{equation*}
as a proxy for the modelling error while noting that the computed error includes both the spatio-temporal approximation errors as well as modelling errors:
\begin{multline*}
  \|c_{v}(T) - \cc(T)\|_{L^2(\Omega_v)} \\
  \leq \| c_v (T) - c_{v,\tau h}(T) \|_{L^2(\Omega_v)} + 
  \|c_{v,\tau h} (T)- \cc_{\tau h} (T)\|_{L^2(\Omega_v)} +
  \|\cc(T)- \cc_{\tau h}(T)\|_{L^2(\Omega_v)} .
\end{multline*}
We here thus presume that with the choice of small mesh size and time step, the approximation errors are negligible compared to the modelling error.

Table \ref{table:rates_vessel_tissue} shows the computed $L^2$ norms in
$\Omega_v$ and $\Omega_s$ along with normalized norms and the
corresponding rates. We observe that the errors decrease with
decreasing $R$ until the radius and mesh size become of comparable
size, and that the modelling error in the surroundings continues to decrease
even when the modelling error in the vessel stagnates.
  \begin{table}[H]
      \centering
    \begin{tabular}{c|cc|cc|cc} 
      \toprule
      $R$ & $E_v$ & rate  & $\tilde{E}_v$ & rate & $E_s$ & rate \\  
      \midrule
      0.1 & 4.405e-04 & -  & 2.485e-03  & - & 4.426e-04 & -  \\ 
          0.05 & 5.244e-05 & 3.07 & 5.918e-04 & 2.07 & 1.375e-04 & 1.69\\ 
          0.025 & 1.394e-05 & 1.91 & 3.146e-04 & 0.91 & 3.636e-05 & 1.92 \\ 
          0.0125 & 7.806e-06 & 0.84 & 3.523e-04 & -0.16 & 9.738e-06 & 1.90 \\ 
          \bottomrule 
    \end{tabular}
    \vspace{0.5em}
    \caption{Numerical example 1: Model errors $E_v = \|c_{v,\tau h} (T)- \cc_{\tau ,
        h}(T)\|_{L^2(\Omega_v)}$ and $\tilde{E}_v = |\Omega_v|^{-1/2}
      \|c_{v,\tau h} (T)-\cc_{\tau, h} (T)\|_{L^2(\Omega_v)}$ in a
      vessel $\Omega_v$ of varying radius $R$, and in the surroundings
      $E_s = \|c_{s,\tau h} (T)- c_{h}(T)\|_{L^2(\Omega_s)}$. }
    \label{table:rates_vessel_tissue} \vspace{-5em}
  \end{table}
\subsection{A coupled 3D-1D-1D solute transport example} 
As a second example, we consider solutions to the coupled 3D-1D-1D models of solute transport and the corresponding 3D-3D-3D model set up in the blood vessel, $\Omega_v$, the perivascular domain $\Omega_p$, and the tissue $\Omega_s$.  We also use backward Euler and continuous linear finite element methods to solve \eqref{eq:coupled_3d3d1d_1}-\eqref{eq:coupled_3d3d1d_3} with solutions denoted by $(c_{\tau h}, \cc_{p,\tau h}, \cc_{v,\tau h})$, and the corresponding 3D-3D-3D model with solutions denoted by $(c_{s, \tau h}, c_{p,\tau h}, c_{v,\tau h})$. We set $\Omega_s = (-1,1)\times (-1,1)\times (-0.5,0.5)$, $\Omega_v$ be a cylinder of radius $R_1$ with centerline $x = 0$, $y=0$, $\Omega_p$ be the annular cylinder around $\Omega_v$ with outer radius $R_2 = 2 R_1$ .  We vary $R_1$ and compute the $L^2$ error between the 3D solutions $c_{i, \tau h}$ and the reduced 1D solutions $\hat{c}_{i,\tau h }$ for $i \in \{p,v\}$.
We keep $\tau = 0.01$ and $T= 0.1$,  $D_v=  D_p = D_s= \xi_v = \xi_p =  1$,  $f_s= f_v = f_p= 0.5$, $c_{v,\tau h} = \hat{c}_{v,\tau h} (0) = 1.0$, $ c_{p,\tau h} = \hat{c}_{p, \tau h} (0)  =  c_{\tau h}(0) = c_{s,\tau h}(0) = 0$,  $\bm{u}_v = (0.5,0,0)$, $\bm{u}_p = (0.1, 0, 0)$, and $\bm{u}_s = (0.05,0,0)$.
For the mesh--size in the variuous domains, we have $(h_{\min}, h_{\max})\vert_{\Omega_v} = (0.011, 0.019)$, $(h_{\min}, h_{\max})\vert_{\Omega_p} = (0.011, 0.024)$, and $(h_{\min}, h_{\max})\vert_{\Omega_s}= (0.017,0.043)$. 

Tables~\ref{table:rates_vessel_pvs} and~\ref{table:rates_pvs_surr} show the computed $L^2$ norm in $\Omega_v$, $\Omega_p$ and $\Omega_s$ along with a normalized norm and the corresponding rates. We observe that the modelling errors all decrease for decreasing radii, though in a non-uniform manner and with uneven rates. The modeling error in the surroundings decreases robustly at rates between 1 and 2. The (non-normalized) modelling error in the vascular domain decreases with similar rates. The modelling error in the perivascular space increases in the first $R$-refinement before decreasing at rates close to $2$. Clearly, further theoretical and numerical studies of the interplay between the modelling and approximation errors are warranted (though outside the scope of the current study).

\begin{figure}[H]
  \begin{minipage}{0.48\linewidth}
    \centering
    \includegraphics[width=\textwidth]{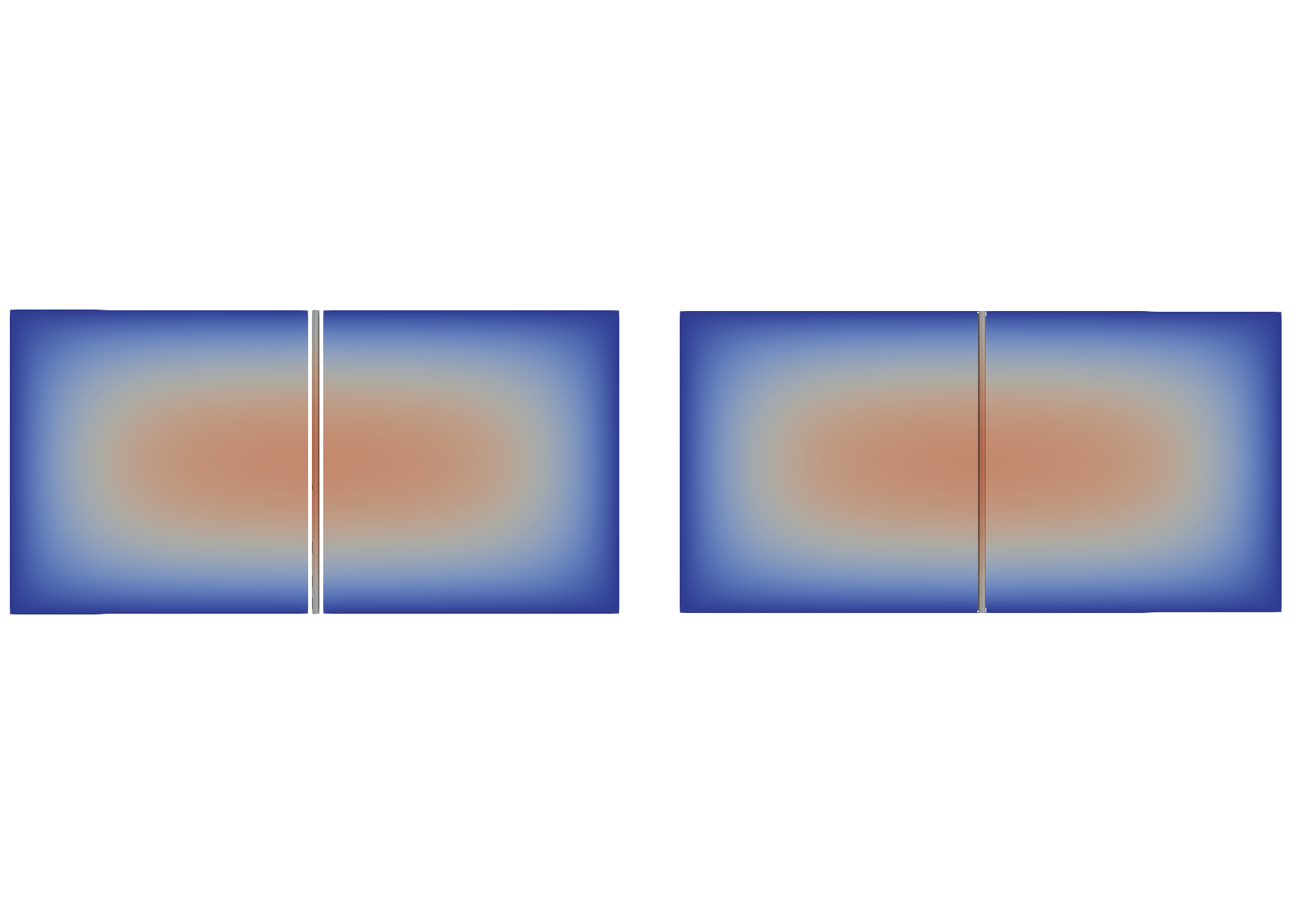}\\
    \vspace{-7em}
    \includegraphics[width=\textwidth]{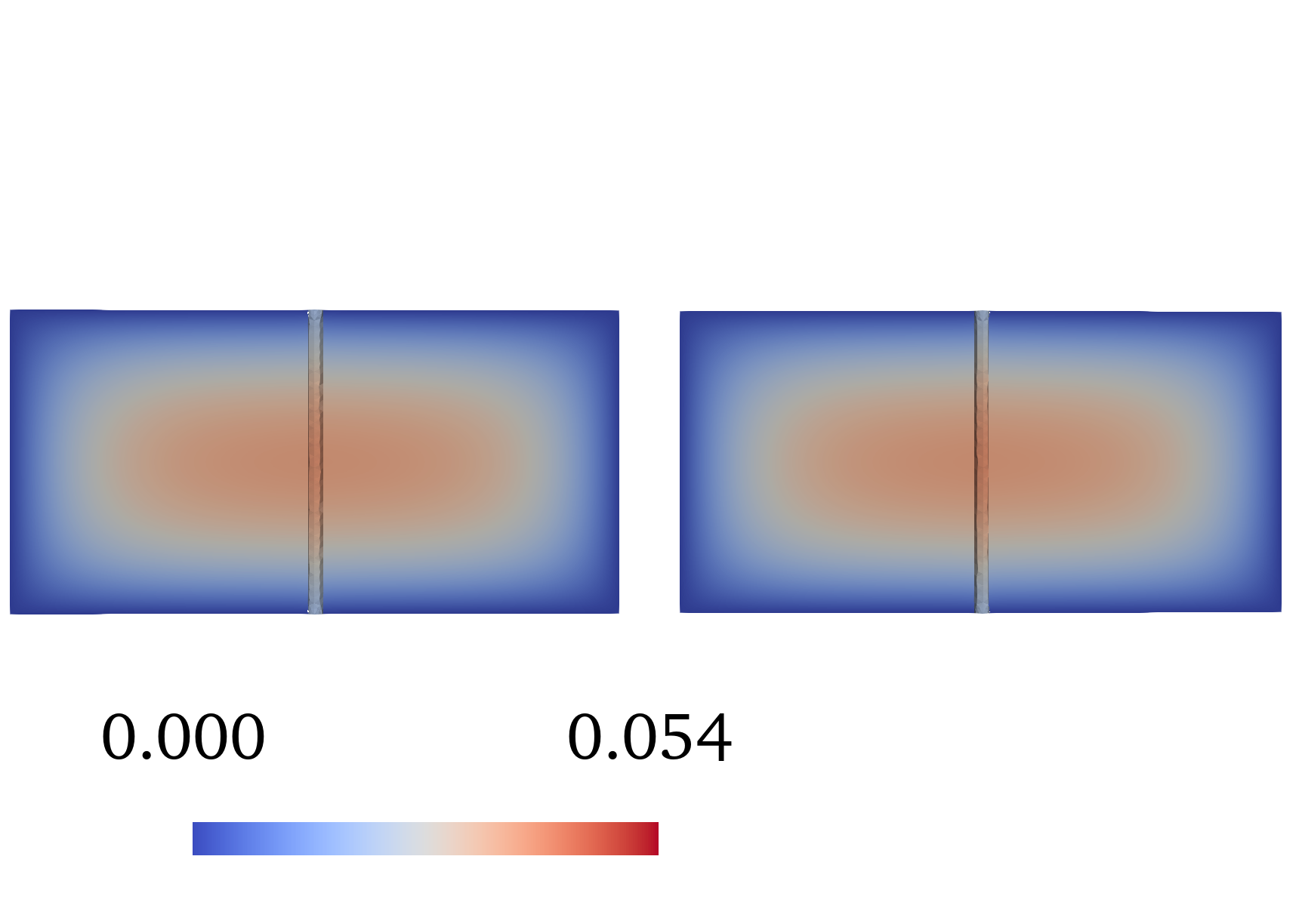}
  \end{minipage}
  \begin{minipage}{0.48\linewidth}
    \centering
    \includegraphics[width=\textwidth]{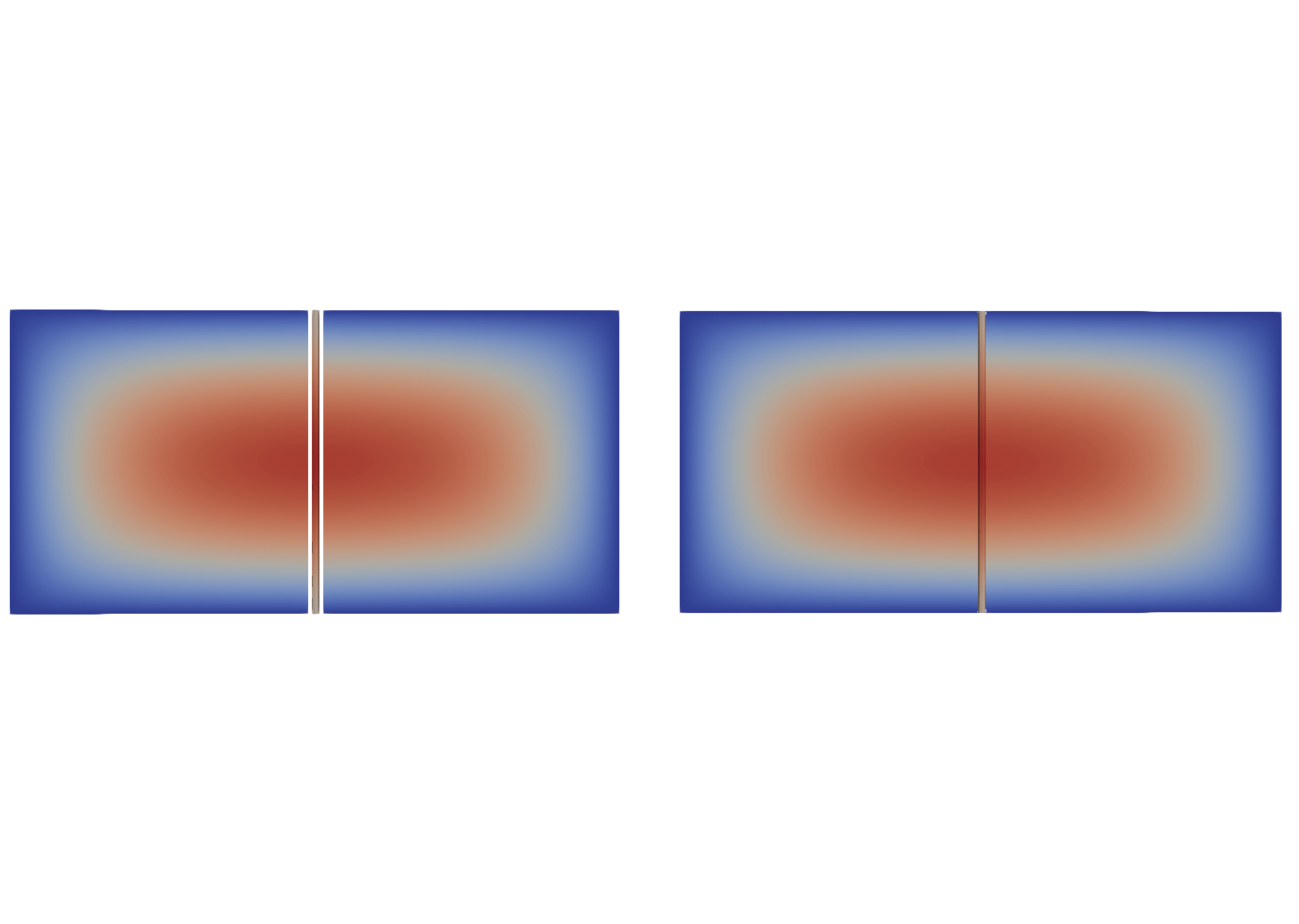} \\ 
    \vspace{-7em}
    \includegraphics[width = \textwidth]{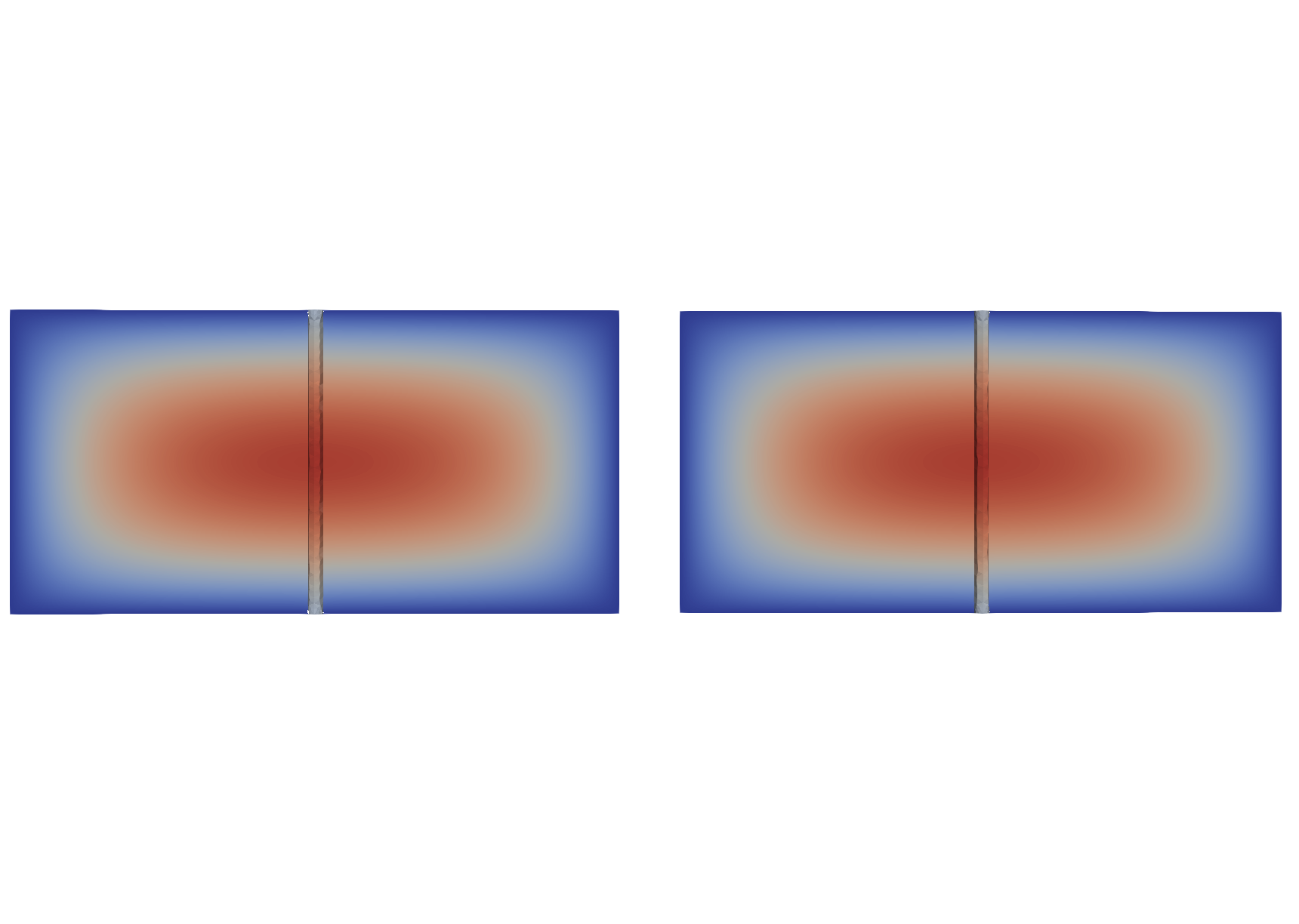} 
  \end{minipage}
  \vspace{-10pt}    
  \caption{Plot of the numerical solutions for the second example with $R = 0.0125$. In each of the four quadrants, the 3D-3D-3D solution $(c_{s, \tau h}, c_{p,\tau h}, c_{v,\tau h})$ is plotted next to the 3D-1D-1D model $(c_{\tau h}, \cc_{p,\tau h}, \cc_{v,\tau h})$  with the 1D solutions extended to their respective cylinder or annulus. Top row shows a slice of $\Omega_s$ and $\Omega$ respectively with the vessel solutions without the PVS domain. The second row shows the PVS solution. (Left) Solutions shown at $t= 0.1$. (Right) Solutions shown at $ t = T=0.2$. }
  \label{fig:second_num_eg}
\end{figure}
\begin{center}
  \begin{table}
  \begin{tabular}{c|cccccccc} 
      \toprule
      $R_1$ & $E_{v,2}$  & rate  & $\tilde{E}_{v,2}$ & rate & $E_{p}$  & rate  & $\tilde{E}_{p}$ & rate   \\
      \midrule
      0.1& 2.195e-03 & -- & 1.239e-02 &-- & 1.698e-03 & -- & 4.790e-03 & --\\ 
      0.05 & 1.902e-04 & 3.53 & 2.146e-03 & 2.53 &1.609e-04 & 3.40 & 9.077e-04 & 2.40   \\
      0.025  &3.243e-05 & 2.55 & 7.319e-04 & 1.55 & 2.950e-05 & 2.48 & 3.329e-04 & 1.48 \\ 
      0.0125 &1.507e-05 & 1.10 & 6.803e-04 & 0.10 & 2.162e-05 & 0.49 & 4.880e-04  & -0.55 
      \\ 
      \bottomrule
  \end{tabular}
  \vspace{0.5em}
  \caption{Numerical example 2: Model errors $E_{v,2} = \|c_{v,\tau h } - E \cc_{v,\tau h} \|_{L^2(\Omega_v)}$, $ \tilde{E}_{v,2} = |\Omega_v|^{-1/2}\|c_{v,\tau h } - E \cc_{v,\tau h} \|_{L^2(\Omega_v)}$,$E_{p} = \|c_{p,\tau h } - E \cc_{p,\tau h} \|_{L^2(\Omega_p)}$, and $ \tilde{E}_{p} = |\Omega_v|^{-1/2}\|c_{p,\tau h } - E \cc_{p,\tau h} \|_{L^2(\Omega_p)}$  of varying radius $R_1$. Note that the smallest $R_1$ is of the order of $h_{\min}$ in $\Omega_v$  }
  \label{table:rates_vessel_pvs}
\end{table}
\end{center}
\begin{center}
\begin{table}[h]
 \centering
\begin{tabular}{c|cccc} 
    \toprule
    $R_2$ & $\|c_{\tau h} (T)- c_{s,\tau  h}(T)\|_{L^2(\Omega_s)}$  & rate  \\  \midrule
    0.2 &9.816e-04 & --   \\
    0.1&1.762e-04 & 2.48  \\ 
    0.05 &6.869e-05 & 1.36  \\
    0.025 & 3.693e-05 & 0.89  \\ 
    \bottomrule
\end{tabular}
  \vspace{0.5em}
\caption{Numerical example 2: Model error in the surrounding domain $\Omega_s$.}
\label{table:rates_pvs_surr}   
\end{table}
\end{center}

\section{Conclusions and outlook}
\label{sec:conclusion}

Understanding solute transport and exchange in the brain vasculature,
perivasculature, and surrounding tissue is critical for unraveling
the brain's delivery and clearance mechanisms. Here, we have presented
a mathematical model for modelling diffusive and convective transport
and exchange in deformable domains, and rigorously analyzed its
modelling characteristics. Future research directions include the
error analysis of conforming and non-conforming finite element
approximations of such models. We easily envision that this framework
can be combined with medical imaging to study brain perivascular
transport and exchange at scale. 
\section*{Acknowledgments}We gratefully acknowledge valuable discussions with Prof.~Barbara Wohlmuth and Dr. Johannes Haubner.


\bibliography{references}

\begin{thebibliography}{10}

\bibitem{abbott2004evidence}
N.~J. Abbott.
\newblock Evidence for bulk flow of brain interstitial fluid: significance for
  physiology and pathology.
\newblock {\em Neurochemistry International}, 45(4):545--552, 2004.

\bibitem{alnaes2015fenics}
M.~Aln{\ae}s, J.~Blechta, J.~Hake, A.~Johansson, B.~Kehlet, A.~Logg,
  C.~Richardson, J.~Ring, M.~E. Rognes, and G.~N. Wells.
\newblock The {FEniCS} project version 1.5.
\newblock {\em Archive of Numerical Software}, 3(100), 2015.

\bibitem{alphonse2015abstract}
A.~Alphonse, C.~M. Elliott, and B.~Stinner.
\newblock An abstract framework for parabolic {PDE}s on evolving spaces.
\newblock {\em Portugaliae Mathematica}, 72(1):1--46, 2015.

\bibitem{arendt2017jl}
W.~Arendt, D.~Dier, and S.~Fackler.
\newblock {JL} {L}ions' problem on maximal regularity.
\newblock {\em Archiv der Mathematik}, 109(1):59--72, 2017.

\bibitem{boron2012medical}
W.~F. Boron and E.~L. Boulpaep.
\newblock {\em Medical Physiology}.
\newblock Elsevier Health Sciences, 2012.

\bibitem{brenner2008mathematical}
S.~C. Brenner and L.~R. Scott.
\newblock {\em The Mathematical Theory of Finite Element Methods}, volume~3.
\newblock Springer, 2008.

\bibitem{brezis2010functional}
H.~Brezis.
\newblock {\em Functional Analysis, Sobolev Spaces and Partial Differential
  Equations}, volume~2.
\newblock Springer, 2010.

\bibitem{brinker2014new}
T.~Brinker, E.~Stopa, J.~Morrison, and P.~Klinge.
\newblock A new look at cerebrospinal fluid circulation.
\newblock {\em Fluids and Barriers of the CNS}, 11(1):1--16, 2014.

\bibitem{vcanic2003mathematical}
S.~{\v{C}}ani{\'c} and E.~H. Kim.
\newblock Mathematical analysis of the quasilinear effects in a hyperbolic
  model blood flow through compliant axi-symmetric vessels.
\newblock {\em Mathematical Methods in the Applied Sciences},
  26(14):1161--1186, 2003.

\bibitem{causemann2022human}
M.~Causemann, V.~Vinje, and M.~E. Rognes.
\newblock Human intracranial pulsatility during the cardiac cycle: a
  computational modelling framework.
\newblock {\em Fluids and Barriers of the CNS}, 19(1):1--17, 2022.

\bibitem{d2007multiscale}
C.~D'Angelo.
\newblock Multiscale modelling of metabolism and transport phenomena in living
  tissues.
\newblock Technical report, EPFL, 2007.

\bibitem{d2012finite}
C.~D'Angelo.
\newblock Finite element approximation of elliptic problems with {D}irac
  measure terms in weighted spaces: applications to one-and three-dimensional
  coupled problems.
\newblock {\em SIAM Journal on Numerical Analysis}, 50(1):194--215, 2012.

\bibitem{d2008coupling}
C.~D'angelo and A.~Quarteroni.
\newblock On the coupling of {1D} and {3D} diffusion-reaction equations:
  application to tissue perfusion problems.
\newblock {\em Mathematical Models and Methods in Applied Sciences},
  18(08):1481--1504, 2008.

\bibitem{daversin2022geometrically}
C.~Daversin-Catty, I.~G. Gjerde, and M.~E. Rognes.
\newblock Geometrically reduced modelling of pulsatile flow in perivascular
  networks.
\newblock {\em Frontiers in Physics}, page 360, 2022.

\bibitem{daversin2020mechanisms}
C.~Daversin-Catty, V.~Vinje, K.-A. Mardal, and M.~E. Rognes.
\newblock The mechanisms behind perivascular fluid flow.
\newblock {\em Plos one}, 15(12):e0244442, 2020.

\bibitem{delfour2011shapes}
M.~C. Delfour and J.-P. Zol{\'e}sio.
\newblock {\em Shapes and Geometries: Metrics,Analysis, Differential Calculus,
  and Optimization}.
\newblock SIAM, 2011.

\bibitem{di2012hitchhikers}
E.~Di~Nezza, G.~Palatucci, and E.~Valdinoci.
\newblock Hitchhiker's guide to the fractional {S}obolev spaces.
\newblock {\em Bulletin des Sciences Math{\'e}matiques}, 136(5):521--573, 2012.

\bibitem{evans2009partial}
L.~C. Evans.
\newblock {\em Partial Differential Equations}, volume~19.
\newblock American Mathematical Society, 2010.

\bibitem{fleischman1986interaction}
G.~J. Fleischman, T.~W. Secomb, and J.~F. Gross.
\newblock The interaction of extravascular pressure fields and fluid exchange
  in capillary networks.
\newblock {\em Mathematical Biosciences}, 82(2):141--151, 1986.

\bibitem{formaggia2001coupling}
L.~Formaggia, J.-F. Gerbeau, F.~Nobile, and A.~Quarteroni.
\newblock On the coupling of {3D} and {1D} {N}avier--{S}tokes equations for
  flow problems in compliant vessels.
\newblock {\em Computer Methods in Applied Mechanics and Engineering},
  191(6-7):561--582, 2001.

\bibitem{gjerde2020singularity}
I.~G. Gjerde, K.~Kumar, and J.~M. Nordbotten.
\newblock A singularity removal method for coupled {1D}--{3D} flow models.
\newblock {\em Computational Geosciences}, 24(2):443--457, 2020.

\bibitem{gjerde2019splitting}
I.~G. Gjerde, K.~Kumar, J.~M. Nordbotten, and B.~Wohlmuth.
\newblock Splitting method for elliptic equations with line sources.
\newblock {\em ESAIM: Mathematical Modelling and Numerical Analysis},
  53(5):1715--1739, 2019.

\bibitem{gong2014approximations}
W.~Gong, G.~Wang, and N.~Yan.
\newblock Approximations of elliptic optimal control problems with controls
  acting on a lower dimensional manifold.
\newblock {\em SIAM Journal on Control and Optimization}, 52(3):2008--2035,
  2014.

\bibitem{guermond2021finite}
J.-L. Guermond and A.~Ern.
\newblock {\em Finite Elements I: Approximation and Interpolation}.
\newblock Springer, 2021.

\bibitem{hannocks2018molecular}
M.-J. Hannocks, M.~E. Pizzo, J.~Huppert, T.~Deshpande, N.~J. Abbott, R.~G.
  Thorne, and L.~Sorokin.
\newblock Molecular characterization of perivascular drainage pathways in the
  murine brain.
\newblock {\em Journal of Cerebral Blood Flow \& Metabolism}, 38(4):669--686,
  2018.

\bibitem{slepc}
V.~Hernandez, J.~E. Roman, and V.~Vidal.
\newblock {SLEPc}: A scalable and flexible toolkit for the solution of
  eigenvalue problems.
\newblock {\em ACM Trans. Math. Softw.}, 31(3):351–362, sep 2005.

\bibitem{hladky2022glymphatic}
S.~B. Hladky and M.~A. Barrand.
\newblock The glymphatic hypothesis: the theory and the evidence.
\newblock {\em Fluids and Barriers of the CNS}, 19(1):1--33, 2022.

\bibitem{hofmann2007geometric}
S.~Hofmann, M.~Mitrea, and M.~Taylor.
\newblock Geometric and transformational properties of {L}ipschitz domains,
  {S}emmes-{K}enig-{T}oro domains, and other classes of finite perimeter
  domains.
\newblock {\em The Journal of Geometric Analysis}, 17(4):593--647, 2007.

\bibitem{kelley2022glymphatic}
D.~H. Kelley, T.~Bohr, P.~G. Hjorth, S.~C. Holst, S.~Hrab{\v{e}}tov{\'a},
  V.~Kiviniemi, T.~Lilius, I.~Lundgaard, K.-A. Mardal, E.~A. Martens, et~al.
\newblock The glymphatic system: Current understanding and modeling.
\newblock {\em Iscience}, page 104987, 2022.

\bibitem{koch2018new}
T.~Koch, K.~Heck, N.~Schr{\"o}der, H.~Class, and R.~Helmig.
\newblock A new simulation framework for soil--root interaction, evaporation,
  root growth, and solute transport.
\newblock {\em Vadose Zone Journal}, 17(1):1--21, 2018.

\bibitem{koch2020modeling}
T.~Koch, M.~Schneider, R.~Helmig, and P.~Jenny.
\newblock Modeling tissue perfusion in terms of 1d-3d embedded mixed-dimension
  coupled problems with distributed sources.
\newblock {\em Journal of Computational Physics}, 410:109370, 2020.

\bibitem{koch2022nonlinear}
T.~Koch, H.~Wu, and M.~Schneider.
\newblock Nonlinear mixed-dimension model for embedded tubular networks with
  application to root water uptake.
\newblock {\em Journal of Computational Physics}, 450:110823, 2022.

\bibitem{koppl20203d}
T.~K{\"o}ppl, E.~Vidotto, and B.~Wohlmuth.
\newblock A {3D}-{1D} coupled blood flow and oxygen transport model to generate
  microvascular networks.
\newblock {\em International Journal for Numerical Methods in Biomedical
  Engineering}, 36(10):e3386, 2020.

\bibitem{koppl2018mathematical}
T.~K{\"o}ppl, E.~Vidotto, B.~Wohlmuth, and P.~Zunino.
\newblock Mathematical modeling, analysis and numerical approximation of
  second-order elliptic problems with inclusions.
\newblock {\em Mathematical Models and Methods in Applied Sciences},
  28(05):953--978, 2018.

\bibitem{kuchta2020assembly}
M.~Kuchta.
\newblock Assembly of multiscale linear {PDE} operators.
\newblock In {\em Numerical Mathematics and Advanced Applications ENUMATH 2019:
  European Conference, Egmond aan Zee, The Netherlands, September 30-October
  4}, pages 641--650. Springer, 2020.

\bibitem{kuchta2021analysis}
M.~Kuchta, F.~Laurino, K.-A. Mardal, and P.~Zunino.
\newblock Analysis and approximation of mixed-dimensional {PDE}s on {3D}-{1D}
  domains coupled with {L}agrange multipliers.
\newblock {\em SIAM Journal on Numerical Analysis}, 59(1):558--582, 2021.

\bibitem{kuchta2019preconditioning}
M.~Kuchta, K.-A. Mardal, and M.~Mortensen.
\newblock Preconditioning trace coupled 3{D}-1{D} systems using fractional
  {L}aplacian.
\newblock {\em Numerical Methods for Partial Differential Equations},
  35(1):375--393, 2019.

\bibitem{kuttler1969inequality}
J.~Kuttler and V.~Sigillito.
\newblock An inequality for a {S}tekloff eigenvalue by the method of defect.
\newblock {\em Proceedings of the American Mathematical Society},
  20(2):357--360, 1969.

\bibitem{lamontagne2022recent}
E.~LaMontagne, A.~R. Muotri, and A.~J. Engler.
\newblock Recent advancements and future requirements in vascularization of
  cortical organoids.
\newblock {\em Frontiers in Bioengineering and Biotechnology}, page 2059, 2022.

\bibitem{laurino2019derivation}
F.~Laurino and P.~Zunino.
\newblock Derivation and analysis of coupled {PDE}s on manifolds with high
  dimensionality gap arising from topological model reduction.
\newblock {\em ESAIM: Mathematical Modelling and Numerical Analysis},
  53(6):2047--2080, 2019.

\bibitem{logg2012automated}
A.~Logg, K.-A. Mardal, and G.~Wells.
\newblock {\em Automated solution of differential equations by the finite
  element method: The FEniCS book}, volume~84.
\newblock Springer Science \& Business Media, 2012.

\bibitem{lohela2022glymphatic}
T.~J. Lohela, T.~O. Lilius, and M.~Nedergaard.
\newblock The glymphatic system: implications for drugs for central nervous
  system diseases.
\newblock {\em Nature Reviews Drug Discovery}, 21(10):763--779, 2022.

\bibitem{malenica2018groundwater}
L.~Malenica, H.~Gotovac, G.~Kamber, S.~Simunovic, S.~Allu, and V.~Divic.
\newblock Groundwater flow modeling in karst aquifers: Coupling 3{D} matrix and
  1{D} conduit flow via control volume isogeometric analysis—experimental
  verification with a 3{D} physical model.
\newblock {\em Water}, 10(12):1787, 2018.

\bibitem{mestre2018flow}
H.~Mestre, J.~Tithof, T.~Du, W.~Song, W.~Peng, A.~M. Sweeney, G.~Olveda, J.~H.
  Thomas, M.~Nedergaard, and D.~H. Kelley.
\newblock Flow of cerebrospinal fluid is driven by arterial pulsations and is
  reduced in hypertension.
\newblock {\em Nature Communications}, 9(1):1--9, 2018.

\bibitem{nance2022drug}
E.~Nance, S.~H. Pun, R.~Saigal, and D.~L. Sellers.
\newblock Drug delivery to the central nervous system.
\newblock {\em Nature Reviews Materials}, 7(4):314--331, 2022.

\bibitem{nicholson2001diffusion}
C.~Nicholson.
\newblock Diffusion and related transport mechanisms in brain tissue.
\newblock {\em Reports on progress in Physics}, 64(7):815, 2001.

\bibitem{nobile2001numerical}
F.~Nobile.
\newblock Numerical approximation of fluid-structure interaction problems with
  application to haemodynamics.
\newblock Technical report, EPFL, 2001.

\bibitem{nordbotten2009model}
J.~M. Nordbotten, D.~Kavetski, M.~A. Celia, and S.~Bachu.
\newblock Model for {CO2} leakage including multiple geological layers and
  multiple leaky wells.
\newblock {\em Environmental Science \& Technology}, 43(3):743--749, 2009.

\bibitem{notaro2016mixed}
D.~Notaro, L.~Cattaneo, L.~Formaggia, A.~Scotti, and P.~Zunino.
\newblock A mixed finite element method for modeling the fluid exchange between
  microcirculation and tissue interstitium.
\newblock {\em Advances in Discretization methods: Discontinuities, Virtual
  Elements, Fictitious Domain Methods}, pages 3--25, 2016.

\bibitem{possenti2019numerical}
L.~Possenti, G.~Casagrande, S.~Di~Gregorio, P.~Zunino, and M.~L. Costantino.
\newblock Numerical simulations of the microvascular fluid balance with a
  non-linear model of the lymphatic system.
\newblock {\em Microvascular Research}, 122:101--110, 2019.

\bibitem{possenti2021mesoscale}
L.~Possenti, A.~Cicchetti, R.~Rosati, D.~Cerroni, M.~L. Costantino, T.~Rancati,
  and P.~Zunino.
\newblock A mesoscale computational model for microvascular oxygen transfer.
\newblock {\em Annals of Biomedical Engineering}, 49:3356--3373, 2021.

\bibitem{rohan2018modeling}
E.~Rohan, V.~Luke{\v{s}}, and A.~Jon{\'a}{\v{s}}ov{\'a}.
\newblock Modeling of the contrast-enhanced perfusion test in liver based on
  the multi-compartment flow in porous media.
\newblock {\em Journal of Mathematical Biology}, 77(2):421--454, 2018.

\bibitem{sauter1999extension}
S.~Sauter and R.~Warnke.
\newblock Extension operators and approximation on domains containing small
  geometric details.
\newblock {\em East West Journal of Numerical Mathematics}, 7:61--77, 1999.

\bibitem{sloots2020cardiac}
J.~J. Sloots, G.~J. Biessels, and J.~J. Zwanenburg.
\newblock Cardiac and respiration-induced brain deformations in humans
  quantified with high-field {MRI}.
\newblock {\em Neuroimage}, 210:116581, 2020.

\bibitem{Stekloff1902}
W.~Stekloff.
\newblock Sur les probl{\'e}mes fondamentaux de la physique math{\'e}matique.
\newblock {\em Annales Scientifiques de l'École Normale Supérieure},
  19:191--259, 1902.

\bibitem{tarasoff2015clearance}
J.~M. Tarasoff-Conway, R.~O. Carare, R.~S. Osorio, L.~Glodzik, T.~Butler,
  E.~Fieremans, L.~Axel, H.~Rusinek, C.~Nicholson, B.~V. Zlokovic, et~al.
\newblock Clearance systems in the brain—implications for {Alzheimer}
  disease.
\newblock {\em Nature Reviews Neurology}, 11(8):457--470, 2015.

\bibitem{thomee2007galerkin}
V.~Thom{\'e}e.
\newblock {\em Galerkin Finite Element Methods for Parabolic Problems},
  volume~25.
\newblock Springer Science \& Business Media, 2007.

\bibitem{vinje2021brain}
V.~Vinje, E.~N. Bakker, and M.~E. Rognes.
\newblock Brain solute transport is more rapid in periarterial than perivenous
  spaces.
\newblock {\em Scientific Reports}, 11(1):1--11, 2021.

\bibitem{vinje2023human}
V.~Vinje, B.~Zapf, G.~Ringstad, P.~K. Eide, M.~E. Rognes, and K.~Mardal.
\newblock Human brain solute transport quantified by glymphatic {MRI}-informed
  biophysics during sleep and sleep deprivation.
\newblock {\em bioRxiv}, pages 2023--01, 2023.

\bibitem{wardlaw2020perivascular}
J.~M. Wardlaw, H.~Benveniste, M.~Nedergaard, B.~V. Zlokovic, H.~Mestre, H.~Lee,
  F.~N. Doubal, R.~Brown, J.~Ramirez, B.~J. MacIntosh, et~al.
\newblock Perivascular spaces in the brain: anatomy, physiology and pathology.
\newblock {\em Nature Reviews Neurology}, 16(3):137--153, 2020.

\bibitem{wheeler2021bioengineering}
M.~L. Wheeler and M.~L. Oyen.
\newblock Bioengineering approaches for placental research.
\newblock {\em Annals of Biomedical Engineering}, pages 1--14, 2021.

\bibitem{zhao2022physiology}
L.~Zhao, A.~Tannenbaum, E.~N. Bakker, and H.~Benveniste.
\newblock Physiology of glymphatic solute transport and waste clearance from
  the brain.
\newblock {\em Physiology}, 37(6):349--362, 2022.

\end{thebibliography}
\bibliographystyle{abbrv}

\appendix
\section{Technical estimates}
 \subsection{Uniform bound on $C_1 , C_2 $ and $N_2$ as defined in Propositions \ref{prop:first_model_err} and \ref{lemma:model_error_3d}} \label{sec:bound_W} We derive a bound  independent of $\epsilon_{\max}$ on  $C_1 + C_2$ under the assumption that $\partial_t c \in L^2(0,T; L^2(\Omega))$ and $\partial_t \cc \in L^2(0,T;L^2_A(\Lambda))$ which follow from maximal regularity, see Proposition \ref{prop:well_posedness_regularity}. From the definitions of $C_1$ and $C_2$, it suffices to obtain a bound on:  
 \begin{align}
 W_1    & =  \|\cc\|_{L^2(0,T;H^1(\Omega_v(t)))}+ \|c\|_{L^2(0,T;H^1(\Omega))} + \|\cc\|_{L^2(0,T;L^2(\partial \Omega_v(t)))} + \|\overline{c}\|_{L^2(0,T;L^2(\Gamma(t)))}, \\   W_2  &=\|c_s\|_{L^2(0,T;H^1(\Omega_s(t)))} +  \|\overline{c_s}\|_{L^2(0,T;L^2(\Gamma(t)))}. 
 \end{align}
 We will make use of the following estimate. With Cauchy-Schwarz inequality, we have that
 \begin{equation}
    b_\Lambda(v,w) \leq   \|\xi^{1/2}v\|_{L^2_{P}(\Lambda)} \| \xi^{1/2}w \|_{L^2_{P}(\Lambda)}, \quad \forall v,w \in L^2_{P} (\Lambda). \label{eq:bounded_b_lambda}
 \end{equation}
 In \eqref{eq:coupled_pb_1}, choose $v = c$. We obtain 
 \begin{align}
 \frac{1}{2} \int_{\Omega} \partial_t c^2 + \tilde{\nu}\|\nabla c\|^2_{L^2(\Omega)}  + \|\xi^{1/2}\bar{c}\|_{L_{P}^2(\Lambda)}^2  = (\mathcal{E}\bm{u}_s c, \nabla c)_{\Omega} +  b_{\Lambda}(\cc,\bar{c}) +(\mathcal{E}f_s, c )_{\Omega}. 
 \end{align}
 With H\"{o}lder's inequality and \eqref{eq:bounded_b_lambda}, we obtain 
 \begin{align} \nonumber
     \frac{1}{2} \int_{\Omega} \partial_t c^2 +  \tilde{\nu} \|\nabla c\|^2_{L^2(\Omega)}  + \|\xi^{1/2}\bar{c}\|_{L_{P}^2(\Lambda)}^2  \leq \|\mathcal{E}\bm{u}_s\|_{L^{\infty}(\Omega)}\|c\|_{L^2(\Omega)}\|\nabla c\|_{L^2(\Omega)} \\ +\nonumber   \|\xi^{1/2}\cc\|_{L^2_{P}(\Lambda)} \|\xi^{1/2}\bar{c}\|_{L^2_{P}(\Lambda)} + \|\mathcal{E}f_s\|_{L^2(\Omega)} \|c\|_{L^2(\Omega)}.  
 \end{align}
 Applying Young's inequality and multiplying by $2$ yields:
 \begin{multline}
      \int_{\Omega} \partial_t c^2 + \tilde{\nu}\|\nabla c\|^2_{L^2(\Omega)}  +  \|\xi^{1/2}\bar{c}\|_{L_{P}^2(\Lambda)}^2  \\ \leq (1+  \tilde{\nu}^{-1} \|\mathcal{E}\bm{u}_s\|^2_{L^{\infty}(\Omega)})\|c\|^2_{L^2(\Omega)} +  \|\xi^{1/2} \cc\|^2_{L^2_{P}(\Lambda)} +  \|\mathcal{E}f_s\|^2_{L^2(\Omega)}. \label{eq:justify_bdW1_0}
 \end{multline}
 Next, choose $v = \cc$ in \eqref{eq:coupled_pb_2}. Since $\avg{\cc} = \bar{\cc} = \cc$, we have 
 \begin{multline}
 \int_{\Lambda} \partial_t (A \cc) \cc  +  \|D_v^{1/2}\partial_s \cc\|_{L^2(\Omega_v)}^2 + \|\xi^{1/2}\cc\|^2_{L^2_{P}(\Lambda)}   = (\cc U, \partial_s \cc)_{\Omega_v} + b_{\Lambda} (\cc, \bar{c}) + (\avg{f_v} , \cc)_{\Omega_v}. 
 \end{multline}
In the above, it is implicitly understood that $\cc = E \cc$. Similarly, with H\"{o}lder's and Young's inequalities, we have 
\begin{multline}
2\int_{\Lambda} \partial_t (A\cc) \cc +  \|D_v^{1/2}\partial_x \cc\|_{L^2(\Omega_v)}^2 +  \|\xi^{1/2}\cc\|^2_{L^2_{P}(\Lambda)} \\ \leq (1+  \|D_v^{-1/2}U\|^2_{L^\infty(\Omega_v)})\|\cc\|^2_{L^2(\Omega_v)} + \|\xi^{1/2}\bar{c}\|_{L_{P}^2(\Lambda)}^2  + \|\avg{f_v}\|_{L^2(\Omega_v)}^2 .  \label{eq:justify_bdW1_1}   
\end{multline}
Observe that with Reynolds transport theorem, we have that 
\begin{align} \label{eq:integ_leibinz}
& \int_{\Lambda} \partial_t(A\cc)\cc =\int_{\Lambda } \partial_t(A\cc^2) - \frac12 \int_{\Lambda} A \partial_t \cc^2 = \partial_t \int_{\Lambda } (A\cc^2) - \frac12 \int_{\Omega_v(t)}\partial_t \cc^2 \\  &= \frac12 \partial_t \int_{\Omega_v(t)} \cc^2 + \frac12 \int_\Lambda \cc^2  \int_{\partial\Theta} \bm{w} \cdot \bm{n}  = \frac12 \partial_t \int_{\Omega_v(t)} \cc^2 + \frac12 \int_\Lambda \cc^2  \int_{\Theta} \nabla \cdot \bm{w}   \nonumber \\ &  \nonumber = \frac12 \partial_t \int_{\Omega_v(t)} \cc^2 + \frac12 \int_{\Omega_v(t)} \cc^2  \nabla \cdot \bm{w}. \nonumber
 \end{align}
Adding \eqref{eq:justify_bdW1_0} and \eqref{eq:justify_bdW1_1} with using \eqref{eq:integ_leibinz},  and integrating from  $0,t$, we obtain 
\begin{align*}
  &\|c(t)\|^2_{L^2(\Omega)} + \|\cc(t)\|_{L^2(\Omega_v)}^2 +  \tilde{\nu}\|\nabla c\|_{L^2(0,t;L^2(\Omega))}^2 +  \|D_v^{1/2}\partial_s \cc\|_{L^2(0,t;L^2(\Omega_v))}^2 \\ & \leq   \|\mathcal{E}f_s\|_{L^2(0,t;L^2(\Omega))}^2 + \|\avg{f_v}\|_{L^2(0,t;L^2(\Omega_v))}^2 +\|c^0\|^2_{L^2(\Omega)} + \|\cc^0\|^2_{L^2(\Omega_v)}  \\ & +
  (1+  \tilde{\nu}^{-1}\|\mathcal{E}\bm{u}_s\|^2_{L^{\infty}(0,T;L^{\infty}(\Omega))})\|c\|^2_{L^2(0,t;L^2(\Omega_v))} \\ & +   (1+ \|D_v^{-1/2} U\|^2_{L^\infty(0,T;L^{\infty}(\Omega_v))} + \|\nabla \cdot \bm{w}\|_{L^{\infty}(0,T;L^\infty(\Omega))} )\|\cc\|^2_{L^2(0,t;L^2(\Omega_v))} .
  \end{align*}
With continuous Gronwall's inequality, we obtain 
\begin{equation}
\|c\|_{L^\infty(0,T;L^2(\Omega))}^2 + \|\cc\|_{L^\infty(0,T;L^2(\Omega_v))}^2 +  \tilde{\nu} \| \nabla  c\|_{L^2(0,T;L^2(\Omega))}^2 +   \|D_v^{1/2} \partial_s \cc\|_{L^2(0,T;L^2(\Omega_v))}^2 \leq K. 
\end{equation}
Here, $K$ depends on the final time $T$, the velocities $\bm{u}_s$ and $U$, the initial conditions, and the source terms, but not on $|\Omega_v|$. Further, from  \eqref{eq:trace_eq_Omegas} and with the above bound, we have 
\begin{equation}
\|\overline{c}\|_{L^2(0,T;L^2_{P}(\Lambda))} \leq \|c\|_{L^2(0,T;L^2(\Gamma))} \leq K_\Gamma \epsilon^{1/6} \|c\|^2_{L^2(0,T;H^1(\Omega
))} \leq  K.  \label{eq:regularity_bounds_interface_0}
 \end{equation} 
We then use the above in \eqref{eq:justify_bdW1_1} after integrating over time. With Gronwall's inequality, we can conclude the bound 
\begin{equation}
 \|\cc\|_{L^2(0,T;L^2(\Gamma))}^2 \leq K.  \label{eq:regularity_bounds_interface_1}
\end{equation}
To obtain a bound on $\|\cc \|_{L^2(0,T;L^2(\partial \Omega_v))}$, we first note that since $\cc$ is uniform on each cross-section:  
\begin{multline}
\|\cc\|^2_{L^2(\partial \Omega_v)} = \|\cc\|_{L^2(\Gamma)}^2 + \int_{\Lambda} \int_{\partial \Theta_1} \cc^2 = \|\cc\|_{L^2(\Gamma)}^2 + \int_{\Lambda}\cc^2 |\partial \Theta_1|   \\ 
= \|\cc\|_{L^2(\Gamma)}^2 + \int_{\Lambda} \frac{|\partial \Theta_1|}{|\partial \Theta_2|}  \int_{\Theta_2}\cc^2 \lesssim  \|\cc\|_{L^2(\Gamma)}^2 \lesssim K . 
\end{multline}
The above follows by the assumption on the radii $R_1$ and $R_2$ \eqref{eq:assump_radii}. 
Thus, $W_1$ is bounded independent of $\epsilon_{\max}$. 
For the original 3D-3D problem, choose $\bm{\phi} = (c_v, c_s)$ in \eqref{eq:weak_form_3d3d} and use \eqref{eq:coercivity_3d3d}. We obtain that 
\begin{multline} \label{eq:uniform_bd_1}
\langle \dot{c_v}, c_v\rangle_{H^{-1}(\Omega_v(t))} + (\nabla\cdot \bm{w} c_v,c_v)_{\Omega_v(t)} + \langle \dot{c_s}, c_s\rangle_{H^{-1}(\Omega_s(t))}  +(\nabla\cdot \bm{w} c_s,c_s)_{\Omega_s(t)} \\ + K_1 (\|c_v\|_{H^1(\Omega_v(t))}^2 + \|c_s\|_{H^1(\Omega_s(t))}^2)  
\leq (K_2+1) (\|c_v\|_{L^2(\Omega_v(t))}^2 +  \|c_s\|_{L^2(\Omega_s(t))}^2 )   + \|f_v\|^2_{L^2(\Omega_v)} + \|f_s\|^2_{L^2(\Omega_s)}. 
\end{multline}
Note that the constants $K_1$ and $K_2$ are independent of $\epsilon$. 
With \cite[Corollary 2.41]{alphonse2015abstract} for $ i \in \{v,s\}$:
\begin{align}
\int_0^t\langle \dot{c_i} , c_i\rangle_{H^{-1}(\Omega_i(t))} = \frac12 (\|c_i(t)\|_{L^2(\Omega_i(t))}^2 - \|c_i(0)\|_{L^2(\Omega_i(0))}^2) - \int_0^t \frac12 (\nabla \cdot \bm{w} c_i, c_i)_{\Omega_i(t)}, 
\end{align} we obtain: 
\begin{align}
\frac{1}{2} (& \|c_v(t)\|^2_{L^2(\Omega_v(t))} + \|c_s(t)\|^2_{L^2(\Omega_s(t))} ) + K_1 (\|c_v\|_{L^2(0,t;H^1(\Omega_v(t)))}^2 + \|c_s\|_{L^2(0,t;H^1(\Omega_s(t)))}^2)  
\\ \nonumber &\leq (K_2+1 + \|\nabla \cdot \bm{w}\|_{L^\infty(0,t;L^\infty(\Omega))}) (\|c_v\|_{L^2(0,t;L^2(\Omega_v(t)))}^2 +  \|c_s\|_{L^2(0,t;L^2(\Omega_s(t)))}^2 )  \\ \nonumber & \,\,\,\,\, +  \frac12 (\|c_v(0)\|^2_{L^2(\Omega_v(0))} + \|c_s(0)\|^2_{L^2(\Omega_s(0))}) + \|f_v\|^2_{L^2(0,t;L^2(\Omega_v))} + \|f_s\|^2_{L^2(0,t;L^2((\Omega_s)))}.
\end{align}
With Gronwall's inequality, we can then conclude that: 
\begin{multline}
   \|c_v\|_{L^\infty(0,T;L^2(\Omega_v(t)))}^2 +  \|c_s\|_{L^\infty(0,T;L^2(\Omega_s(t)))}^2 \\ + \| c_v\|_{L^2(0,T;H^1(\Omega_v(t)))}^2 + \|c_s\|_{L^2(0,T;H^1(\Omega_s(t)))}^2  
    \leq K.\label{eq:original_pb_1}
    \end{multline}
Similar to the arguments above, we use trace estimate \eqref{eq:trace_eq_Omegas} to  conclude that 
\begin{equation}
\|\overline{c_s}\|_{L^2(0,T;L^2(\Gamma))} \leq \|c_s\|_{L^2(0,T;L^2(\Gamma))}  \leq K.  \label{eq:interface_bound_original3d3d0}
\end{equation}
Thus, $W_2$ is bounded which concludes the argument that $W$ is bounded independent of  $\epsilon_{\max}$. 
Then, we use the above bound in the energy inequality for $c_v$ (which can be obtained by choosing $\bm{\phi}=(c_v, 0)$ in \eqref{eq:weak_form_3d3d}) to obtain that 
\begin{equation}
    \|c_v\|_{L^2(0,T;L^2(\Gamma))}  \leq K.  \label{eq:interface_bound_original3d3d}
    \end{equation}
  This shows that $N_2$, as defined in Proposition~\ref{lemma:model_error_3d} is bounded independent of $\epsilon_{\max}$. 
\subsection{H\"{o}lder continuity bound} \label{appendix:handling_b} We verify that the terms involving $b_{\Lambda}$ in $\mathcal{A}_{\Lambda}(t,\cdot, \cdot)$ satisfy the bound stated in \eqref{eq:Holder_cont_bound}. We focus on the following term as it is the most intricate: 
\begin{align}\label{eq:dealing_with_b}
  b_{\Lambda}( \overline{c}(t_1)&, \overline{v}(t_1))  - b_{\Lambda}(\overline{c}(t_2), \overline{v}(t_2)) = \int_{\Gamma(t_2)} \xi(t_2)c  \overline{v}(t_2) - \int_{\Gamma(t_1)} \xi(t_1)c  \overline{v}(t_1) \\ & =  \int_{\Gamma(t_2)} (\xi(t_2) - \xi(t_1)) c  \overline{v}(t_2)  + \left( \int_{\Gamma(t_2)} \xi(t_1)c  \overline{v}(t_2)  - \int_{\Gamma(t_1)} \xi(t_1)c  \overline{v}(t_2) \right) \nonumber \\ & + \int_{\Gamma(t_1)} \xi(t_1)c \left(  \overline{v}(t_2)  - \overline{v}(t_1) \right) = T_1 + T_2 +T_3. \nonumber
\end{align} 
For the first term, the H\"{o}lder continuity bound follows form the assumption on $\xi$, the observation that $\|\overline{v}\|_{L^2(\Gamma(t_2))} \leq \|v\|_{L^2(\Gamma(t_2))}$,  and from \eqref{eq:trace_constant_indep_t}. Indeed, we have that 
\begin{align}
  \int_{\Gamma(t_2)} (\xi(t_2) - \xi(t_1)) c  \overline{v}(t_2) & \leq \|\xi(t_2) - \xi(t_2)\|_{L^{\infty}(\Lambda)} \|c\|_{L^2(\Gamma(t_2))} \|\overline{v}\|_{L^2(\Gamma(t_2))} \\ & \leq C |t_2 - t_1|^{\beta} \|c\|_{H^{1}(\Omega)} \|v\|_{H^1(\Omega)}. \nonumber
\end{align}
To handle the second term in \eqref{eq:dealing_with_b}, observe that the specific domains of section \ref{sec:geometry_specific} allow us to write 
\begin{multline*}
 T_2 = \int_{\Lambda} \int_{0}^{2\pi} (R_2(t_2) - R_2(t_1)) \xi(t_1) c \overline{v}(t_2) \leq C |t_2 - t_1|^{\beta} \|c\|_{L^{2}(\Gamma(0))} \|\overline{v}(t_2)\|_{L^2(\Gamma(0))} \\ \leq C |t_2 - t_1|^{\beta} \|c\|_{H^{1}(\Omega)} \|v\|_{H^1(\Omega)}.
\end{multline*}
In the above, we assumed without loss of generality that  $\partial \Theta_2(0)$  has radius $1$, and we use that $R_2$ is H\"older continuous. 
For $T_3$, we have for  a.e in $\Lambda$ 
\begin{align*}
& |\overline{v}(t_2) - \overline{v}(t_1) |  =  \left(\frac{1}{P(t_2) }  - \frac{1}{P(t_1)} \right) \int_0^{2\pi} R_2(t_2) v  + \frac{1}{P(t_1)} \int_0^{2\pi} (R_2(t_2) - R_2(t_1)) v  \\ &
\leq \left \|\frac{1}{P(t_2)} -\frac{1}{P(t_1) }\right \|_{L^{\infty}(\Lambda)} \|R_2(t_2)\|_{L^2(\partial \Theta_2(0))} \|v\|_{L^2(\partial \Theta_2(0))} \\ &  \quad + \left \|\frac{1}{P(t_1)}\right \|_{L^{\infty}(\Lambda)} \|R_2(t_2) - R_2(t_1)\|_{L^2(\partial  \Theta_2(0))} \|v\|_{L^2(\partial \Theta_2(0))}.
\end{align*}
Hence, $T_3$ is bounded as follows: 
\begin{equation*}
T_3 \leq  \|c\|_{L^2(\Gamma(t_1))} \|\overline{v}(t_2) - \overline{v}(t_1)\|_{L^2(\Gamma(t_1))} \leq C |t_2 -t_1|^{\beta} \|c\|_{H^1(\Omega) } \|v\|_{H^1(\Omega)}, 
\end{equation*}
where we used the H\"{o}lder continuity assumptions of $P^{-1}$ and on $R_2$ with follow from the smoothness assumptions on the domain velocity $\bm{w} \in C^2$. The other terms with $b_{\Lambda}$ can be handled similarly. 

\end{document}